\documentclass[11pt,a4paper]{article}
\usepackage{amsmath}
\usepackage[colorlinks=true]{hyperref}
\def\RiskOpt{\mathrm{RiskOpt}}
\def\Cov{{\hbox{\rm Cov}}}
\def\cV{{\cal V}}
\def\cA{{\cal A}}
\def\cI{{\cal I}}
\def\rX{{ x}}
\def\rV{{v}}
\def\cH{{\cal H}}

\def\cP{{\cal P}}

\def\cS{{\cal S}}
\def\cM{{\cal M}}

\def\Col{{\hbox{\rm Col}}}
\def\cQ{{\cal Q}}
\def\cE{{\cal E}}

\def\cZ{{\cal Z}}
\def\cB{{\cal B}}
\def\cY{{\cal Y}}
\def\Risksigma{\hbox{\rm Risk}_{\cH}}

\def\Sh{\hbox{\tiny\rm Sh}}
\def\cW{{\cal W}}
\def\cT{{\cal T}}
\def\cN{{\cal N}}
\def\rank{{\hbox{\rm Rank}}}
\def\cX{{\cal X}}

\def\Riskopt{\Risk_{\hbox{\tiny\rm opt}}}

\def\cl{\mathop{\hbox{\rm cl}}}
\def\inter{\mathop{\hbox{\rm int}}}

\def\bT{{\mathbf{T}}}
\def\bK{{\mathbf{K}}}
\newcommand{\half}{ \mbox{\small$\frac{1}{2}$}}
\newcommand{\four}{\hbox{\small$\frac{1}{4}$}}
\usepackage{amsfonts}
\usepackage{epsfig}

\def\Risk{{\hbox{\rm Risk}}}

\def\Ker{\mathop{\hbox{\rm Ker}}}

\def\Z{{\cal Z}}

\usepackage{amssymb}

\def\cR{{\cal R}}
\def\Diag{{\hbox{\rm Diag}}}
\usepackage{graphicx}
\def\bR{{\mathbf{R}}}

\def\Opt{{\hbox{\rm Opt}}}

\def\T{{\cal T}}

\def\R{{\cal R}}

\def\M{{\cal M}}
\def\X{{\cal X}}

\def\N{{\cal N}}
\def\cC{{\cal C}}
\def\S{{\cal S}}

\def\bE{{\mathbf{E}}}
\def\bS{{\mathbf{S}}}
\def\Prob{\hbox{\rm Prob}}

\newtheorem{lemma}{Lemma}[section]
\newtheorem{corollary}{Corollary}[section]
\newtheorem{proposition}{Proposition}[section]
\newtheorem{remark}{Remark}[section]
\newtheorem{theorem}{Theorem}[section]

\def\qed{$\Box$}

\def\e{{\rm e}}
\def\Tr{{\hbox{\rm Tr}}}

\newcommand{\aic}[2]{{\color{blue}  #2}}

\newcommand{\be}{\begin{eqnarray}}
\newcommand{\ee}[1]{\label{#1}\end{eqnarray}}
\newcommand{\nn}{\nonumber \\}
\newcommand{\ese}{\end{eqnarray*}}
\newcommand{\bse}{\begin{eqnarray*}}
\newcommand{\rf}[1]{~(\ref{#1})}
\newcommand{\wh}[1]{{\widehat{#1}}}
\newcommand{\hide}[1]{{}}
\def\qed{\hfill$\Box$\par}

\def\RiskPinormopt{\hbox{\rm RiskOpt}_{\Pi,\|\cdot\|}}
\def\qed{\hfill$\Box$\par}
\renewcommand{\Pi}{{\cQ}}
\renewcommand{\lll}{{\preccurlyeq}}
\title{Near-Optimality of Linear Recovery
 from Indirect Observations}
\author{
Anatoli Juditsky
\thanks{LJK, Universit\'e Grenoble Alpes, 700 Avenue Centrale 38041 Domaine Universitaire
de Saint-Martin-d'H\`{e}res, France,
{\tt anatoli.juditsky@univ-grenoble-alpes.fr}}
\and Arkadi Nemirovski
\thanks{Georgia Institute
 of Technology, Atlanta, Georgia
30332, USA, {\tt nemirovs@isye.gatech.edu}\newline
The first author was supported by the LabEx PERSYVAL-Lab (ANR-11-LABX-0025) and PGMO grant 2016-2032H. Research of
the second author was supported by NSF grants   CCF-1523768 and CMMI-1262063.}}
\date{}
\begin{document}
\maketitle

\begin{abstract}
We consider the problem of recovering linear image $Bx$ of a signal $x$ known to belong to a given convex compact set $\X$ from indirect observation
$\omega=Ax+\xi$ of $x$ corrupted by random noise $\xi$ with finite covariance matrix. It is shown that under some assumptions on $\X$ (satisfied, e.g., when $\X$ is the intersection of $K$ concentric ellipsoids/elliptic cylinders, or the unit ball of the spectral norm in the space of matrices) and on the norm $\|\cdot\|$ used to measure the recovery error (satisfied, e.g., by $\|\cdot\|_p$-norms, $1\leq p\leq 2$, on $\bR^m$ and by the nuclear norm on the space of matrices), one can build, in a computationally efficient manner,
a ``seemingly good'' {\em linear in observations estimate.} Further, in the case of zero mean Gaussian observation noise and general mappings $A$ and $B$, this estimate is
near-optimal among all (linear and nonlinear) estimates in terms of the maximal over $x\in \X$ expected $\|\cdot\|$-loss. These results form an essential extension of {classical results \cite{Pinsker1980,donoho1990minimax}} and of the recent work \cite{AnnStat}, where the assumptions on $\X$ were more restrictive, and the norm $\|\cdot\|$ was assumed to be the Euclidean one.
\par
This arXiv paper slightly strengthens the journal publication Juditsky, A., Nemirovski, A. Near-Optimality of Linear Recovery from Indirect Observations, {\sl Mathematical Statistics and Learning} {\bf 1:2} (2018), 171-225.
\end{abstract}
\section{Introduction}

Broadly speaking, what follows contributes to a long line of research (see, e.g., \cite{Don95,donoho1990minimax,efromovich1996sharp,efromovich2008nonparametric,IbrHas1981,Tsybakov,wasserman2006all} and references therein) started by the pioneering works \cite{kuks1,kuks2} and \cite{Pinsker1980} and aimed at building efficiently and analysing performance of linear estimates of signals from noisy observations. Specifically, we consider the classical estimation problem as follows: given a ``sensing matrix'' $A\in\bR^{m\times n}$ and an indirect noisy observation
\be
\omega=Ax+\xi
\ee{eq1obs}
of unknown deterministic ``signal'' $x$ known to belong to a given ``signal set'' $\cX\subset \bR^n$, we are interested to recover the linear image $Bx$ of the signal, where $B\in\bR^{\nu\times n}$ is a given matrix. We assume that the observation noise $\xi$ is random with unknown (and perhaps depending on $x$) distribution belonging to some family $\cP$ of Borel probability distributions on $\bR^m$ associated with a given nonempty convex compact subset $\Pi$ of the set of positive {\sl definite} $m\times m$ matrices. In this context, ``associated'' means that the  non-centered  covariance matrix\footnote{For the sake of brevity and with some terminology abuse, in the sequel, we refer to $\bE_{\xi\sim P}\{\xi\xi^T\}$ as to covariance matrix of $\xi\sim P$. Note that within the proposed approach we do not need the observation noise to be centered, except for the case of repeated observations, where we explicitly  request for the expectation of the noise to vanish (cf. Section \ref{reccovmatr} and Remark \ref{rem:multi}).} $\Cov[P]:=\bE_{\xi\sim P}\{\xi\xi^T\}$ of a distribution $P\in\cP$ is $\succeq$-dominated by some matrix from $\Pi$:
\begin{equation}\label{domination}
P\in \cP\Rightarrow \exists Q\in\Pi: \Cov[P]\preceq Q. \footnotemark
\end{equation}
\footnotetext{{Here and below,} $U\succeq V$ ($U\succ V$) means that $U,V$ are symmetric matrices of the same size and $U-V$ is positive semidefinite (resp., positive definite); $V\preceq U$ ($V\prec U$) means exactly the same as $U\succeq V$, (resp., $U\succ V$).}
We  quantify {the risk of} a candidate estimate -- a Borel function $\widehat{x}(\cdot):\bR^m\to\bR^\nu$ -- by its worst-case, under the circumstances, expected $\|\cdot\|$-error defined as
$$
\Risk_{\Pi,\|\cdot\|}[\widehat{x}|\cX]=\sup_{x\in \cX,P\in\cP}\bE_{\xi\sim P}\{\|Bx-\widehat{x}(Ax+\xi)\|\};
$$
here $\|\cdot\|$ is a given norm on $\bR^\nu$.
\par
We assume that signal set $\cX$ is a special type symmetric w.r.t. the origin convex compact set (a {\sl spectratope} to be defined in Section \ref{spec:1}), and require from the norm $\|\cdot\|_*$ conjugate to $\|\cdot\|$ to have a spectratope as the unit ball.\footnote{Obviously, any result of this type should impose some restrictions on $\X$ -- it is well known that linear estimates are ``heavily sub-optimal'' on some simple signal domains \cite{Nem85,donn95,dj98} (e.g., $\|\cdot\|_1$-ball).} This allows, e.g., for $\cX$ to be the (bounded) intersection of finitely many centered at the origin ellipsoids/elliptic cylinders/$\|\cdot\|_p$-balls ($p\in[2,\infty]$), or the (bounded) solution set of a system of two-sided Linear Matrix Inequalities
$$
\{x\in\bR^n: -L_k\preceq R_k[x]\preceq L_k, k\leq K\}\, \eqno{[R_k[x]: \hbox{\ linear in $x$ symmetric matrices}]}
$$
As for the norm $\|\cdot\|$, it can be $\|\cdot\|_p$-norm on $\bR^\nu$, $1\leq p\leq 2$, or the nuclear norm on the space $\bR^\nu=\bR^{u\times v}$ of matrices.
\par
 An important property of spectratopes is that they allow for precise concentration inequalities for random (Rademacher and Gaussian) vectors, see \cite{LP78,P78,Buch,Tropp111,Tropp112} and references therein. It plays a crucial role in what follows due to several important implications.
\begin{itemize}
\item
It allows for a tight computationally efficient upper bounding of the maximum of a quadratic
form over a spectratope (Proposition \ref{propmaxqf}). The latter allow to efficiently upper-bound the maximal over a spectratope risk of linear estimation (i.e., estimate of the form $\widehat{x}_H(\omega)=H^T\omega$), and thus leads to a computationally efficient scheme for building ``presumably good'' {\sl linear} estimates with guaranteed risk (Proposition \ref{summaryprop}).
\item
It is also decisive in demonstration of near-optimality of the resulting estimates (cf. \cite{Pinsker1980}): when the family $\cP$ of distributions contains all normal distributions $\{\cN(0,Q):Q\in\Pi\}$, it allows to tightly lower-bound the minimax risk of estimation over spectratopes via the Bayesian risk of estimating a random Gaussian signal, and to show that the presumably good  linear estimates are ``near-optimal'' (optimal up to logarithmic factors) among {\sl all} estimates, linear and nonlinear alike.
\end{itemize}
An ``executive summary'' of our main result -- Proposition \ref{newoptimalityprop} -- is as follows:
\begin{quote}
{\sl Given a spectratope $\cX$ and assuming that the unit ball $\cB_*$ of the norm conjugate to $\|\cdot\|$ is a spectratope as well, the efficiently computable optimal solution $H_*$ to an explicitly posed convex optimization problem yields a near-optimal linear estimate $\widehat{x}_{H_*}(\omega)=H_*^T\omega$, specifically,
$$
\Risk_{\cQ,\|\cdot\|}[\widehat{x}_{H_*}|\cX]\leq C\ln\big({\cal L}\big)\;
\RiskPinormopt[\cX],  \eqno{(*)}
$$
where $C$ is an absolute constant, ${\cal L}$ is polynomial in the (naturally defined) sizes of the spectratopes $\cX$, $\cB_*$, and $\RiskPinormopt$ is the ``true'' minimax optimal risk associated with zero mean Gaussian observation noises with covariance matrices from $\cQ$:
$$
\RiskPinormopt[\cX]=\sup_{Q\in\cQ}\inf_{\widehat{x}}\max_{x\in\cX} \bE_{\xi\sim\cN(0,Q)}\{\|Bx-\widehat{x}(Ax+\xi)\|\},
$$
the infimum being taken over all estimates, linear and nonlinear alike. }
\end{quote}
It should be stressed that the ``nonoptimality factor'' in $(*)$ is logarithmic and is completely independent of $B$ and of the sensing matrix $A$ -- the entities ``primarily responsible'' for the minimax optimal risk.
\par
{The above result constitutes an important extension to} the approach developed in \cite{AnnStat}, the progress as compared to \cite{AnnStat} being as follows:
\begin{itemize}
\item \cite{AnnStat} dealt with the case of $\cP=\{\cN(0,Q)\}$, i.e., the observation noise was assumed to be zero mean Gaussian with known covariance matrix, while now we allow for $\cP$ to be a general family of probability distributions with covariance matrices $\succeq$-dominated by matrices from a given convex compact set $\Pi\subset\inter \bS^m_+$; \footnote{From now on, $\bS^k$ stands for the space of symmetric $k\times k$ matrices, and $\bS^k_+$ is the cone of positive semidefinite matrices from $\bS^k$.}
\item present results apply to an essentially wider family of signal sets: spectratopes as compared to {\em ellitopes} considered in \cite{AnnStat}; ellitopes are also spectratopes, see Section {\ref{spec:1}}, but not vice versa. For instance, the intersection of centered at the origin ellipsoids/elliptic cylinders/$\|\cdot\|_p$-balls, $p\in[2,\infty]$, is an ellitope (and thus a spectratope), whereas the (bounded) solution set of a finite system of two-sided LMI's is a spectratope, but not an ellitope;
\item The analysis in \cite{AnnStat} was limited to the case of $\|\cdot\|_2$-losses, while now we allow for a much wider family of norms quantifying the recovery error.
\par
Note that, in addition to observations with random noise, in what follows we also address observations with  ``uncertain-but-bounded'' and ``mixed'' (combined) noise. In the latter case $\xi$, instead of being random, is selected, perhaps in adversarial manner, from a given spectratope -- the situation which was not considered in \cite{AnnStat} at all.
\end{itemize}
These results can also be considered
as a new contribution to the line of research initiated by \cite{donoho1990minimax}, where it is proved (Proposition 4 -- Theorem 7) that {\sl if $\cX$ is convex, orthosymmetric and quadratically convex} (that is, $\cX=\{x\in\bR^n:\exists t\in\cT: x_i^2\leq t_i,i\leq n\}$ with convex compact $\cT\subset\bR^n_+$), {\sl observations are direct: $\omega=x+\xi$, $\xi\sim\cN(0,I_n)$, $Bx=x$,  and $\|\cdot\|=\|\cdot\|_2$, the risk of an efficiently computable linear estimate is within factor 1.25 of the minimax optimal risk.} The suboptimality guarantees provided by the latter result are {\em essentially better} than those of Proposition \ref{newoptimalityprop} in the current paper. However, it is also {\em essentially more restrictive} in its scope -- an orthosymmetric convex and quadratically convex set is a very special case of an ellitope, the observations should be direct, and  $\|\cdot\|$ should be $\|\cdot\|_2$.
\par
Note that linear estimators can be efficiently built and optimized for some signal domains which are not spectratopes,
e.g., when $\cX$ is given as a convex hull of a finite set, e.g., $\cX$ is $\|\cdot\|_1$-ball (in this case, the smallest risk linear estimate can be ``heavily nonoptimal'' among all estimates). In general, however, optimizing risk over just linear estimates in a computationally efficient fashion can be problematic.
Beyond the scope of spectratopes in the role of signal sets and unit balls of the norms conjugate to those in which the recovery error is measured, the only known to us general situation where ``presumably good'' linear estimation is computationally tractable and results in (nearly) minimax optimal estimates is that where the recovery error is measured in $\|\cdot\|_\infty$. In the latter case, the breakthrough papers \cite{ikh88,Don95} (see also \cite{jn09}) {imply} that whenever $\cX$ is a computationally tractable convex compact set and the observation noise is Gaussian, an efficiently computable linear estimate is $\|\cdot\|_\infty$-minimax optimal within the factor $O(1)\sqrt{\ln(\nu)}$.

\par
The main body of the paper is organized as follows. We start with describing the family of sets we work with -- the spectratopes (Section \ref{spec:1}), and derive the crucial for the rest of the paper result on tight upper-bounding the maximum of a quadratic form over a spectratope (Section \ref{upperboundquadform}).  Next we explain how to build in a computationally efficient fashion ``presumably good'' linear estimates in the case of  stochastic (Section \ref{sect1}) and uncertain-but-bounded (Section \ref{sec:u-b-b-n}) observation noise and establish near-optimality of these estimates.  All technical proofs are relegated to Section \ref{sect:Proofs}. Appendix \ref{app:Scalc} lists principal rules of calculus of spectratopes. Appendix \ref{Processing} contains implementation details for the illustrative example
(covariance matrix estimation) presented in Section \ref{reccovmatr}. Finally, Appendix \ref{conicd} contains an ``executive summary'' of conic duality, our principal working horse.

\section{Preliminaries}\label{sect01}
We start with describing the main geometric object we intend to work with -- a {\sl spectratope}.
\subsection{Spectratopes}\label{spec:1}
\paragraph{A {basic} spectratope} is a set $\X\subset\bR^n$ given by {\sl {basic} spectratopic representation} -- representation of the form
\begin{equation}\label{sspectratope}
\X=\left\{x\in \bR^n: \exists t\in\cT: R_k^2[x]\preceq t_k I_{d_k},1\leq k\leq K\right\}
\end{equation}
where
\begin{enumerate}
\item[($S_1$)] $R_k[x]=\sum_{i=1}^n x_iR^{ki}$ are symmetric $d_k\times d_k$ matrices linearly depending on $x\in\bR^n$ (i.e., ``matrix coefficients'' $R^{ki}$ belong to $\bS^{d_k}$)
\item[($S_2$)] $\cT\in\bR^K_+$ is a {\sl monotonic} set, meaning that  $\cT$ is a convex compact subset of $\bR^K_+$ which contains a positive vector
and is monotone:
$$
0\leq t'\leq t\in\cT\Rightarrow t'\in\cT.{\,\footnotemark}
$$
\footnotetext{The inequalities between vectors are understood componentwise.}
\item[($S_3$)] Whenever $x\neq0$, it holds $R_k[x]\neq0$ for at least one $k\leq K$.
\end{enumerate}
\begin{remark}\label{remspect}
{\rm By Schur Complement Lemma, the set {\rm (\ref{sspectratope})} given by the data satisfying ($S_1$), ($S_2$) can be represented as
$$
\X=\left\{x\in \bR^n: \exists t\in\cT: \left[\begin{array}{c|c}t_kI_{d_k}&R_k[x]\cr\hline
R_k[x]&I_{d_k}\cr\end{array}\right]\succeq0,\,k\leq K\right\}
$$
By the latter representation, $\cX$  is nonempty, closed, convex, symmetric w.r.t. the origin, and contains a neighbourhood of the origin (the latter is due to the fact that $\cT$ contains a strictly positive vector). This set is bounded if and only if the data, in addition to ($S_1$), ($S_2$), satisfies ($S_3$).}
\end{remark}
\paragraph{A spectratope} $\X\subset \bR^p$ is a set represented as a linear image of a {basic} spectratope:
\begin{equation}\label{spectra}
\X=\{x\in\bR^p:\exists (y\in\bR^n,t\in\cT): x=Py,\;R_k^2[y]\preceq t_k I_{d_k},1\leq k\leq K\},
\end{equation}
where $P$ is a $p\times n$ matrix, and $R_k[\cdot]$, $\cT$ are as in ($S_1$)--($S_3$). We call the quantity
$$
D=\sum_{k=1}^Kd_k
$$
the {\sl size} of spectratope $\cX$.
\paragraph{Example 1: Ellitopes.} An {\sl ellitope} was defined in \cite{AnnStat} as a set $\X\subset\bR^n$ representable as
\begin{equation}\label{ellitope}
\X=\{x\in\bR^n:\exists(y\in\bR^N,t\in\cT):\;x=Py,\;y^TS_ky\leq t_k,\,k\leq K\},
\end{equation}
where $S_k\succeq0$, $\sum_KS_k\succ0$, and $\T$ satisfies ($S_2$). Basic examples of ellitopes are:
\begin{itemize}
\item bounded intersections of centered at the origin ellipsoids/elliptic cylinders: whenever $S_k
\succeq0$ and $\sum_kS_k\succ0$,
$$
\bigcap\limits_{k=1}^K\{x\in\bR^n:\;x^TS_kx\leq 1\}=\{x\in\bR^n:\;\exists t\in\cT=[0,1]^K: \;x^TS_kx\leq t_k,\,1\leq k\leq K\}.
$$
\item $\|\cdot\|_p$-balls, $2\leq p\leq\infty$:
$$
\{x\in\bR^n:\|x\|_p\leq1\}=\{x\in\bR^n: \exists t\in\T:=\{t\geq0,\|t\|_{p/2}\leq1\}: x^TS_kx:=x_k^2\leq t_k,\,k\leq n\}.
$$
\end{itemize}
It is immediately seen that an ellitope  (\ref{ellitope}) is a spectratope as well. Indeed, let  $S_k=\sum_{j=1}^{r_k} s_{kj}s_{kj}^T$, $r_k=\rank(S_k)$,  be a dyadic representation of the positive semidefinite matrix $S_k$, so that
$$
y^TS_ky=\sum_j(s_{kj}^Ty)^2\;\;\forall y,
$$
and let
$$
\widehat{\T}=\left\{\{t_{kj}\geq0,1\leq j\leq r_k,1\leq k\leq K\}:\exists t\in \T: \sum_jt_{kj}\leq t_k,\:k\leq K\right\}, \,\, R_{kj}[y]=s_{kj}^Ty\in\bS^1=\bR.
$$
We clearly have
$$
\X=\{x\in\bR^n:\exists (\{t_{kj}\}\in \widehat{\T},y): x=Py,\; R_{kj}^2[y]\preceq t_{kj}I_1\,\forall k,j\}
$$
and the right hand side is a {valid} spectratopic representation of $\X$. {Note that the spectratopic size of $\X$ is $D=\sum_{k=1}^K r_k$.}
\paragraph{Example 2: ``Matrix box.''} Let $L$ be a positive definite $d\times d$ matrix. Then the ``matrix box''
$$
\begin{array}{rcl}
\X&=&\{X\in \bS^d:-L\preceq X\preceq L\}=\{X\in\bS^d: -I_d\preceq L^{-1/2}XL^{-1/2}\preceq I_d\}\\
&=&\{X\in\bS^d: R^2[X]:=[L^{-1/2}XL^{-1/2}]^2\preceq I_d\}\\
\end{array}
$$
is a {basic} spectratope (augment  $R_1[\cdot]:=R[\cdot]$ with $K=1$, $\cT=[0,1]$). As a result, a {\sl bounded} set  $\X\subset\bR^n$ given by a system of ``two-sided'' Linear Matrix Inequalities, specifically,
$$
\cX=\{x\in\bR^n:\exists t\in\cT:\;-\sqrt{t_k}L_k\preceq S_k[x]\preceq \sqrt{t_k}L_k,\,1\leq k\leq K\}
$$
where $S_k[x]$ are symmetric $d_k\times d_k$ matrices linearly depending on $x$, $L_k\succ0$ and $\cT$ satisfies ($S_2$),
is a {basic} spectratope:
$$
\cX=\{x\in\bR^n:\exists t\in\cT:\;R_k^2[x]\preceq t_kI_{d_k},\,1\leq k\leq K\}\eqno{[R_k[x]=L_k^{-1/2}S_k[x]L_k^{-1/2}]}.
$$
More examples of spectratopes can be built using their ``calculus.'' It turns out that nearly all basic operations with sets preserving convexity, symmetry w.r.t. the origin, and boundedness (these are ``built-in'' properties of spectratopes), {such as} taking finite intersections, direct products, arithmetic sums, linear images, and inverse linear images under linear embeddings, as applied to spectratopes, yield spectratopes as well. Furthermore, a spectratopic representation of the result of such an operation is readily given by spectratopic representations of the operands; see Appendix \ref{app:Scalc} for principal calculus rules.
\par

\subsection{Upper-bounding quadratic form on a spectratope}\label{upperboundquadform}
We are about to establish the first crucial in our context fact about spectratopes --
the possibility to tightly upper-bound an (indefinite) quadratic form over a spectratope. To proceed, we need some definitions.
\paragraph{Linear maps associated with a spectratope.}
We associate with a basic spectratope (\ref{sspectratope}), ($S_1$)--($S_3$) the following entities:
\begin{enumerate}
\item[1.] Linear mappings
\[
Q\mapsto \cR_k[Q]=\sum_{i,j}Q_{ij}R^{ki}R^{kj}:\bS^n\to\bS^{d_k}
\]
\end{enumerate}
As is immediately seen, we have
\begin{equation}\label{khin201}
\cR_k[yy^T]\equiv R_k^2[y],
\end{equation}
implying that $\cR_k[Q]\succeq 0$ whenever $Q\succeq0$, whence $\cR_k[\cdot]$ is $\succeq$-monotone:
\begin{equation}\label{khin203}
Q'\succeq Q\Rightarrow \cR_k[Q']\succeq\cR_k[Q].
\end{equation}
Besides this, if $\xi$ is a random vector taking values in $\bR^n$ with covariance matrix $Q$, we have
\begin{equation}\label{khin202}
\bE_\xi\{R_k^2[\xi]\}=\bE_\xi\{\R_k[\xi\xi^T]\}=\R_k[\bE_\xi\{\xi\xi^T\}]={\R_k}[Q],
\end{equation}
where the first equality is given by (\ref{khin201}).
\begin{enumerate}
\item[2.] Linear mappings $\Lambda_k\mapsto\cR_k^*[\Lambda_k]:\bS^{d_k}\to\bS^n$ given by
\begin{equation}\label{cRkstar}
\left[\cR_k^*[\Lambda_k]\right]_{ij}=\half\Tr(\Lambda_k[R^{ki}R^{kj}+R^{kj}R^{ki}]),\,1\leq i,j\leq n.
\end{equation}
\end{enumerate}
It is immediately seen that $\cR_k^*[\cdot]$ is the adjoint of $\cR_k[\cdot]$:
\begin{equation}\label{khin23}
\forall (\Lambda_k\in\bS^{d_k},Q\in\bS^n): \langle \Lambda_k,\cR_k[Q]\rangle=\Tr(\Lambda_k\cR_k[Q])=\Tr(\cR_k^*[\Lambda_k]Q)=\langle\cR_k^*[\Lambda_k],Q\rangle,
\end{equation}
where $\langle A,B\rangle=\Tr(AB)$ is the Frobenius inner product of symmetric matrices. Besides this,  we have{\footnote{note that when $\Lambda_k\succeq0$ and $Q=yy^T$, the first quantity in (\ref{khin23}) is nonnegative by (\ref{khin201}), and therefore (\ref{khin23}) states that $y^T\cR^*_k[\Lambda_k]y\geq0$ for every $y$, implying that $\cR^*_k[\Lambda_k]\succeq0$.
}}
\begin{equation}\label{khin2201}
\Lambda_k\succeq0\Rightarrow \R^*_k[\Lambda_k]\succeq0;
\end{equation}
\begin{enumerate}
\item[3.] The linear space  $\Lambda^K=\bS^{d_1}\times...\times\bS^{d_K}$ of all ordered collections $\Lambda=\{\Lambda_k\in\bS^{d_k}\}_{k\leq K}$ along with the linear mapping
\[
\Lambda\mapsto \lambda[\Lambda]:=[\Tr(\Lambda_1);...;\Tr(\Lambda_K)]: \Lambda^K\to\bR^K.
\]
\end{enumerate}
Besides this, for a monotonic set $\cT\subset\bR^K$ we define
\begin{itemize}
\item the {\sl support function} of $\cT$
\[
\phi_{\cT}(g)=\max_{t\in \cT} g^Tt,
\]
which clearly is a convex positively homogeneous, of degree 1, nonnegative real-valued function on $\bR^K$. Since $\cT$ contains positive vectors, $\phi_{\cT}$ is coercive on $\bR^K_+$, meaning that $\phi_{\cT}(\lambda^s)\to+\infty$ along every sequence $\{\lambda^s\geq0\}$ such that $\|\lambda^s\|\to\infty$;
\item the conic hull of $\cT$
\[
\bK[\cT]=\cl\{[t;s]\in\bR^{K+1}: s>0,s^{-1}t\in\cT\}
\]
which clearly is a regular cone in $\bR^{K+1}$ (i.e., it is closed, convex, and pointed with a nonempty interior) such that
\[
\cT=\{t:[t;1]\in\bK[\cT]\}.
\]
It is immediately seen that the cone $(\bK[\cT])_*$ dual to $\bK[\cT]$ can be described as follows:
\[
(\bK[\cT])_*:=\{[g;r]\in\bR^{K+1}:[g;r]^T[t;s]\geq0\,\forall [t;s]\in\bK[\cT]\}=\{[g;r]\in\bR^{K+1}:r\geq \phi_{\cT}(-g)\}.
\]
\end{itemize}
\begin{proposition}\label{propmaxqf} Let $C$ be a symmetric $p\times p$ matrix, let ${\X}\subset\bR^p$ be given by spectratopic representation {\rm \rf{spectra}},
$$
\Opt=\max_{x\in \cX} x^TCx
$$
and let
\begin{equation}\label{optkhinpro}
\begin{array}{c}
\Opt_*=\min\limits_{\Lambda=\{\Lambda_k\}_{k\leq K}}\left\{\phi_\T(\lambda[\Lambda]):\Lambda_k\succeq0,k\leq K,P^TCP\preceq \sum_k\R_k^*[\Lambda_k]\right\}\\
\left[\lambda[\Lambda]=[\Tr(\Lambda_1);...;\Tr(\Lambda_K)]\right]\\
\end{array}
\end{equation}
Then {\rm (\ref{optkhinpro})} is solvable, and
\begin{equation}\label{boundskhin}
\Opt\leq\Opt_*\leq {2}\max[\ln({2}D),1]\Opt,
\end{equation}
where $D=\sum_kd_k$ is the size of the spectratope $\cX$.
\end{proposition}
To explain where the result of the proposition comes from, let us prove right now its easy part -- the first inequality in (\ref{boundskhin}); the remaining, essentially less trivial, part of the proof is provided in Section \ref{proofofpropmaxqf}.  Let $\Lambda$ be a feasible solution to the optimization problem in (\ref{optkhinpro}), and let $x\in \X$, so that $x=Py$ for some $y$ such that $R_k^2[y]\preceq t_kI_{d_k}$, $k\leq K$, for properly selected $t\in\cT$. We have
\bse
x^TCx&=&y^T[P^TCP]y\,{\underbrace{\leq}_{(a)}}
\,\sum_ky^T\cR^*_k[\Lambda_k]y=\sum_k\Tr(\cR^*_k[\Lambda_k]yy^T)\,
{\underbrace{=}_{(b)}}\,\sum_k\Tr(\Lambda_k\cR_k[yy^T])\\&
{\underbrace{=}_{(c)}}&\sum_k\Tr(\Lambda_{k}R^2_k[y])
{\underbrace{\leq}_{(d)}}\sum_k\Tr({\Lambda}_kt_kI_{d_k})=
\sum_kt_k\Tr(\Lambda_k)=\lambda^T[\Lambda]t{\underbrace{\leq}_{(e)}}\phi_{\cT}(\lambda[\Lambda]),
\ese
where $(a)$ is due to the fact that $\Lambda$ is feasible for the optimization problem in (\ref{optkhinpro}), $(b)$ is by (\ref{khin23}), $(c)$ is by (\ref{khin201}), $(d)$ is due to $\Lambda_K\succeq0$ and $R_k^2[y]{\preceq} t_kI_{d_k}$, and $(e)$ is by the definition of $\phi_\cT$. The bottom line is that the value of the objective of the {optimization} problem in  (\ref{optkhinpro}) at every feasible solution to this problem upper-bounds $\Opt$, implying the first inequality in (\ref{boundskhin}). Note that the derivation we have carried out is nothing but a minor modification of the standard semidefinite relaxation scheme.
\begin{remark}\label{remb}{\rm Proposition \ref{propmaxqf} has some history. When $\cX$ is an intersection of centered at the origin ellipsoids/elliptic cylinders, it was established in \cite{NRT}; matrix analogy of the latter result can be traced back to \cite{NOrtConstr}, see also \cite{ManChoSo}. The case when $\cX$ is  a general-type ellitope {\rm (\ref{ellitope})} was considered in \cite{AnnStat}, with tightness guarantee slightly better than in {\rm (\ref{boundskhin})}, namely,
$$
\Opt\leq\Opt_*\leq 4\ln(5K)\Opt.
$$
Note that in the case where $\X$ is an ellitope {\rm (\ref{ellitope})}, Proposition \ref{propmaxqf} results in a worse than $O(1)\ln(K)$ ``nonoptimality factor'' $O(1)\ln(\sum_{k=1}^K\rank(S_k))$. We remark that passing from ellitopes to spectratopes requires replacing elementary bounds on deviation probabilities used in \cite{NRT,AnnStat} by a more powerful tool  -- {matrix concentration inequalities}, see \cite{Tropp111,Tropp112} and references therein.}
\end{remark}

\section{Near-optimal linear estimation under random noise}\label{sect1}
\subsection{Situation and goal}\label{sitgoal}\label{sect3.1}
Given $\nu\times n$ matrix $B$, consider the problem of estimating linear image $Bx$ of unknown deterministic signal $x$ known to belong to a given set ${\X}
\subset\bR^n$ via noisy observation
\begin{equation}
\label{eq1obs3}
\omega=Ax+\xi
\end{equation}
where $A$ is a given $m\times n$ matrix $A$ and $\xi$ is random observation noise. In some signal processing applications, the distribution of noise is fixed and is part of the data of the estimation problem. In order to cover some interesting applications (cf. Section \ref{reccovmatr}), we allow for ``ambiguous'' noise distributions; all we know in advance is that this distribution belongs to a family $\cP$ of Borel probability distributions on $\bR^m$ associated, in the sense of (\ref{domination}),  with a given convex compact subset $\Pi$ of the interior of the cone $\bS^m_+$ of positive semidefinite $m\times m$  matrices. Actual distribution of noise in (\ref{eq1obs3}) is somehow selected from $\cP$ by nature (and may, e.g., depend on $x$).\par
In the sequel, for a Borel probability distribution $P$ on $\bR^m$ we write $P\lll \Pi$ to express the fact that $\Cov[P]$ is $\succeq$-dominated by a matrix from $\Pi$:
$$
\{P\lll \Pi\}\Leftrightarrow \{\exists {Q}\in\Pi: \Cov[P]\preceq{Q}\}.
$$
From now on we assume that {\sl all matrices from $\Pi$ are positive definite.}\par
Given $\Pi$ and a norm $\|\cdot\|$ on $\bR^\nu$,  we quantify the {\sl risk} of a candidate estimate -- a Borel function $\widehat{x}(\cdot):\bR^m\to\bR^\nu$ -- by its $(\Pi,\|\cdot\|)$-risk on $\X$ defined as
\begin{equation}\label{normrisk}
\Risk_{\Pi,\|\cdot\|}[\widehat{x}|\X] =\sup_{x\in\X,P\lll\Pi} \bE_{\xi\sim P}\left\{\|\widehat{x}(Ax+\xi)-Bx\|\right\},
\end{equation}
where $\|\cdot\|$ is some norm on $\bR^\nu$.
Our focus is on {\sl linear estimates} -- estimates of the form
$$
\widehat{x}_H(\omega)=H^T\omega
$$
given by $m\times\nu$ matrices $H$, and our current goal is to demonstrate that under some restrictions on the signal domain $\X$,  a ``good'' linear estimate yielded by an optimal solution to an efficiently solvable convex optimization problem is near-optimal in terms of its risk among {\sl all} estimates, linear and nonlinear alike.

We assume here that set $\X$ is a spectratope (cf. \rf{spectra}). Ideally, to compute a ``good'' linear estimate one would look for $H$ which minimizes the risk $\Risk_{\Pi,\|\cdot\|}[\widehat{x}_H|\X]$ in $H$.
This risk is generally difficult to compute even when $\X$ is a spectratope, and with our approach, we minimize in $H$ an efficiently computable upper bound on the risk rather than {the} risk itself. In order this bound to be tight -- good enough to allow to build a near-optimal linear estimate, we have to impose further restrictions, specifically,
we make from now on the following
 \begin{quote}
 {\bf Assumption $\mathbf{A}$:} {\em The unit ball $\cB_*$ of the norm $\|\cdot\|_*$ {\sl conjugate} to the norm $\|\cdot\|$ in the definition {\rm\rf{normrisk}} of the estimation risk {\em is a spectratope:}
\begin{equation}\label{ber6spec}
\begin{array}{rcl}
\cB_*&=&\{z\in\bR^\nu: \exists  y\in\cY:z=My\},\\
\cY~&:=&\{y\in\bR^q:\exists r\in\R: S_\ell^2[y]\preceq r_\ell I_{f_\ell},\,1\leq\ell\leq L\},
\end{array}
\end{equation}
where the right hand side data are as required in a spectratopic representation.}
\end{quote}
{Examples of norms $\|\cdot\|$ satisfying Assumption $\mathbf{A}$ include $\|\cdot\|_q$-norms on $\bR^\nu$, $1\leq q\leq 2$ ({conjugates of  the norms $\|\cdot\|_p$ with $1/p+1/q=1$,} see Example 1 in Section {\ref{spec:1}}). Another example is {\sl nuclear norm} $\|\cdot\|_{\Sh,1}$ on the space $\bR^\nu=\bR^{p\times q}$ of $p\times q$ matrices, $\|V\|_{\Sh,1}=\sum \sigma_i(V)$ -- the sum of singular values of a matrix $V$.  The  conjugate of the nuclear norm is the spectral norm $\|\cdot\|_{\Sh,\infty}$ on $\bR^\nu=\bR^{p\times q}$, and the unit ball of the latter norm is a basic spectratope (cf. Example 2 in Section {\ref{spec:1}}):
$$
\{Z\in\bR^{p\times q}:\,\|Z\|_{\Sh,\infty}\leq1\}=\{Z:\,\exists r\in\cR=[0,1]: S^2[Z]\preceq tI_{p+q}\},\,\,S[Z]=\left[\begin{array}{c|c}&Z^T\cr\hline Z&\cr\end{array}\right].
$$}
{It} is immediately seen that the case when $\cX$ is a spectratope (\ref{spectra}) can be reduced to the one where $\cX$ is a {\sl basic} spectratope -- to this end it suffices to replace matrices $A$ and $B$ with $AP$ and $BP$, respectively, and to treat $y$ rather than $x=Py$ {as the signal underlying observation (\ref{eq1obs3})}, see (\ref{spectra}). We assume that this reduction has been carried out in advance, so that from now on our signal set will be
\[
\X=\{x\in\bR^n:\,\exists t\in\cT: \,R_k^2[x]\preceq t_k I_{d_k},1\leq k\leq K\}.
\]

\subsection{Building linear estimate}\label{newbeyond}
Observe
 that
 the $(\Pi,\|\cdot\|)$-risk of the linear estimate $\widehat{x}_H(\omega)=H^T\omega$, $H\in\bR^{m\times \nu}$, can be upper-bounded as follows:
\be
\Risk_{\Pi,\|\cdot\|}[\widehat{x}_H(\cdot)|\cX]&=&\sup_{x\in \X,P\lll\Pi}\bE_{\xi\sim P}\{\|H^T(Ax+\xi)-Bx\|\}\nn
&\leq& \sup_{x\in \X}\|H^TAx-Bx\|+\sup_{P\lll\Pi}\bE_{\xi\sim P}\{\|H^T\xi\|\}\nn
&\leq& \Phi_{\cX}({H})+\Psi_\Pi(H),
\ee{ber7spec}
where
\bse
{\Phi_{\cX}({H})}&=&\max_x\left\{\|{(B-H^TA)}x\|:\,x\in \cX\right\},\quad
\Psi_\Pi(H)=\sup_{P\lll\Pi}\bE_{\xi\sim P}\left\{\|H^T\xi\|\right\}.\\
\ese
While {$\Phi_{\cX}({H})$} and $\Psi_\Pi(H)$ are convex functions of $H$, these functions can be difficult to compute.\footnote{For instance, computing $\Psi_{\cX}({H})$ reduces to maximizing the convex function $\|(B-H^TA)x\|$ over $x\in \cX$, which is computationally intractable even when $\X$ is as simple as the unit box, and $\|\cdot\|$ is the Euclidean norm.} {In such a case, matrix $H$ of a} ``good'' linear estimate $\widehat{x}_H$ {which is also efficiently computable can be chosen as} a minimizer of the sum of  efficiently computable convex upper bounds on
{$\Phi_{\cX}$} and  $\Psi_\Pi$. 

\subsubsection{Upper-bounding ${\Phi_{\cX}(\cdot)}$}
With Assumption $\mathbf{A}$ in force, let us consider the direct product spectratope
\[
\begin{array}{rcl}
\Z&:=&\cX\times \cY=\{[x;y]\in\bR^n\times\bR^q: \exists s=[t;r]\in\T\times\R:\\
&&\multicolumn{1}{r}{ R_k^2[x]\preceq t_kI_{d_k},\,1\leq k\leq K,S_\ell^2[y]\preceq r_\ell I_{f_\ell},\,1\leq \ell\leq L\}}\\
&=&\{w=[x;y]\in\bR^n\times\bR^q: \exists s=[t;r]\in\S=\T\times\R:U_i^2[w]\preceq s_i I_{g_i},\\
&&\multicolumn{1}{r}{1\leq i\leq I=K+L\}}\\
\end{array}
\]
with $U_{i}[\cdot]$ readily given by $R_k[\cdot]$ and $S_\ell[\cdot]$.  Given a $\nu\times n$ matrix $V$ and setting
$$
W[V]={1\over2}\left[\begin{array}{c|c}&V^TM\cr\hline M^TV&\cr\end{array}\right]
$$
it clearly holds
$$
\max_{x\in \cX}\|Vx\|=\max_{x\in \cX,z\in\cB_*}z^TVx=\max_{x\in\cX,y\in\cY}y^TM^TVx=\max_{w\in \Z}w^TW[V]w.
$$
Applying Proposition \ref{propmaxqf}, we arrive at the following
result (cf. Proposition \ref{propbnoiseopt}):
{  \begin{corollary}\label{corrbyproductspec} In the just defined situation, the efficiently computable convex function
\begin{equation}\label{ber11spec}
\begin{array}{rcl}
{\overline\Phi_{\cX}(H)}&=&\min\limits_{\Lambda,\Upsilon}\bigg\{\phi_\T(\lambda[\Lambda])+\phi_{\R}(\lambda[\Upsilon]):\Lambda=\{\Lambda_k\in\bS^{d_k}_+\}_{k\leq K},\Upsilon=\{\Upsilon_\ell\in\bS^{f_\ell}_+\}_{\ell\leq L},\\
&&\multicolumn{1}{r}{\left[\begin{array}{c|c}\sum_k\R_k^*[\Lambda_k]&{1\over2} {(B-H^TA)}^TM\cr\hline {1\over2} M^T{(B-H^TA)}&\sum_\ell\S_\ell^*[\Upsilon_\ell]\cr\end{array}\right]\succeq0\bigg\}}
\end{array}
\end{equation}
 is a tight upper bound on $\Phi_\cX(\cdot)$, namely,
\[
\begin{array}{c}
\forall H\in\bR^{m\times\nu}:\;{\Phi_{\cX}({H})\leq \overline\Phi_{\cX}({H})}\leq 2\max[\ln(2D),1]\,{\Phi_{\cX}({H})},\;\;
D=\sum_kd_k+\sum_\ell f_\ell.
\end{array}
\]
Recall, that here
\be\begin{array}{rcl}
{[\R_k^*[\Lambda_k]]_{ij}}&=&{1\over 2}\Tr(\Lambda_k[R_k^{ki}R_k^{kj}+R_k^{kj}R_k^{ki}]),
\hbox{\ where\ } R_k[x]=\sum_ix_iR^{ki},\\
{[\S_\ell^*[\Upsilon_\ell]]_{ij}}&=&{1\over 2}\Tr(\Upsilon_\ell[S_\ell^{\ell i}S_\ell^{\ell j}+S_\ell^{\ell j}S_\ell^{\ell i}]),
\hbox{\ where\ } S_\ell[y]=\sum_iy_iS^{\ell i},
\end{array}
\ee{LamYi}
are the mappings {\rm \rf{cRkstar}} associated with $R_k$ and $S_\ell$,
\[
\phi_\T(\lambda)=\max\limits_{t\in\cT}\lambda^Tt,\;\; \phi_\R(\lambda)=\max\limits_{r\in\R}\lambda^Tr, \;\;\mbox{and}\;\;
\lambda[\{\Xi_1,...,\Xi_N\}]=[\Tr(\Xi_1);...;\Tr(\Xi_N)].
\]
\end{corollary}}

\subsubsection{Upper-bounding $\Psi_\Pi(\cdot)$} {Our next } observation is as follows (for proof, see Section \ref{lemnewsimpleproof}):
\begin{lemma}\label{newsimplelemma} Let $Y$ be a $m\times\nu$ matrix, $Q\in \bS^m_+$, and $P$ be a probability distribution on $\bR^m$ with $\Cov[P]\preceq Q$.
Let, further, $\|\cdot\|$ be a norm on $\bR^\nu$ with the unit ball $\cB_*$ of the conjugate norm $\|\cdot\|_*$ given by {\rm (\ref{ber6spec})}. Finally, let $\Upsilon=\{\Upsilon_\ell\in\bS^{f_\ell}_+\}_{\ell\leq L}$ and a  matrix $\Theta\in\bS^m$ satisfy the constraint
\begin{equation}\label{theconstraint}
\left[\begin{array}{c|c}\Theta&{1\over 2}YM\cr \hline {1\over 2}M^TY^T&\sum_\ell\cS_\ell^*[\Upsilon_\ell]\cr\end{array}\right]\succeq0
\end{equation}
(for notation, see {\rm (\ref{ber6spec}) and (\ref{LamYi})}).
Then
\begin{equation}\label{thendearmy}
\bE_{\xi\sim P}\{\|Y^T\xi\|\}\leq\Tr(Q\Theta)+\phi_\cR(\lambda[\Upsilon]).
\end{equation}
\end{lemma}

We have the following immediate consequence of Lemma \ref{newsimplelemma}.
{ \begin{corollary}\label{consenspec} Let
\begin{equation}\label{Gamma}
\Gamma(\Theta)=\max\limits_{Q\in \Pi}\Tr(Q\Theta)
\end{equation}
and
\begin{equation}\label{Psiplus}
\overline{\Psi}_\Pi(H)=\min\limits_{\{\Upsilon_\ell \}_{\ell\leq L},\Theta\in\bS^m}\left\{\Gamma(\Theta)+\phi_\cR(\lambda[\Upsilon]):
\Upsilon_\ell\succeq0\,\forall  \ell,\;
\left[\begin{array}{c|c}\Theta&{1\over 2}HM\cr \hline {1\over 2}M^TH^T&\sum_\ell \cS_\ell^*[\Upsilon_\ell ]\cr\end{array}\right]\succeq0
\right\}
\end{equation}
Then $\overline{\Psi}_\Pi(\cdot):\;\bR^{m\times\nu}\to\bR$ is efficiently computable convex upper bound on $\Psi_\Pi(\cdot)$.
\end{corollary}}
Indeed, given Lemma \ref{newsimplelemma}, the only non-evident part of the corollary is that $\overline{\Psi}_\Pi(\cdot)$ is a well-defined real-valued function, which is readily given by Lemma \ref{verylastlemma}, see Section \ref{sectechlem}.
\begin{remark}\label{rembalance} {\rm When $\Upsilon=\{\Upsilon_\ell \}_{\ell\leq L}$, $\Theta$ is a feasible solution to the right hand side problem in {\rm (\ref{Psiplus})} and $s>0$, the pair
$\Upsilon'=\{s\Upsilon_\ell \}_{\ell\leq L}$, $\Theta'=s^{-1}\Theta$ also is a feasible solution; since $\phi_\cR(\cdot)$ and $\Gamma(\cdot)$ are positively homogeneous of degree 1, we conclude that $\overline{\Psi}_\Pi$ is in fact the infimum
of the function
$$
2\sqrt{\Gamma(\Theta)\phi_\cR(\lambda[\Upsilon])} =\inf_{s>0}\left[s^{-1}\Gamma(\Theta)+s\phi_\cR(\lambda[\Upsilon])\right]
$$
over $\Upsilon,\Theta$ satisfying the constraints of the problem {\rm (\ref{Psiplus})}.\par
In addition, for every feasible solution $\Upsilon=\{\Upsilon_\ell \}_{\ell\leq L}$, $\Theta$ to the problem {\rm (\ref{Psiplus})} with  $\M[\Upsilon]:=\sum_\ell \cS_\ell^*[\Upsilon_\ell ]\succ0$, the pair
$\Upsilon,\;\widehat{\Theta}={\four}HM\M^{-1}[\Upsilon]M^TH^T$ is feasible for the problem as well and $0\preceq\widehat{\Theta}\preceq \Theta$ (Schur Complement Lemma),  so that $\Gamma(\widehat{\Theta})\leq\Gamma(\Theta)$. As a result,
\begin{equation}\label{finalrepr}
\overline{\Psi}_\Pi(H)=\inf_\Upsilon\left\{\begin{array}{r}\four\Gamma(HM\M^{-1}[\Upsilon]M^TH^T)+\phi_\cR(\lambda[\Upsilon]):\\
\Upsilon=\{\Upsilon_\ell \in\bS^{f_\ell }_+\}_{\ell\leq L},\M[\Upsilon]\succ0
\end{array}\right\}.
\end{equation}
}
\end{remark}

\paragraph{Illustration.} Consider the case when $\|u\|=\|u\|_p$ with $p\in[1,2]$, and let us apply the just described scheme for upper-bounding $\Psi_\Pi(\cdot)$, assuming $\Pi=\{V\in\bS^m_+:\,V\preceq Q\}$ for some given $Q\succ0$, so that $\Gamma(\Theta)=\Tr(Q\Theta)$, $\Theta\succeq0$. The unit ball of the norm conjugate to $\|\cdot\|$, that is, the norm $\|\cdot\|_q$, $q={p\over p-1}\in[2,\infty]$, is the {basic} spectratope (in fact, ellitope)
$$
\cB_*=\{y\in\bR^\mu:\exists r\in\cR:=\{\bR^\nu_+:\|r\|_{q/2}\leq1\}: S_\ell^2[y]\leq r_\ell, \,1\leq \ell\leq L=\nu\}, \,\,S_\ell[y]=y_\ell.
$$
As a result, $\Upsilon$'s  from Remark \ref{rembalance} are collections of $\nu$ positive semidefinite $1\times 1$ matrices, and we can identify them with $\nu$-dimensional nonnegative vectors $\upsilon$, resulting in $\lambda[\Upsilon]=\upsilon$ and $\cM[\Upsilon]=\Diag\{\upsilon\}$.  Besides this, for {\sl nonnegative} $\upsilon$ we clearly have $\phi_{\cR}(\upsilon)=\|\upsilon\|_{p/(2-p)}$. The optimization problem in (\ref{finalrepr}) now reads
$$
\overline{\Psi}_\Pi(H)=\inf\limits_{\upsilon\in\bR^\nu}\left\{\four\Tr(Q^{1/2}H\Diag^{-1}\{\upsilon\}H^TQ^{1/2})
+\|\upsilon\|_{p/(2-p)}:\upsilon>0\right\}.
$$
After setting  $a_\ell=\|\Col_\ell[Q^{1/2}H]\|_2$, (\ref{finalrepr}) becomes
 $$
 \overline{\Psi}_\Pi(H)=\inf_{\upsilon>0}\left\{{1\over 4}\sum_i {a_i^2\over \upsilon_i} + \|\upsilon\|_{p/(2-p)}\right\},
$$
resulting in
$\overline\Psi_Q(H)=\|[a_1;...;a_\nu]\|_p$. Recalling what are $a_i$'s, we end up with
\be
\Psi_\Pi(H)\leq \overline\Psi_\Pi(H):=\left\|\left[\|\Col_1[Q^{1/2}H]\|_2;\ldots;\,\|\Col_\nu[Q^{1/2}H]\|_2\right]\right\|_p.
\ee{simpleboundPsi}
Note that the bound \rf{simpleboundPsi} can be easily improved when $\xi\sim\cN(0,Q)$. Indeed, in this case $\eta=H^T\xi$ is normal with components $\eta_i\sim \N(0, a_i^2)$, $a_i=\|\Col_\ell[Q^{1/2}H]\|_2$, and therefore
$$
\begin{array}{rcl}
\Psi_\Pi(H)&=&\bE\{\|\eta\|_p\}\leq \left[\bE_{\eta}\{\|\eta\|_p^p\}\right]^{1/p}={\sqrt{2}\,\Gamma\big({p+1\over 2}\big)\over \pi^{{1\over 2p}}}\left[\sum_ia_i^p\right]^{1/p}\\
&=&{\sqrt{2}\,\Gamma\big({p+1\over 2}\big)\over \pi^{{1\over 2p}}}\|[a_1;...;a_\nu]\|_p\,=:\widetilde{\Psi}_\Pi(H)\;\big[\leq \|[a_1;...;a_\nu]\|_p\big].\\
\end{array}
$$
For instance, when $p=1$ the bound $\widetilde{\Psi}_\Pi(H)$ becomes exact and equals $\sqrt{{2\over\pi}}\sum_ia_i=\sqrt{{2\over \pi}}\Psi_\Pi(H)$.
\par
\subsubsection{Putting things together: building presumably good linear estimate}
 Corollaries \ref{corrbyproductspec} and \ref{consenspec} {imply} the following recipe for building a ``presumably good'' linear estimate:
\begin{proposition}\label{summaryprop} In the situation of Section \ref{sitgoal} and under Assumption {\bf A}, consider the convex optimization problem (for notation, see {\rm (\ref{LamYi})} and {\rm (\ref{Gamma})})
\be
\Opt&=&\min\limits_{H,\Lambda,\Upsilon,\Upsilon',\Theta}\bigg\{\phi_{\cT}(\lambda[\Lambda])+\phi_{\cR}(\lambda[\Upsilon])+\phi_\cR(\lambda[\Upsilon'])+\Gamma(\Theta):\nn
&&\left.\begin{array}{r}\Lambda=\{\Lambda_k\succeq0,k\leq K\},\;\Upsilon=\{\Upsilon_\ell\succeq0,\ell\leq L\},\;\Upsilon'=\{\Upsilon'_\ell\succeq0,\ell\leq L\}, 
\\
\left[\begin{array}{c|c}\sum_k\R_k^*[\Lambda_k]&\half [B^T-A^TH]M\cr\hline\half M^T[B-H^TA]&\sum_\ell\S_\ell^*[\Upsilon_\ell]\cr\end{array}\right]\succeq0,\\
\left[\begin{array}{c|c}\Theta&\half HM\cr \hline \half M^TH^T&\sum_\ell \cS_\ell^*[\Upsilon'_\ell ]\cr\end{array}\right]\succeq0
\end{array}\right\}
\ee{777eq1}
 The problem is solvable, and the $H$-component $H_*$ of its optimal solution yields linear estimate $\widehat{x}_{H_*}(\omega)=H_*^T\omega$ such that
\begin{equation}\label{777eq2}
\Risk_{\Pi,\|\cdot\|}[\widehat{x}(\cdot)|\cX]\leq \Opt.
\end{equation}
\end{proposition}
The only claim in Proposition \ref{summaryprop} which is not an immediate consequence of Corollaries \ref{corrbyproductspec}, \ref{consenspec} is that problem (\ref{777eq1}) is solvable; this claim is readily given by the fact that the objective clearly is coercive on the feasible set (recall that $\Gamma(\Theta)$ is coercive on $\bS^m_+$ due to $\Pi\subset\inter \bS^m_+$ and that $y\mapsto My$ is an onto mapping, since $\cB_*$ is full-dimensional).
\begin{remark}\label{rem:multi}
{\rm 
In some applications, observations (\ref{eq1obs}) have additional structure, namely, $\omega$ is a $T$-element sample: $\omega=[\bar{\omega}_1;...;\bar{\omega}_T]$ with components
$$
\bar{\omega}_t=\bar{A}x+\xi_t,\;\;t=1,...,T,
$$
and $\xi_t$ are i.i.d. observation noises with {zero mean} distribution $\bar{P}$ satisfying $\bar{P}\lll\bar{\Pi}$ for some convex compact set $\bar{\Pi}\subset\inter \bS^{\bar{m}}_+$. In other words, we deal with {\sl repeated observations}, where for $m=T\bar m$,
\begin{equation}\label{repeated}
A=[\underbrace{\bar{A};...;\bar{A}}_{T}]\in\bR^{m\times n}\hbox{\ for some\ }\bar{A}\in\bR^{\bar{m}\times n}, \;
\Pi=\{Q=\Diag\{\underbrace{\bar{Q},...,\bar{Q}}_{T}\},\bar{Q}\in\bar{\Pi}\}.
\end{equation}
It can be easily verified (see Section \ref{prfrepeated}) that in the case of repeated observations the optimization problem (\ref{777eq1}) responsible for the presumably good linear estimate reduces to similar problem {\sl with size independent of $T$}:
\begin{proposition}\label{proprepeated} In the case of repeated observations and under Assumption {\bf A}, the linear estimate of $Bx$ yielded by an optimal solution to problem {\rm (\ref{777eq1})} can be computed as follows. Consider the convex optimization problem
\begin{equation}\label{weconsideragain}
\begin{array}{rcl}
\overline{\Opt}&=&\min\limits_{\bar{H},\Lambda,\Upsilon,\Upsilon',\bar{\Theta}}\bigg\{\phi_{\cT}(\lambda[\Lambda])+\phi_{\cR}(\lambda[\Upsilon])+\phi_\cR(\lambda[\Upsilon'])+{1\over T}\overline{\Gamma}(\bar{\Theta}):\nn
&&\left.\begin{array}{r}\Lambda=\{\Lambda_k\succeq0,k\leq K\},\;\Upsilon=\{\Upsilon_\ell\succeq0,\ell\leq L\},\;\Upsilon'=\{\Upsilon^\prime_\ell\succeq0,\ell\leq L\},
\\
\left[\begin{array}{c|c}\sum_k\R_k^*[\Lambda_k]&\half [B^T-\bar{A}^T\bar{H}]M\cr\hline\half M^T[B-\bar{H}^T\bar{A}]&\sum_\ell\S_\ell^*[\Upsilon_\ell]\cr\end{array}\right]\succeq0,\\
\left[\begin{array}{c|c}\bar{\Theta}&\half \bar{H}M\cr \hline \half M^T\bar{H}^T&\sum_\ell \cS_\ell^*[\Upsilon^\prime_\ell ]\cr\end{array}\right]\succeq0
\end{array}\right\}
\\
\end{array}
\end{equation}
where
$$
\overline{\Gamma}(\bar{\Theta}) =\max\limits_{\bar{Q}\in\bar{\Pi}}\Tr(\bar{Q}\bar{\Theta}).
$$
The problem is solvable, and the estimate in question is yielded by the $\bar{H}$-component $\bar{H}_*$ of the optimal solution according to
$$
\widehat{x}([\bar{\omega}_1;...;\bar{\omega}_T])={1\over T}\bar{H}_*^T\sum_{t=1}^T\bar{\omega}_t,
$$
and the  upper bound, provided by Proposition \ref{summaryprop}, on the risk $\Risk_{\Pi,\|\cdot\|}[\widehat{x}(\cdot)|\cX]$ of this estimate is $\overline{\Opt}$.
\end{proposition}
}
\end{remark}
\subsection{Near-optimality in Gaussian case}\label{secnearopt}
The bound \rf{777eq2} for the risk of the linear estimate $\widehat{x}_{H_*}(\cdot)$ constructed in \rf{777eq1} can be compared to the minimax optimal risk of recovering $Bx$, $x\in\cX$, from observations corrupted by zero mean Gaussian noise with covariance matrix from $\Pi$; formally, this minimax optimal risk is defined as
\be
\RiskPinormopt[\cX]=\sup_{Q\in \Pi}\inf\limits_{\widehat{x}(\cdot)}\left[\sup_{x\in\cX}\bE_{\xi\sim\cN(0,Q)}\{\|Bx-\widehat{x}(Ax+\xi)\|\}\right]
\ee{gaussrisk}
where the infimum is taken over all estimates.

{ \begin{proposition}\label{newoptimalityprop} Under the premise {and in the notation of Proposition \ref{summaryprop},
we have }
\begin{equation}\label{eqhurra}
\Risk_{\Pi,\|\cdot\|}[\widehat{x}_{H_*}|\cX]\leq \Opt\leq 64\sqrt{(2\ln F+10\ln 2)(2\ln D+10\ln 2)}\RiskPinormopt[\cX]
\end{equation}
where
\begin{equation}\label{797eqFD}
D=\sum_kd_k,\,\,F=\sum_\ell f_\ell.
\end{equation}
Thus, the upper bound $\Opt$ on the risk $\Risk_{\Pi,\|\cdot\|}[\widehat{x}_{H_*}|\cX]$ of the presumably good linear estimate $\widehat{x}_{H_*}$ yielded by an optimal solution to optimization problem {\em \rf{777eq1}} is within logarithmic in the sizes of spectratopes $\cX$ and $\cB_*$ factor  from the Gaussian minimax risk  $\RiskPinormopt[\cX]$.
\end{proposition}}\noindent For the proof, see Section \ref{newproofs}. The key component of the proof is the following important by its own right fact (for proof, see Section \ref{lemtightnewproof}):
{ \begin{lemma}\label{lemtightnew}
Let $Y$ be an $N\times\nu$ matrix, let $\|\cdot\|$ be a norm on $\bR^\nu$ such that the unit ball $\cB_*$ of the conjugate norm is
the spectratope {\em(\ref{ber6spec})}, and let $\zeta\sim\cN(0,Q)$ for some positive semidefinite $N\times N$ matrix $Q$. Then the best upper bound on ${\psi_Q}(Y):=\bE\{\|Y^T\zeta\|\}$ yielded by Lemma \ref{newsimplelemma}, that is, the optimal
value $\Opt[Q]$ in the convex optimization problem (cf. { \em(\ref{Psiplus})})
\begin{equation}\label{797newexer1}
\begin{array}{rcl}\Opt[Q]&=&\min\limits_{\Theta,\Upsilon}\bigg\{\phi_\cR(\lambda[\Upsilon])+\Tr(Q\Theta):\,\Upsilon=\{\Upsilon_\ell\succeq0,1\leq \ell\leq L\},\,\Theta\in \bS^N,\\
&&~~~~~~~~~~~~~~~~~~~~~~~~~~~\left[\begin{array}{c|c}\Theta&\half YM\cr \hline \half M^TY^T&\sum_\ell\cS_\ell^*[\Upsilon_\ell]\cr\end{array}\right]\succeq0\bigg\}\end{array}
\end{equation}
(for notation, see Lemma \ref{newsimplelemma} and 
{\em (\ref{LamYi})}
satisfies the identity
{\small\begin{equation}\label{797newexer11}
\begin{array}{l}
\forall (Q\succeq0):\\
\Opt[Q]=\overline{\Opt}[Q]:=\min\limits_{G,\Upsilon=\{\Upsilon_\ell,\ell\leq L\}}\bigg\{\phi_\cR(\lambda[\Upsilon])+\Tr(G):\Upsilon_\ell\succeq0,
\\
\multicolumn{1}{r}{\left[\begin{array}{c|c}G&{1\over 2}Q^{1/2}YM\cr \hline {1\over 2}M^TY^TQ^{1/2}&\sum_\ell\cS_\ell^*[\Upsilon_\ell]\cr\end{array}\right]\succeq0\bigg\},}\\
\end{array}
\end{equation}}\noindent
and is a tight bound on $\psi_Q(Y)$,
namely,
 \begin{equation}\label{797newexer2}
{\psi_Q(Y)}\leq\Opt[Q]\leq 22\sqrt{2\ln F+10\ln 2}\psi_Q(Y),
 \end{equation}
 where $F=\sum_\ell f_\ell$ is the size of the spectratope {\em (\ref{ber6spec})}. Besides this,
 for all $\varkappa\geq 1$ one has
\begin{equation}\label{797newexer2proba}
\Prob_\zeta\left\{\|Y^T\zeta\|\geq {\Opt[Q]\over 4{\varkappa}}
\right\}\geq \beta_{\varkappa}:=1-{\e^{3/8}\over 2}-2F\e^{-\varkappa^2/2}.
 \end{equation}
In particular, when selecting $\varkappa=\sqrt{2\ln F+10\ln 2}$, we obtain
 \begin{equation}\label{797newexer2proba1}
\Prob_\zeta\left\{\|Y^T\zeta\|\geq {\Opt[Q]\over 4\sqrt{2\ln F+10\ln 2}}
\right\}{\geq\beta_\varkappa}=0.2100>\mbox{${3\over 16}$}.
 \end{equation}
\end{lemma}}

\subsection{Illustration: covariance matrix estimation via indirect observations}
\label{reccovmatr}
{Suppose that } we observe a sample
\begin{equation}\label{mbox}
 \{\eta_t=A\xi_t\}_{t\leq T}
\end{equation}
where $A$ is a given $m\times n$ matrix, and $\xi_1,...,\xi_T$ are sampled, independently of each other,
from zero mean Gaussian distribution with unknown covariance matrix ${\vartheta}$ known to satisfy
\begin{equation}\label{eqmatrbox}
\gamma{\vartheta}_*\preceq {\vartheta}\preceq {\vartheta}_*,
\end{equation}
where $\gamma\geq0$ and ${\vartheta}_*\succ0$ are given. Our goal is to recover the linear image $\cB(\theta)$ of $\theta$, and the norm in which recovery error is measured satisfies Assumption $\mathbf{A}$.

For the covariance estimation problem to fit the framework presented in the previous section we reformulate it as follows.
\paragraph{1.}
We represent the set $\{{\vartheta}\in\bS^n_+:\gamma{\vartheta}_*\preceq{\vartheta}\preceq {\vartheta}_*\}$ as the image of the basic spectratope (matrix box)
$$
\cV=\{v\in\bS^n:\|v\|_{\Sh,\infty}\leq1\} \eqno{[\hbox{$\|\cdot\|_{\Sh,\infty}$: the spectral norm}]}
$$
under affine mapping: we set
 $
 {\vartheta}_0={1+\gamma\over 2}{\vartheta}_*,\,\, \sigma={1-\gamma\over 2},
 $
and treat the matrix
$$
\begin{array}{c}
v=\sigma^{-1}{\vartheta}_*^{-1/2}({\vartheta}-{\vartheta}_0){\vartheta}_*^{-1/2}\quad
\left[\Leftrightarrow {\vartheta}={\vartheta}_0+\sigma {\vartheta}_*^{1/2}v{\vartheta}_*^{1/2}\right]
\end{array}
$$
as the signal underlying our observations.
Note that a priori information \rf{eqmatrbox} on ${\vartheta}$ reduces to $v\in\cV$.
\paragraph{2.}  We
 pass from observations $\eta_t$ to ``lifted'' observations
$\eta_t\eta_t^T\in\bS^m$, so that
$$
\bE\{\eta_t\eta_t^T\}=\bE\{A\xi_t\xi_t^TA^T\}=A{\vartheta} A^T=A\underbrace{({\vartheta}_0+\sigma A{\vartheta}_*^{1/2}v{\vartheta}_*^{1/2})}_{{\vartheta}[v]}A^T,
$$
and treat as ``actual'' observations the matrices
$$
\omega_t=\eta_t\eta_t^T-A{\vartheta}_0A^T,\;\;1\leq t\leq T.
$$
We have\footnote{In our current considerations, we need to operate with linear mappings acting from $\bS^p$ to $\bS^q$. We treat $\bS^k$ as Euclidean space equipped with the Frobenius inner product
$
\langle u,v\rangle=\Tr(uv)
$
and denote linear mappings from $\bS^p$ into $\bS^q$ by capital calligraphic letters, like $\cA$, $\cQ$, etc.  Thus, $\cA$ in (\ref{observationseq}) denotes the linear mapping which, on a close inspection, maps matrix $v\in\bS^n$ into the matrix $\cA v=A[{\vartheta}[v]-{\vartheta}[0]]A^T$.}
\begin{equation}\label{observationseq}
\omega_t=\cA v+ \zeta_t \hbox{\ with\ }
\cA v=\sigma A{\vartheta}_*^{1/2}v{\vartheta}_*^{1/2}A^T
\hbox{\ and\ }
\zeta_t=\eta_t\eta_t^T-A{\vartheta}[v]A^T.
\end{equation}
Observe that random matrices $\zeta_1,...,\zeta_T$ are i.i.d. with zero mean and covariance mapping $\cC[v]$ (that of the random matrix-valued variable $\zeta=\eta\eta^T-\bE\{\eta\eta^T\}$, $\eta\sim\cN(0,A{\vartheta}[v]A^T)$) which satisfies (see Section \ref{proof:covmap} for the derivation)
\begin{equation}\label{conclude}
\forall v\in\cV:\cC[v]\preceq Q,\,\,\langle e,Q h\rangle=2\Tr(\vartheta_*A^ThA\vartheta_*A^TeA),\;\;e,h\in\bS^m.
\end{equation}
We have represented the problem of interest in the form described in Section \ref{sitgoal} and have specified all required data.
\paragraph{Numerical illustration.} Here we report on preliminary numerical experiments with the estimation problem stated above. They are restricted to the {\sl diagonal case} where $A\in\bR^{n\times n}$ is nonsingular,  $\cB(\theta)=B\theta B^T$ with $B\in\bR^{n\times n}$, $\vartheta_*=I_n$, and $\gamma=0$; our goal is to recover $\cB(\theta)\in\bS^n$ in the Frobenius (``Frobenius norm case''), or in the nuclear norm (``nuclear norm case'').
Our first observation is that the case of square nonsingular $A$ reduces immediately to the case  of direct observations $A=I_n$; to this end it suffices to treat, as observations,  vectors $\xi_t=A^{-1}\eta_t$, see (\ref{mbox}). It is easily seen that the estimate given by Proposition \ref{proprepeated} is ``intelligent enough'' to recognize this possibility. Furthermore, the case of $B\in\bR^{n\times n}$ reduces to the case of diagonal $B$: if $B=UDV^T$ is the singular value decomposition of $B$, with our $\vartheta_*$ and choice of the norm, we lose nothing when replacing $B$ with $D$, and our design again recognizes this possibility. Therefore, from the start, in our experiment we assume that $A=I_n$ and $\cB(\vartheta)=B\vartheta B^T$ with diagonal $B$, and
$$
\vartheta_0=\half I_n,\,\,\sigma=\half ,\,\,\vartheta[v]=\half I_n+\half v,\,\,\cA v=\half v.
$$
Thus, the estimation problem in question is reduced to that of recovering the matrix
$$
B\vartheta[v]B^T=\half B^2+\half Bv B
$$
from observations (\ref{observationseq}) stemming from a signal $v$ known to satisfy
$
v\in\cV=\{v\in\bS^n: v^2\preceq I_n\}.$ The outlined setup, as compared to the general one, simplifies dramatically optimization problem (\ref{weconsideragain}) (for details, see Appendix \ref{Processing}) and allows to run experiments with $n$ in the range of hundreds.
 \par
\begin{figure}[ht]
\begin{tabular}{cc}
\includegraphics[scale=0.55]{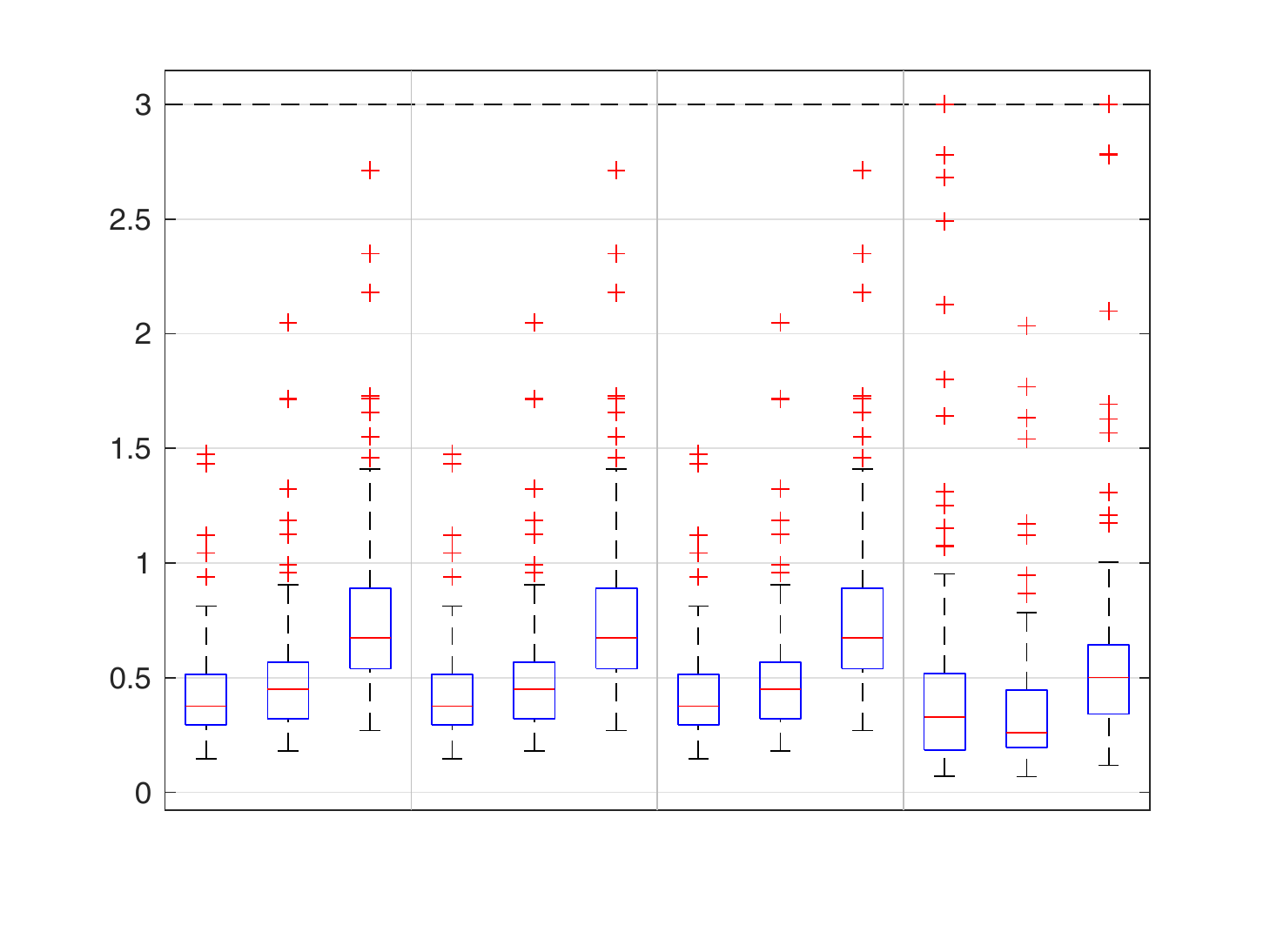}&\includegraphics[scale=0.55]{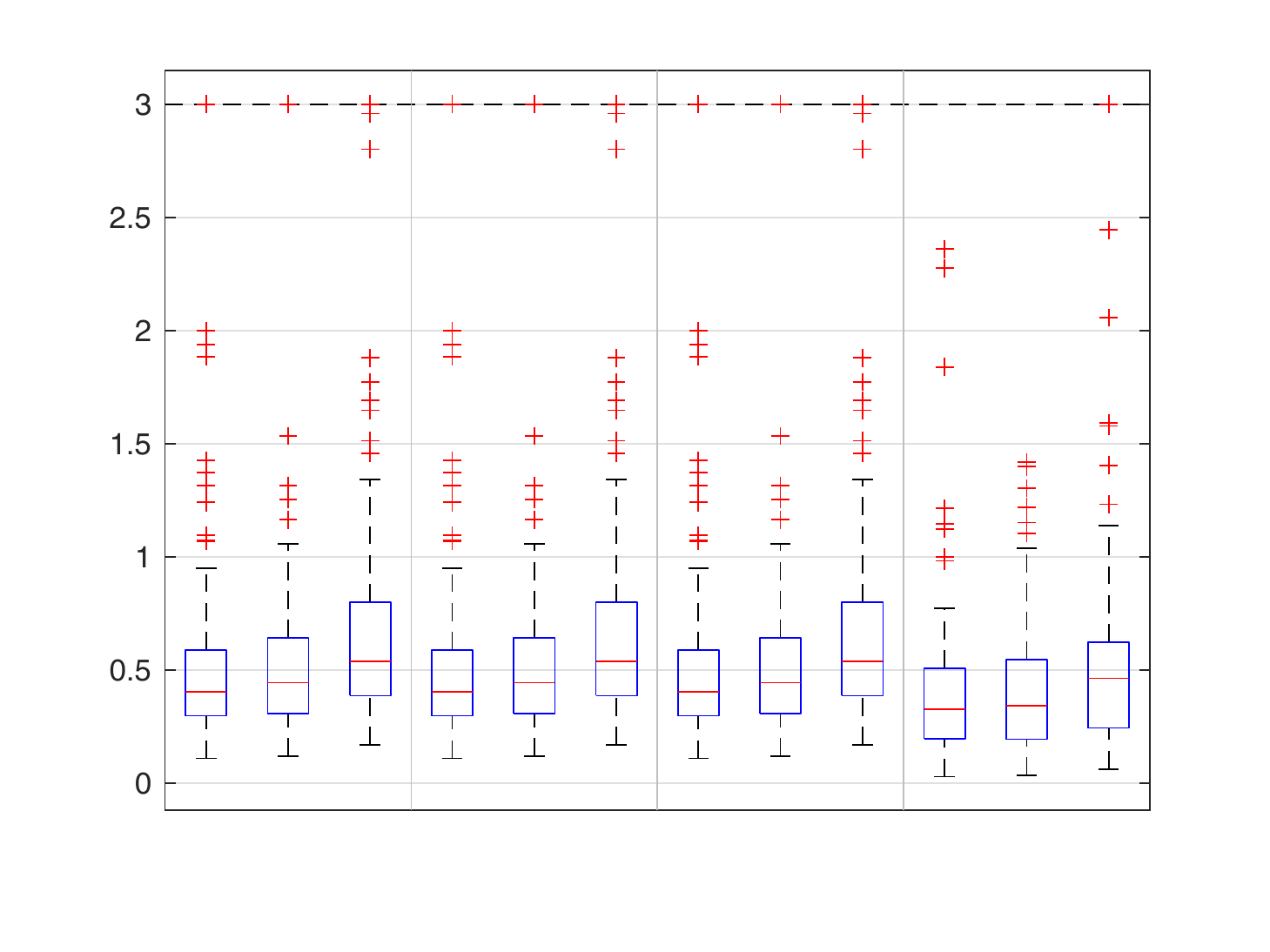}\\
$(a)$&$(b)$
\end{tabular}
\caption{\small  Boxplots for ratios of empirical errors of the linear estimation to the error of the MLE; 100 realisations of randomized experiments. $(a)$: Frobenius norm case, $(b)$: nuclear norm case.}
\label{fig:mlr-nf}
\end{figure}

 In our simulations, we use $T\in\{32,128,512\}$
 and diagonal matrix $B$ with diagonal entries $B_{ii}=i^{-\beta}$, with
 $\beta$ running through $\{0,1,2,3\}$. For every combination of $T$ and $\beta$ from the just outlined ranges,  we compute, in the Frobenius and the nuclear norm cases, the linear estimate and (the upper bound on) its risk as given by Proposition \ref{proprepeated}. Next, we run $K=100$ simulations and record the actual recovery errors as yielded by the linear estimate and the Maximum Likelihood estimate (MLE) in the role of the reference point.\footnote{It is immediately seen that with our setup the ML estimate is as follows: given observations $\eta_t$, $1\leq t\leq T$ (see (\ref{mbox}) and recall that in our case $A=I_n$), we compute the empirical covariance matrix $\wh C={1\over T}\sum_{t=1}^T\eta_t\eta_t^T$. The MLE $\widehat{\vartheta}$ of the covariance matrix $\vartheta$ of $\eta_t$'s is obtained by keeping the eigenvectors of $\wh{C}$ intact and projecting the eigenvalues of $\wh{C}$ onto $[0,1]$. The resulting MLE of $\cB(\vartheta)=B\vartheta B$ is $B\widehat{\vartheta}B$.} In our experiments, the covariance matrices underlying observations were generated as random rotations of diagonal matrices with diagonal entries drawn, independently of each other, from the uniform distribution on $[0,1]$.
 \begin{figure}[ht]
\begin{tabular}{cc}
\includegraphics[scale=0.55]{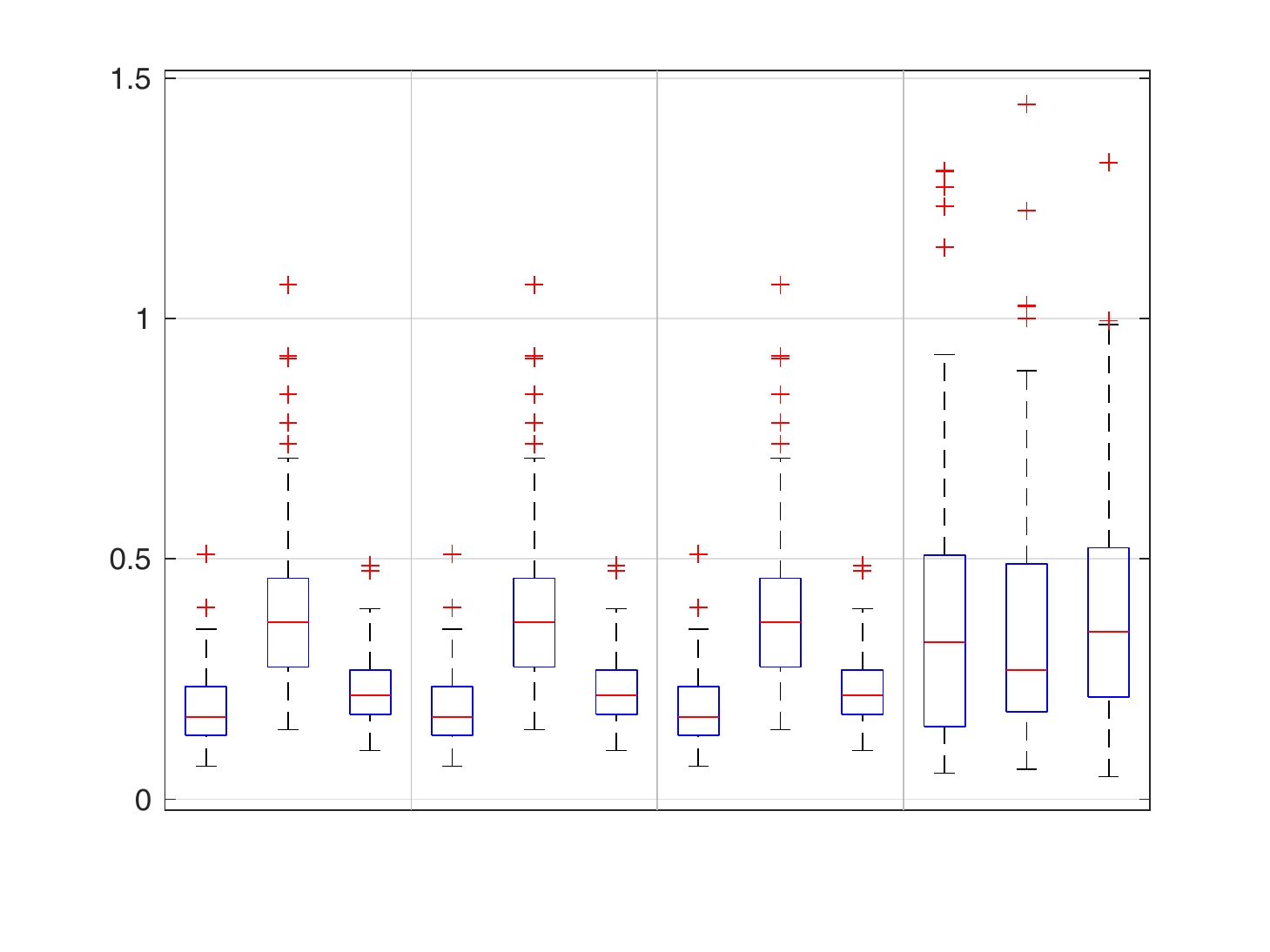}&\includegraphics[scale=0.55]{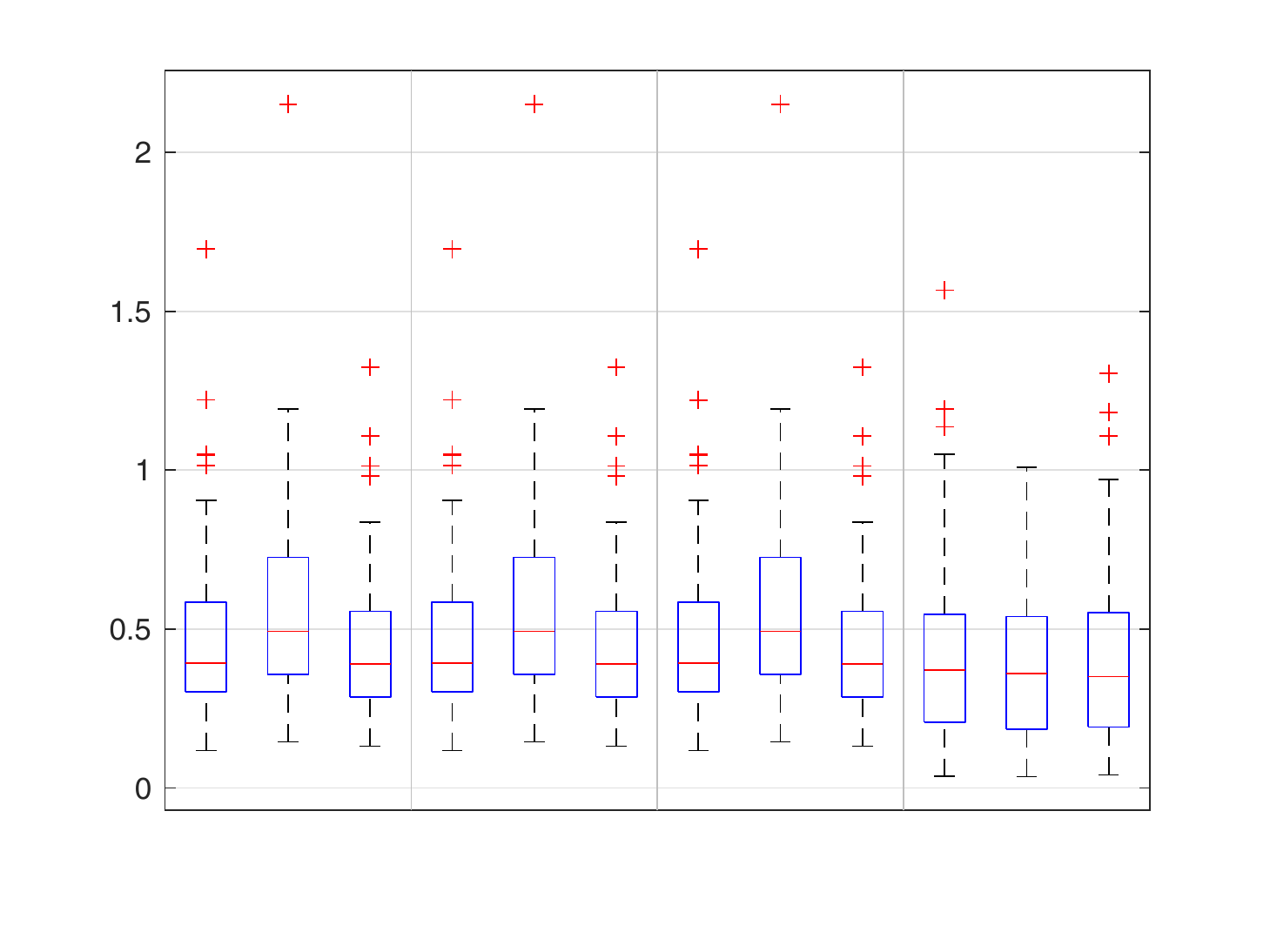}\\
$(a)$&$(b)$
\end{tabular}
\caption{\small  Ratios of empirical errors of the linear estimation to the upper risk bounds. Boxplots for 100 realisations of randomized experiments. $(a)$: Frobenius norm case, $(b)$: nuclear norm case.}
\label{fig:risk-nf}
\end{figure}
 The results of experiments with covariance matrix of size $n=128$ are presented in Figures \ref{fig:risk-nf} and \ref{fig:mlr-nf}. The first figure displays the boxplots of the ratios of actual errors of the linear estimators to the theoretical upper risk bounds for the Frobenius norm case (plot $(a)$) and nuclear norm case (plot $(b)$). Each of four boxplot groups corresponds, from left to right, to $\beta=0,\,1,\,2$ and $3$; three boxplots inside each group correspond to the observation sample lengths $T=32,\,64$ and $128$. Boxplots for ratios of errors of linear estimation to those of MLE for each simulation are displayed in Figure \ref{fig:mlr-nf}. The ``ordering'' of boxplots is the same as in Figure \ref{fig:risk-nf}; for better readability of the plots the data is ``clipped'' at the level $3$.  We see that no estimate ``uniformly dominates'' the other one, and that the linear estimate outperforms the MLE when the number $T$ of observations is relatively low.\footnote{The fact that the relative to linear estimation performance of MLE improves as $T$ grows is completely natural --  {the latter} estimate is asymptotically optimal as $T\to\infty$.}.

{\section{Linear estimation under ``uncertain-but-bounded'' noise}\label{sec:u-b-b-n}
We present here another application of the result of Proposition \ref{propmaxqf} -- construction of a linear estimate of a signal in the case of uncertain but bounded perturbation $\xi$ in the observation \rf{eq1obs}.
\subsection{Problem statement}\label{sec:ubbn-1}
Consider an estimation problem where, given an observation
\[
\omega=Ax+\xi
\]
of unknown signal $x$, known to belong to a given signal set $\X$, one wants to recover linear image $Bx$ of $x$. Here $A$ and $B$ are given $m\times n$ and $\nu\times n$ matrices.
Suppose that  all we know about $\xi$ is that it belongs to a given compact set $\cH$ (``uncertain-but-bounded observation noise''). In the situation in question, given a norm $\|\cdot\|$ on $\bR^\nu$, we quantify the accuracy of  a candidate estimate $\omega\mapsto\widehat{x}(\omega)$ by its maximal on $\X$ risk
\[
\Risksigma[\widehat{x}|\X]=\sup\limits_{x\in X,\,\xi\in\cH} \|Bx-\widehat{x}(Ax+\xi)\|
\]
(``$\cH$-risk''). \par
This is a standard problem of {\em optimal recovery} (see, e.g., \cite{MicRiv77,MicRiv84}). It is well known that when $\cH$ and $\X$ are convex compact sets,
when specifying $\widehat{x}(\omega)$ as (any) point from $\{x\in X: \omega-Ax\in\cH\}$, we get a minimax optimal, within factor 2, estimate, see also \cite{TWW83, TWW88}.} We are about to show that when $\cX$ and $\cH$ are spectratopes, and the unit ball of the norm $\|\cdot\|_*$ conjugate to $\|\cdot\|$ is a {basic} spectratope, an efficiently computable {\em linear in observation} estimate $\widehat{x}_H=H\omega$ is near-optimal in terms of its $\cH$-risk.\footnote{{In} the case where an ``efficient description'' of the sets $\cH$, $\X$ and the norm $\|\cdot\|$ is available, a minimax optimal, within factor 2, nonlinear estimate can be computed efficiently. On the other hand, its risk is generally hard to compute. Note that the linear estimate we discuss here, which comes with a ``reasonably tight'' upper bound on its risk, can be of ``numerical'' interest in the situation where estimates are to be computed repeatedly for different observations sharing common problem data -- sets $\cH$, $\X$ and the norm $\|\cdot\|$.}
\par
Our initial observation is that the situation in question reduces straightforwardly to that where there is no observation noise at all. Indeed, let $\cY=\cX\times\cH$; then $\cY$ is a spectratope, and we lose nothing when assuming that the signal underlying observation $\omega$ is $y=[x;\xi]\in \cY$:
$$
\omega=Ax+\xi=\bar{A}y,\,\bar{A}=[A,I_m],
$$
while the entity to be recovered is
$$
Bx=\bar{B}y,\,\,\bar{B}=[B,0_{\nu\times m}].
$$
With these conventions, {the observation noise vanishes, while} the $\cH$-risk of a candidate estimate $\widehat{x}(\cdot):\bR^m\to\bR^\nu$ becomes the quantity
$$
\Risk_{\|\cdot\|}[\widehat{x}|\cX\times\cH]=\sup\limits_{y=[x;\xi]\in\cX\times\cH}\|\bar{B}y-\widehat{x}(\bar{A}y)\|.
$$
To streamline the notation, let us assume that the outlined reduction has already been carried out, so the problem of interest reads: given an observation
\[
\omega=Ax\in\bR^m,
\]
estimate the linear image $Bx\in\bR^\nu$ of an unknown signal $x$ known to belong to a given spectratope $\cX$.
The risk of a candidate estimate $\widehat{x}$ is defined as
$$
\Risk_{\|\cdot\|}[\widehat{x}|\cX]=\sup_{x\in\cX} \|Bx-\widehat{x}(Ax)\|,
$$
and the norm $\|\cdot\|$ is such that the
unit ball $\cB_*$ of the norm $\|\cdot\|_*$ {\sl conjugate} to $\|\cdot\|$  is a {\em basic spectratope:}
\[
\cB_*:=\{u\in\bR^\nu:\,\|u\|_*\leq1\}=\{u\in\bR^\nu:\,\exists r\in\cR:\, S_\ell^2[u]\preceq r_\ell I_{f_\ell},\,1\leq\ell\leq L\},
\]
where the right hand side data are as required in a spectratopic representation.
By the same reasoning as in Section \ref{sect3.1}, we lose nothing when assuming from now on that {\sl the signal set is a {basic} spectratope:}
\[
\cX=\{x\in\bR^n:\, \exists t\in\cT: R_k^2[x]\preceq t_kI_{d_k},\,1\leq k\leq K\}.
\]
\subsection{Near-optimality of linear estimation}\label{euclcase}
Let $\widehat{x}_H(\omega)=H^T\omega$ be a linear estimate. We have
\bse
\Risk_{\|\cdot\|}[\widehat{x}_H|\cX]&=& \max\limits_{x\in\cX}\|(B-H^TA)x\|\\
&=&\max\limits_{[u;x]\in \cB_*\times\cX}[u;x]^T\left[\begin{array}{c|c}&\half (B-H^TA)\cr\hline\half (B-H^TA)^T&\cr\end{array}\right][u;x].\\
\ese
Applying Proposition \ref{propmaxqf}, we arrive at item (i) of the following proposition (cf. Corollary \ref{corrbyproductspec}):
\begin{proposition}\label{propbnoiseopt} In the situation of this section, consider the convex optimization problem
\be
\overline\Opt&=&\min\limits_{{H,\Upsilon=\{\Upsilon_\ell\},\Lambda=\{\Lambda_k\}}}
\bigg\{\phi_{\cR}(\lambda[\Upsilon])+\phi_{\cT}(\lambda[\Lambda]):\Upsilon_\ell\succeq0,\;\Lambda_k\succeq0,\;\forall (\ell,k)\nn
&&~~~~~~~~~~~~~~~~~~~~~~~~~~~~~~~~~~~~~~\left.\begin{array}{l}
\left[\begin{array}{c|c}\sum_\ell\cS^*_\ell[\Upsilon_\ell]&\half(B-H^TA)\cr\hline\half(B-H^TA)^T&\sum_k\cR^*_k[\Lambda_k]\cr\end{array}\right]\succeq0\\
\end{array}\right\},
\ee{considerprob}
where $\cR^*_k[\cdot]$ and $\cS^*_\ell[\cdot]$ are induced by $R_k[\cdot]$, resp., $S_k[\cdot]$, as explained in Section   \ref{spec:1}.
\item[ {\rm\ (i)}]
The problem is solvable, and the risk of the linear estimate $\widehat{x}_{H_*}(\cdot)$ yielded by the $H$-component of an optimal solution to {\rm \rf{considerprob}} does not exceed $\overline\Opt$.
\item[ {\rm\ (ii)}]The linear estimate $\widehat{x}_{H_*}$  is near-optimal in terms of its $\cH$-risk:
\begin{equation}\label{newupper}
\Risk_{\|\cdot\|}[\widehat{x}_{H_*}|\cX]\leq \overline\Opt\leq {2}\ln({2}D)\Riskopt[\cX],\;\;\;D=\sum_kd_k+\sum_\ell f_\ell,
\end{equation}
where $\Riskopt[\cX]$ is the minimax optimal risk:
$$
\Riskopt[\cX]=\inf\limits_{\widehat{x}}\Risk_{\|\cdot\|}[\widehat{x}|\cX],
$$
where $\inf$ is taken w.r.t. all possible estimates.
\end{proposition}
For proof, see Section \ref{proofofpropbnoiseopt}.
{
\subsection{Numerical illustration}\label{sec:num1}
The construction\footnote{In short, the idea of the construction is as follows.
We first note that the maximal norm $\|Bx\|$ for $x$ in the intersection of $\X$ and of the kernel of $A$, i.e., the optimal value
of the problem
\[
\max_x\left\{\|Bx\|:\,x\in \X,\,Ax=0\right\}=\max_{x,u}\left\{u^TBx,\,u\in \cB_*,\, x\in \X,\,Ax=0\right\},\eqno{(*)}
\]
lower bounds the minimax risk. Then we use semidefinite relaxation to compute a feasible solution $[\bar u;\bar x]$ to $(*)$  and
use the value $\bar u^TB\bar x$ to lower bound $\Riskopt[\cX]$. The reader is referred to the proofs of Propositions \ref{propmaxqf} and \ref{propbnoiseopt} for details.}  from the proof of Proposition \ref{propbnoiseopt}.ii can be used to lower bound numerically the minimax risk $\Risk_{\|\cdot\|}[\widehat{x}_{H_*}|\cX]$, and we can compare the resulting lower bound on the ``true'' minimax risk with the upper bound (\ref{considerprob}) on the risk of the linear estimate yielded by our approach, thus quantifying numerically {its} conservatism.
\par
We have conducted two experiments of the outlined type. In both experiments  the signal set $\X$ is the box
\[
{\X}=\{x\in\bR^n:\, j|x_j|\leq1,\,1\leq j\leq n\}\eqno{[K=n,\;R_k=k^2e_ke_k^T,k=1,..., K,\T=[0,1]^K]},
\]
$B$ is the $n\times n$ identity matrix, and ${n\over 2}\times n$ sensing matrix $A$ is a randomly rotated matrix with singular values $\lambda_j$, $1\leq j\leq n$, forming a geometric progression, with $\lambda_1=1$ and $\lambda_{n/2}=0.01$.
In the first experiment the ``dual-norm spectratope'' $\cB_*$ is a random parallelotope
 \[
 \cB^{(P)}_*=\{u\in \bR^n:\,|\eta_i^Tx|\leq 1,\,1\leq i\leq n\}. \eqno{[L=n,\,\R=[0,1]^L]}
 \]
 In the second experiment  $\cB_*$ is a random ``matrix box''
 \[
 \cB^{(M)}_*=\left\{u\in \bR^n:\,\left\|{\sum}_{i=1}^n S^ix_i\right\|_{\Sh,\infty}\leq 1\right\}, \eqno{[L=1,\,\R=[0,1]]}
 \]
 where $S^i\in \bR^{d\times d}$ are random symmetric matrices ($n=d^2/2$ in the reported experiments).
 With a natural implementations of the outlined bounding scheme we arrive at simulation results presented in Figure \ref{fig:r1}. Observe that in all experiments (100 random problems for each problem dimension) the suboptimality factor does not exceed $1.9$, while its theoretical estimation  as in \rf{newupper} varies in the interval $[9.7,\,22.2]$.

\begin{figure}[ht]
\begin{tabular}{cc}
\includegraphics[scale=0.55]{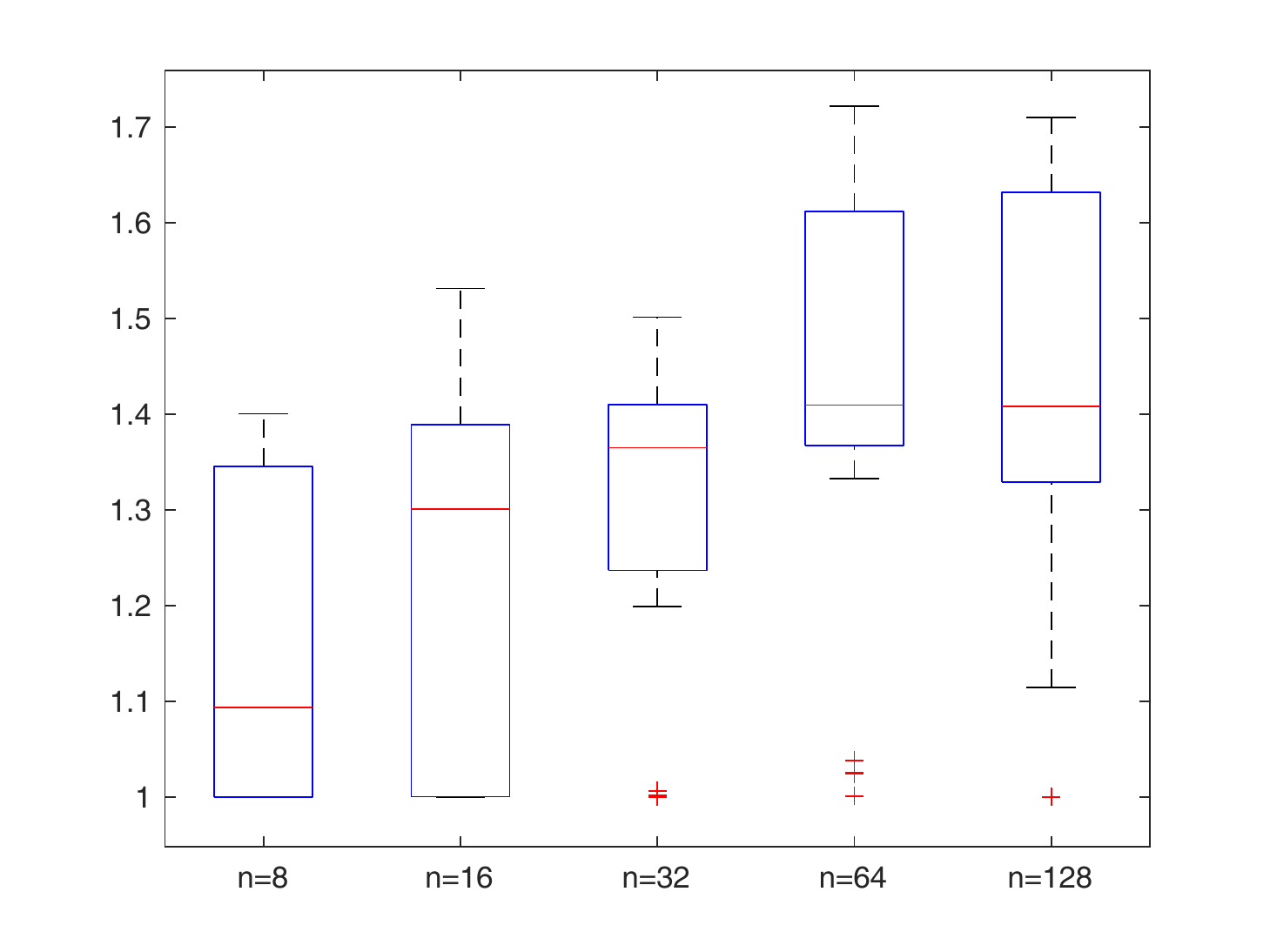}&\includegraphics[scale=0.55]{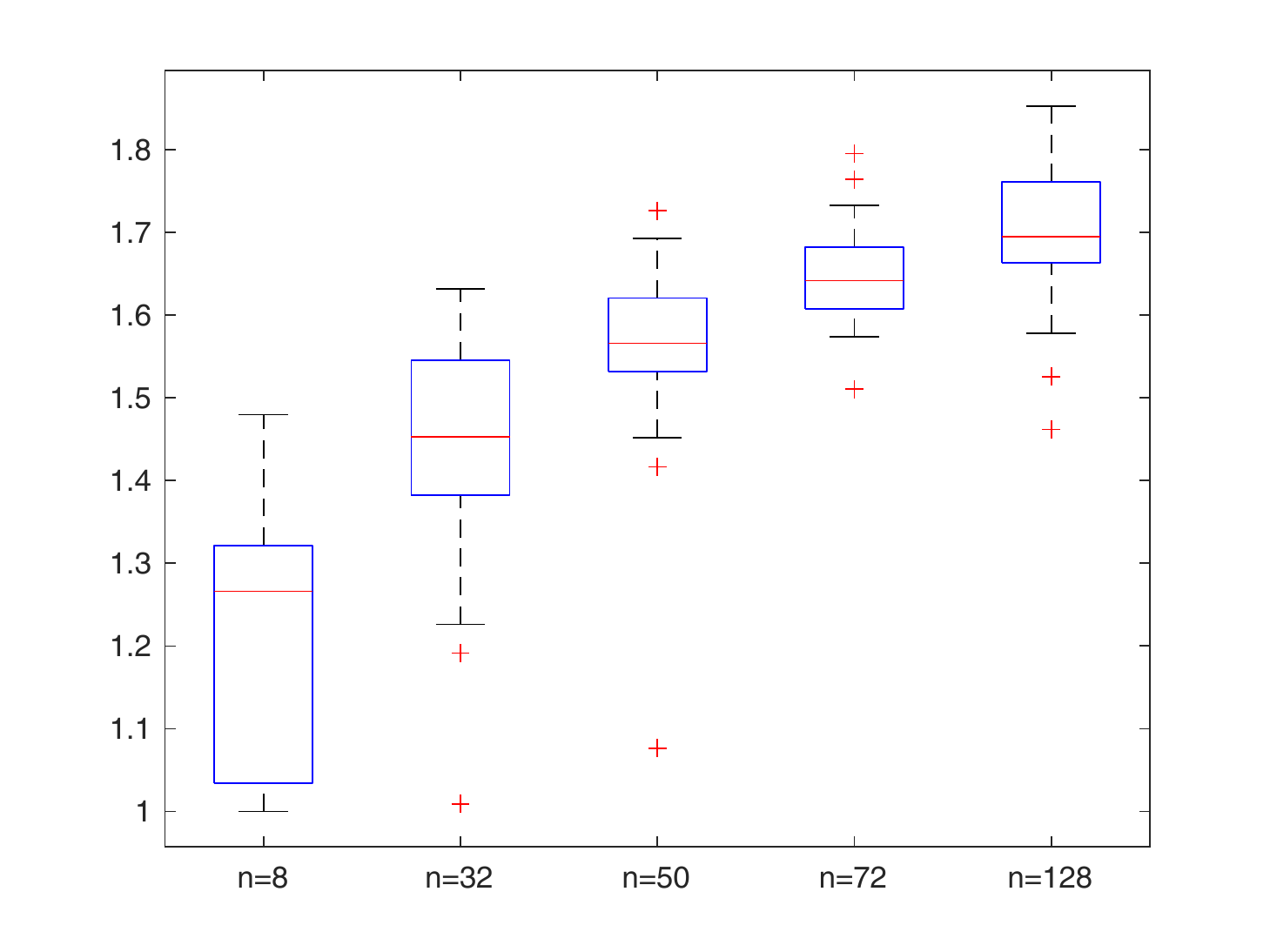}\\
$(a)$&$(b)$
\end{tabular}
\caption{\small Suboptimality factors as functions of the problem dimension. Boxplots for 100 realisations of randomized experiments. $(a)$: random parallelotope $\cB^{(P)}_*$, $(b)$: random ``matrix box'' $\cB^{(M)}_*$.}
\label{fig:r1}
\end{figure}
}

\begin{remark}\label{smixed}
{\rm {Note that Propositions \ref{summaryprop} and \ref{newoptimalityprop} also apply in} the following ``mixed'' observation scheme:
\[
\omega=Ax+\xi+\eta,
\]
where, as above, $A$ is a given $m\times n$ matrix, $x$ us unknown deterministic signal known to belong to a given signal set $\X$, $\xi$ is a random noise with distribution known to belong to a family $\cP$
of Borel probability distributions on $\bR^m$ satisfying (\ref{domination}) for a given convex compact set $\Pi\subset \inter \bS^m_+$,  and $\eta$ is ``uncertain-but-bounded'' perturbation known
to belong to a given set $\cH$. As before, our goal is to recover $Bx\in\bR^\nu$ via observation $\omega$.
Given a norm $\|\cdot\|$ on $\bR^\nu$, we can quantify the performance of a candidate estimate $\omega\mapsto\widehat{x}(\omega): \bR^m\to\bR^\nu$ by its risk
$$
\Risk_{\Pi,\cH,\|\cdot\|}[\widehat{x}|\cX]=\sup_{x\in \cX,P\lll\Pi,\eta\in \cH}\bE_{\xi\sim P} \{\|Bx-\widehat{x}(Ax+\xi+\eta)\|\}.
$$
Observe that the estimation problem associated with this ``mixed'' observation scheme straightforwardly reduces to similar problem  for random observation scheme, by the same trick we have used in Section \ref{sec:ubbn-1} to eliminate the observation noise. Indeed, let us treat $x^+:=[x;\eta]\in\cX^+:=\cX\times\cH$ and $\cX^+$ as the new signal/signal set underlying our observation, and denote $\bar{A}x^+=Ax+\eta$, $\bar{B}x^+=Bx$, where $\bar{A}=[A,I_m]$ and $\bar{B}=[B,0_{\nu\times m}]$. With these conventions, the ``mixed'' observation scheme becomes
$$
\omega=\bar{A}x^++\xi,
$$
and for every candidate estimate $\widehat{x}(\cdot)$ it clearly holds
$$
\Risk_{\Pi,\cH,\|\cdot\|}[\widehat{x}|\cX]=\Risk_{\Pi,\|\cdot\|}[\widehat{x}|\cX^+].
$$
In other words, we are now in the situation of Section \ref{sect3.1}; assuming that $\cX$ and $\cH$ are spectratopes, so is $\cX^+$, meaning that all results of Section \ref{sect1} on constructing linear estimates and their near-optimality are applicable in our present setup.
}
\end{remark}


\section{Proofs}\label{sect:Proofs}
\subsection{Technical lemma}\label{sectechlem}
In the sequel, we repeatedly use the following technical fact:
\begin{lemma}\label{verylastlemma} Given {basic} spectratope {\rm (\ref{sspectratope})},  a positive definite $n\times n$ matrix $Q$  and setting $\Lambda_k=\cR_k[Q]$, we get a collection of positive semidefinite matrices  such that $\sum_k\cR_k^*{[}\Lambda_k{]}$ is positive definite. As a corollary, whenever $M_k$, $k\leq K$, are positive definite matrices, the matrix $\sum_k\cR_k^*[M_k]$ is positive definite. In addition, the set
$$
{\cW}=\{Q\in\bS^n:Q\succeq0,\exists t\in\cT: \cR_k[Q]\preceq t_kI_{d_k},k\leq K\}
$$
is nonempty convex compact set containing a neighbourhood of the origin.
\end{lemma}
\noindent{\bf Proof.} Let us prove the first claim, Assuming the opposite, we can find a nonzero vector $y$ such that $\sum_ky^T\cR_k^*[\Lambda_k]y\leq0$, whence
$$
0\geq \sum_ky^T\cR_k^*[\Lambda_k]y=\sum_k\Tr(\cR_k^*[\Lambda_k][yy^T])=\sum_k\Tr(\Lambda_k\cR_k[yy^T])
$$
(we have used (\ref{khin23}) and (\ref{khin201})).
Since $\Lambda_k=\cR_k[Q]\succeq0$ due to $Q\succeq0$, see (\ref{khin203}), it follows that $\Tr(\Lambda_k\cR_k[yy^T])=0$ for all $k$. Now, the linear mapping $\cR_k[\cdot]$ is $\succeq$-monotone, and $Q$ is positive definite, implying that $Q\succeq r_kyy^T$ for some $r_k>0$, whence $\Lambda_k\succeq r_k\cR_k[yy^T]$. Therefore, $\Tr(\Lambda_k\cR_k[yy^T])=0$ implies that $\Tr(\cR_k^{2}[yy^T])=0$, that is, $\cR_k[yy^T]=R_k^2[y]=0$.
Since $R_k[\cdot]$ takes values in $\bS^{d_k}$, we get $R_k[y]=0$ for all $k$, which is impossible due to $y\neq 0$ and property ($S_3$), see Section \ref{spec:1}.
\par
The second  claim is an immediate consequence of the first one. Indeed, when $M_k$ are positive definite, we can find $\gamma>0$ such that $\Lambda_k\preceq\gamma M_k$ for all $k\leq K$; invoking (\ref{khin2201}), we conclude that $\cR_k^*[\Lambda_k]\preceq\gamma\cR_k^*[M_k]$, whence $\sum_k\cR^*_k[M_k]$ is positive definite along with $\sum_k\cR^*_k[\Lambda_k]$.\par Finally, the only unevident component in the last claim of the lemma is that $\cW$ is bounded. To see that it is the case, let us fix a collection $\{M_k\}$ of positive definite matrices $M_k\in\bS^{d_k}$, and let us set $M=\sum_k\cR^*_k[M_k]$, so that $M\succ0$ by already proved part of the lemma.  For $Q\in{\cW}$,
we have $\cR_k[Q]\preceq t_kI_{d_k}$, $k\leq K$, for properly selected $t\in\cT$,  so that
$$
\Tr(QM)=\sum_k\Tr(Q\cR^*_k[M_k])=\sum_k\Tr(\cR_k[Q]M_k)\leq\sum_kt_k\Tr(M_k)
$$
(we have used (\ref{khin23})), and the concluding quantity does not exceed properly selected $C<\infty$ (since $\cT$ is compact). Thus, ${\cW}\subset \{Q:Q\succeq0,\Tr(QM)\leq C\}$, whence ${\cW}$ is bounded due to $M\succ0$. \qed
\subsection{Proof of Proposition \ref{propmaxqf}}\label{proofofpropmaxqf}
\subsubsection{Preliminaries: {matrix concentration}}
 We are about to use
 {the following deep matrix concentration result, see \cite[Theorem 4.6.1]{Tropp112}:
\begin{theorem}\label{Khintchin}
Let $Q_i\in \bS^n$, $1\leq i\leq I$, and let $\xi_i$, $i=1,...,I$, be independent Rademacher ($\pm1$ with probabilities $1/2$) or $\cN(0,1)$ random variables. Then for all ${s}\geq 0$ one has
\[
\Prob\left\{\left\|{\sum}_{i=1}^I\xi_i Q_i\right\|{>} {s}\right\}\leq 2n\exp\left\{-{{s^2}\over 2v_Q}\right\}
\]
where $\|\cdot\|$ is the spectral norm, and
$
v_Q=\left\|\sum_{i=1}^IQ_i^2\right\|.
$
\end{theorem}
}
We also need the following immediate consequence of the theorem:
\begin{lemma}\label{corcorspect} Given spectratope {\rm (\ref{sspectratope})}, let $Q\in\bS^n_+$ be such that
\begin{equation}\label{khincond}
\cR_k[Q]\preceq \rho t_kI_{d_k},\,1\leq k\leq K,
\end{equation}
for some $t\in\cT$ and some $\rho\in(0,{1}]$. Then
$$
\Prob_{\xi\sim\cN(0,Q)}\{\xi\not\in\cX\}\leq \min\left[{2}D\e^{-{1\over{2}\rho}},1\right],\,\,D:=\sum_{k=1}^Kd_k.
$$
\end{lemma}
\noindent{\bf Proof.}
When setting $\xi=Q^{1/2}\eta$, $\eta\sim \N(0,I_n)$,
we have
\[
R_k[\xi]=R_k[Q^{1/2}\eta]=:\sum_{i=1}^n\eta_i\bar{R}^{ki}=\bar{R}_k[\eta]
\]
with
\[
\sum_i[\bar{R}^{ki}]^2=\bE_{\eta\sim\N(0,I_n)}\left\{\bar{R}_k^2[\eta]\right\}=
\bE_{\xi\sim\cN(0,Q)}\left\{R_k^2[\xi]\right\}=
\cR_k[Q]\preceq \rho t_kI_{d_k}
\]
 due to (\ref{khin202}).
Hence, by Theorem \ref{Khintchin} as applied with $Q_i=\bar{R}^{ki}$, $s=\sqrt{t_k}$, we get $v_Q\leq \rho t_k$ and therefore
\[
\Prob_{\xi\sim \cN(0,Q)}\{\|R_k[\xi]\|^2\geq t_k\}=\Prob_{\eta\sim\cN(0,I_n)}
\{\|\bar{R}_k[{\eta}]\|^2\geq t_k\}\leq 2d_k\e^{-{1\over 2\rho}}.
\]
We conclude that
\[
\Prob_{\xi\sim\cN(0,Q)}\{\xi\not\in\cX\}\leq\Prob_{\xi\sim \cN(0,Q)}\{\exists k: \|R_k[\xi]\|^2>t_k\}\leq 2D\e^{-{1\over 2\rho}}.\eqno{\hbox{\qed}}
\]
\subsubsection{Proving Proposition \ref{propmaxqf}}

\paragraph{1$^o$.} Under the premise of Proposition \ref{propmaxqf}, let us set $\bar{C}=P^TCP$, and
consider the conic problem
\begin{equation}\label{khinpprappeq10}
\Opt_{\#}=\max_{Q,t}\bigg\{\Tr(\bar{C}Q): Q\succeq 0, \R_k[Q]\preceq t_kI_{d_k}\,\forall k\leq K, \underbrace{[t;1]\in\bK[\cT]}_{\Leftrightarrow t\in\cT}\bigg\}.
\end{equation}
Since $\cT$ contains positive vectors, this problem is strictly feasible. Besides this, the feasible set of the problem is bounded by Lemma \ref{verylastlemma} and since $\cT$ is compact.   Thus, problem (\ref{khinpprappeq10}) is strictly feasible with bounded feasible set and thus is solvable along with its conic dual, both problems sharing a common optimal value (Conic Duality Theorem, see Appendix \ref{conicd}):
\bse
\Opt_{\#}&=&\min\limits_{\Lambda=\{\Lambda_k\}_{k\leq K},[g;s],L}\left\{s: \begin{array}{l}\Tr([\sum_k\R_k^*[\Lambda_k]-L]Q)-\sum_k[\Tr(\Lambda_k)+g_k]t_k
=\Tr(\bar{C}Q)\;\;\forall (Q,t),
\\
\Lambda_k\succeq0\,\forall k,\;L\succeq0,\;s\geq \phi_\T(-g)\\
\end{array}\right\}\\
&&\hbox{[recall that the cone dual to $\bK[\cT]$ is $\{[g;s]:s\geq\phi_{\cT}(-g)\}$]}\\
&=&\min\limits_{\Lambda,[g;s],L}\left\{s: \begin{array}{l}\sum_k\R_k^*[\Lambda_k]-L=\bar{C},g=-\lambda[\Lambda],\\
\Lambda_k\succeq0\,\forall k,\;L\succeq0,\;s\geq \phi_\T(-g)\\
\end{array}\right\}\\
&=&\min\limits_{\Lambda}\left\{\phi_\T(\lambda[\Lambda]): \sum_k\R_k^*[\Lambda_k]\succeq\bar{C},\; \Lambda_k\succeq0\,\forall k\right\}=\Opt_*
\ese
We see that (\ref{optkhinpro}) is solvable along with conic dual to problem (\ref{khinpprappeq10}), and
\[
\Opt_{\#}=\Opt_*.
\]
\paragraph{2$^o$.} Problem
(\ref{khinpprappeq10}), as we already know, is solvable; let $Q_*,t^*$ be an optimal solution to the problem. Next, let us set $R_*=Q_*^{1/2}$, $\widehat{C}=R_*\bar{C}R_*$, and let $\widehat{C}=UDU^T$ be the eigenvalue decomposition of $\widehat{C}$, so that the matrix $D=U^TR_*\bar{C}R_*U$ is diagonal, and the trace of this matrix is $\Tr(R_*\bar{C}R_*)=\Tr(\bar{C}Q_*)=\Opt_{\#}=\Opt_*$. Now let $V=R_*U$, and let $\xi=V\eta$, where $\eta\sim \R$ (i.e., $\eta$ is $n$-dimensional random (Rademacher) vector with independent entries taking values $\pm1$ with probabilities $1/2$).
We have
\begin{equation}\label{khintrace}
\xi^T\bar{C}\xi=\eta^T[V^T\bar{C}V]\eta=\eta^T[U^TR_*\bar{C}R_*U]\eta=\eta^TD\eta\equiv\Tr(D)=\Opt_*,
\end{equation}
(recall that $D$ is diagonal) and
$$
\bE_\xi\{\xi\xi^T\}=\bE_\eta\{V\eta\eta^TV^T\}=VV^T=R_*UU^TR_*=R_*^2=Q_*.
$$
From the latter relation,
\begin{equation}\label{khin44}
\bE_\xi\left\{R_k^2[\xi]\right\}=\bE_\xi\left\{\cR_k[\xi\xi^T]\right\}=\cR_k[\bE_\xi\{\xi\xi^T\}]=\cR_k[Q_*]\preceq t^*_kI_{d_k},1\leq k\leq K.
\end{equation}
On the other hand, with properly selected symmetric matrices $\bar{R}^{ki}$ we have
$$
\bar{R}_k[y]:=R_k[Vy]=\sum_{i}\bar{R}^{ki}y_i
$$
identically in $y\in\bR^n$, whence
$$
\bE_\xi\left\{R_k^2[\xi]\right\}=\bE_\eta\left\{R_k^2[V\eta]\right\}=\bE_\eta\left\{\left[{\sum}_i\eta_i\bar{R}^{ki}\right]^2\right\}=\sum_{i,j}\bE_{\eta}\{\eta_i\eta_j\}\bar{R}^{ki}\bar{R}^{kj}=\sum_i[\bar{R}^{ki}]^2.
$$
This combines with (\ref{khin44}) to imply that
\begin{equation}\label{khin45}
\sum_i[\bar{R}^{ki}]^2\preceq t^*_kI_{d_k},\,1\leq k\leq K.
\end{equation}

\paragraph{3$^o$.} Let us fix $k\leq K$. Applying Theorem \ref{Khintchin} with $Q_i=\bar{R}^{ki}$ and $s=\sqrt{t_k^*/\rho}$, we derive from (\ref{khin45}) that
\[
\Prob_{\eta\sim \R}\{\|\bar{R}_k[\eta]\|^2>t_k^*/\rho\}\leq 2d_k\e^{-{1\over 2\rho}},
\]
and recalling
the relation between $\xi$ and $\eta$, we arrive at
\begin{equation}\label{xietarho}
\Prob\{\xi:\|R_k[\xi]\|^2>t_k^*/\rho\}\leq {2}d_k{\e^{-{1\over 2\rho}}}\;\;\forall \rho\in(0,{1}].
\end{equation}
Now let us set $\bar{\rho}={1\over {2}\max[\ln({2}D),1]}$, and let $\rho\in(0,\bar{\rho})$. For this $\rho$, the sum over $k\leq K$ of the right hand sides in inequalities (\ref{xietarho}) is $<1$, implying that there exists a realization $\bar{\xi}$ of $\xi$ such that
$$
\|R_k[\bar{\xi}]\|^2\leq t_k^*/\rho,\,\,\forall k,
$$
or, equivalently,
$$
\bar{x}:=\rho^{1/2}P\bar{\xi}\in\cX,
$$
and
$$
\Opt\geq   \bar{x}^TC\bar{x}=\rho\xi^T\bar{C}\xi=\rho\Opt_*
$$
(the concluding equality is due to (\ref{khintrace})). The resulting inequality holds true for every $\rho\in(0,\bar{\rho})$, and we arrive at the right inequality in (\ref{boundskhin}). \qed
\subsection{Derivation of relation \rf{conclude} of Section \ref{reccovmatr}}\label{proof:covmap}
Let us $\succeq$-upper-bound the covariance mapping $\cC[v]$ of $\zeta=\eta\eta^T-\bE\{\eta\eta^T\}$, $\eta\sim\cN(0,A{\vartheta}[v]A^T)$. Observe that $\cC[v]$ is a symmetric linear mapping of $\bS^m$ into itself given by
\[\langle h,\cC[v]h\rangle =\bE\{\langle h,\zeta\rangle^2\}=\bE\{\langle h,\eta\eta^T\rangle^2\}-\langle h,\bE\{\eta\eta^T\}\rangle^2,\;\;h\in\bS^m.
\]
Given $v\in\cV$, setting $\theta=\vartheta[v]$, so that $0\preceq\theta\preceq\vartheta_*$, and denoting $\cH(h)=\theta^{1/2}A^ThA\theta^{1/2}$, we obtain
$$
\begin{array}{rcl}
\langle h,\cC[v]h\rangle&=&\bE_{\xi\sim\cN(0,\theta)}\{\Tr^2(hA\xi\xi^TA^T)\}-\Tr^2(h\bE_{\xi\sim\cN(0,\theta)}\{A\xi\xi^TA^T\})\\
&=&\bE_{\chi\sim\cN(0,I_n)}\{\Tr^2(hA\theta^{1/2}\chi\chi^T\theta^{1/2}A^T)\}-\Tr^2(hA\theta A^T)\\
&=&\bE_{\chi\sim\cN(0,I_n)}\{(\chi^T\cH(h)\chi)^2\}-\Tr^2(\cH(h)).
\end{array}
$$
We have
$\cH(h)=U\Diag\{\lambda\}U^T$ with orthogonal $U$, so that for $\bar{\chi}=U^T\chi$ we get
$$
\begin{array}{l}
\bE_{\chi\sim \cN(0,I_n)} \{(\chi^T\cH(h)\chi)^2\}-\Tr^2(\cH(h))=\bE_{\bar{\chi}\sim\cN(0,I_n)}\{(\bar{\chi}^T\Diag\{\lambda\}\bar{\chi})^2\}-
(\sum_i\lambda_i)^2\\
=\bE_{\bar{\chi}\sim\cN(0,I_n)}\{(\sum_i\lambda_i\bar{\chi}_i^2)^2\}-(\sum_i\lambda_i)^2=
\sum_{i\neq j}\lambda_i\lambda_j+3\sum_i\lambda_i^2-(\sum_i\lambda_i)^2=2\sum_i\lambda_i^2=2\Tr([\cH(h)]^2).\\
\end{array}
$$
Thus,
$$
\begin{array}{rcl}
\langle h,\cC[v]h\rangle&=&2\Tr([\cH(h)]^2)=2\Tr(\theta^{1/2}A^ThA\theta A^ThA\theta^{1/2})\\
&\leq& 2\Tr(\theta^{1/2}A^ThA\vartheta_* A^ThA\theta^{1/2}) \hbox{\ [since $0\preceq\theta\preceq\vartheta_*$]}\\
&=& 2\Tr(\vartheta_*^{1/2}A^ThA\theta A^ThA\vartheta_*^{1/2})\leq 2\Tr(\vartheta_*^{1/2}A^ThA\vartheta_*A^ThA\vartheta_*^{1/2})\\
&=&2\Tr(\vartheta_*A^ThA\vartheta_*A^ThA),
\end{array}
$$
what implies \rf{conclude}.  \qed
\subsection{Proof of Lemma \ref{newsimplelemma}}\label{lemnewsimpleproof}
 In the case of (\ref{theconstraint}), we have
$$
\begin{array}{rcl}
\|Y^T\xi\|&=&\max\limits_{z\in\cB_*}z^TY^T\xi=\max\limits_{y\in\cY}y^TM^TY^T\xi
\underbrace{\leq}_{\mbox{\small by \rf{theconstraint}}} \max\limits_{y\in\cY}\left[\xi^T\Theta\xi+\sum_\ell y^T\cS^*_\ell [\Upsilon_\ell ]y\right]\hbox{\ }\\
&=&\max\limits_{y\in\cY}\left[\xi^T\Theta\xi+\sum_\ell \Tr(\cS^*_\ell [\Upsilon_\ell ]yy^T)\right]
\underbrace{=}_{\mbox{\small by (\ref{khin201}) and (\ref{khin23})}}\max\limits_{y\in\cY}\left[\xi^T\Theta\xi+\sum_\ell \Tr(\Upsilon_\ell S_\ell ^2[y])\right]\\
&\underbrace{=}_{\mbox{\small by (\ref{ber6spec})}}&\xi^T\Theta\xi+\max\limits_{y,r}\left\{\sum_\ell \Tr(\Upsilon_\ell S_\ell ^2[y]):\,S_\ell ^2[y]\preceq r_\ell I_{f_\ell },\ell\leq L, r\in\cR\right\}\\
&\underbrace{\leq}_{\mbox{\small by $\Upsilon_\ell \succeq0$}}&\xi^T\Theta\xi +\max\limits_{r\in\cR}\sum_\ell \Tr(\Upsilon_\ell )r_\ell
\leq \xi^T\Theta\xi+\phi_\cR(\lambda[\Upsilon]).
\end{array}
$$
Taking expectation of  both sides of the resulting inequality w.r.t.  distribution $P$ of $\xi$ and taking into account that $\Tr(\Cov[P]\Theta)\leq\Tr(Q\Theta)$ due to $\Theta\succeq0$ (by (\ref{theconstraint})) and $\Cov[P]\preceq Q$, we get (\ref{thendearmy}). \qed
\subsection{Proof of Proposition \ref{proprepeated}}\label{prfrepeated}
In the case of (\ref{repeated}),  problem (\ref{777eq1}) reads
\begin{equation}\label{Pweq1}
\begin{array}{rcl}
\Opt&=&\min\limits_{H=[\widetilde{H}_1;...;\widetilde{H}_T],\Lambda,\Upsilon,\Upsilon',\Theta=[\theta^{t\tau}]_{1\leq t,\tau\leq T}}
\bigg\{\phi_{\cT}(\lambda[\Lambda])+\phi_{\cR}(\lambda[\Upsilon])+\phi_\cR(\lambda[\Upsilon'])+\sum_{t=1}^T\bar{\Gamma}(\theta^{tt}):\\
&&\left.\begin{array}{r}\Lambda=\{\Lambda_k\succeq0,k\leq K\},\;\Upsilon=\{\Upsilon_\ell\succeq0,\ell\leq L\},\;\Upsilon'=\{\Upsilon^\prime_\ell\succeq0,\ell\leq L\},
\\
\left[\begin{array}{c|c}\sum_k\R_k^*[\Lambda_k]&\half [B^T-\bar{A}^T\sum_{t=1}^T\widetilde{H}_t]M\cr\hline\half M^T[B-[\sum_{t=1}^T\widetilde{H}_t^T]\bar{A}]&\sum_\ell\S_\ell^*[\Upsilon_\ell]\cr\end{array}\right]\succeq0,\\
\left[\begin{array}{ccc|c}\theta^{1,1}&\cdots&\theta^{1,T}&\half \widetilde{H}_1M\cr
 \vdots&\ddots&\vdots&\vdots\cr
 \theta^{T,1}&\cdots&\theta^{T,T}&\half \widetilde{H}_TM\cr\hline
 \half M^T\widetilde{H}_1^T&\cdots&\half M^T\widetilde{H}_T^T&\sum_\ell \cS_\ell^*[\Upsilon^\prime_\ell]\cr\end{array}\right]\succeq0
\end{array}\right\},\\
\bar{\Gamma}(\theta)&=&\max_{\bar{Q}\in \bar{\Pi}}\Tr(\bar{Q}\theta),\\
\end{array}
\end{equation}
where $\widetilde{H}_t$ are $\bar{m}\times\nu$ matrices, and $\theta^{t\tau}=[\theta^{\tau t}]^T$, $1\leq t,\tau\leq T$, form a partition of $\Theta\in\bS^{\bar{m}T}$ into $\bar{m}\times\bar{m}$ blocks. Problem (\ref{Pweq1}) clearly admits a group of symmetries: a permutation $\sigma$ of $\{1,...,T\}$ induces the transformation on the space of decision variables which keeps $\Lambda,\Upsilon,\Upsilon'$ intact and maps $\widetilde{H}_t$ into $\widetilde{H}_{\sigma(t)}$, and $\theta^{t\tau}$ into $\theta^{\sigma(t)\sigma(\tau)}$; this transformation preserves  the feasible set and keeps intact the value of the objective. Since the problem is convex and solvable, it admits a ``symmetric'' optimal solution -- one with $\widetilde{H}_t=\widetilde{H}$, and $\theta^{tt}=\theta$, $1\leq t\leq T$. From the concluding semidefinite constraint in (\ref{Pweq1}) it follows that
$$
\left[\begin{array}{c|c}\theta&\half \widetilde{H}M\cr\hline
\half \widetilde{H}^TM^T&\sum_\ell \cS_\ell^*[\Upsilon^\prime_\ell]\cr\end{array}\right]\succeq0
$$
where $\Upsilon'$ stems from the symmetric solution in question. We conclude that $\Opt\geq \overline{\Opt}$, where
\begin{equation}\label{Pweq2}
\begin{array}{rcl}
\overline{\Opt}&=&\min_{\widetilde{H}\in\bR^{\bar{m}\times\nu},\Lambda,\Upsilon,\Upsilon',\theta}
\bigg\{\phi_{\cT}(\lambda[\Lambda])+\phi_{\cR}(\lambda[\Upsilon])+\phi_\cR(\lambda[\Upsilon'])+T\bar{\Gamma}(\theta):\\
&&\left.\begin{array}{r}\Lambda=\{\Lambda_k\succeq0,k\leq K\},\;\Upsilon=\{\Upsilon_\ell\succeq0,\ell\leq L\},\;\Upsilon'=\{\Upsilon^\prime_\ell\succeq0,\ell\leq L\},
\\
\left[\begin{array}{c|c}\sum_k\R_k^*[\Lambda_k]&\half [B^T-T\bar{A}^T\widetilde{H}]M\cr\hline\half M^T[B-T\widetilde{H}^T\bar{A}]&\sum_\ell\S_\ell^*[\Upsilon_\ell]\cr\end{array}\right]\succeq0,\\
\left[\begin{array}{c|c}\theta&\half \widetilde{H}M\cr\hline
 \half M^T\widetilde{H}^T&\sum_\ell \cS_\ell^*[\Upsilon^\prime_\ell]\cr\end{array}\right]\succeq0
\end{array}\right\}.
\end{array}
\end{equation}
 It is immediately seen that a feasible solution $\widetilde{H}\in\bR^{\bar{m}\times\nu},\Lambda,\Upsilon,\Upsilon',\theta$ to the optimization problem in (\ref{Pweq2}) gives rise to a feasible solution to (\ref{Pweq1}) with the same value of the objective, specifically, the solution
 $H=[\widetilde{H};...;\widetilde{H}],\Lambda,\Upsilon,\Upsilon',\Theta=[\theta^{t\tau}=\theta]_{t,\tau\leq T}$. We conclude that $\Opt=\overline{\Opt}$, and an optimal solution $\widetilde{H}\in\bR^{\bar{m}\times\nu},\Lambda,\Upsilon,\Upsilon',\theta$ to the optimization problem in (\ref{Pweq2}) gives rise to a symmetric optimal solution to (\ref{Pweq1}). The associated linear estimate is
 $$
 \widetilde{H}^T\sum_{t=1}^T\omega_t,
 $$
 and its risk is upper-bounded by $\Opt=\overline{\Opt}$. It remains to note that the optimization problem in (\ref{weconsideragain}) is obtained from the optimization problem in (\ref{Pweq2}) by  substituting $\widetilde{H}=T^{-1}\bar{H}$ and $\theta=T^{-2}\bar{\Theta}$. \qed

{
\subsection{Proof of Lemma \ref{lemtightnew}}\label{lemtightnewproof}
\bigskip\par\noindent\textbf{1$^o$.} Let us verify (\ref{797newexer11}). When $Q\succ0$, passing from variables $(\Theta,\Upsilon)$ in problem (\ref{797newexer1}) to the variables $(G=Q^{1/2}\Theta Q^{1/2},\Upsilon)$, the problem becomes exactly the optimization problem in (\ref{797newexer11}), implying that $\Opt[Q]=\overline{\Opt}[Q]$ when $Q\succ 0$. As it is easily seen, both sides in this equality are continuous in $Q\succeq0$, and (\ref{797newexer11}) follows.
\bigskip\par\noindent\textbf{2$^o$.}  Let us set $\zeta=Q^{1/2}\eta$ with $\eta\sim\cN(0,I_N)$ and $Z=Q^{1/2}Y$.
Let us show that when $\varkappa\geq1$ one has
\be
\Prob_{\eta}\{\|Z^T\eta\|\geq \bar\delta\}\geq \beta_\varkappa:=1-{\e^{3/8}\over 2}-2F\e^{-\varkappa^2/2},\;\;\;\bar\delta={\Opt[Q]\over 4\varkappa},
\ee{nlower1}
where
\begin{equation}\label{7newexer1}
\begin{array}{ll}
[\overline{\Opt}[Q]=]\quad
\Opt[Q]:=\min\limits_{\Theta,\Upsilon=\{\Upsilon_\ell,\ell\leq L\}}&\bigg\{\phi_\cR(\lambda[\Upsilon])+\Tr(\Theta):\\
&\multicolumn{1}{r}{\Upsilon_\ell\succeq0,
\left[\begin{array}{c|c}\Theta&\half ZM\cr \hline \half M^TZ^T&\sum_\ell\cS_\ell^*[\Upsilon_\ell]\cr\end{array}\right]\succeq0\bigg\}}\\
\end{array}
\end{equation}
\bigskip\par\noindent\textbf{3$^o$.} Let us represent $\Opt[Q]$ as the optimal value of a conic problem. Setting
$$
\bK=\cl\{[r;s]:s>0,r/s\in\cR\},
$$
we ensure that
$$
\cR=\{r:[r;1]\in\bK\},\,\,\bK_*=\{[g;s]:s\geq \phi_\cR(-g)\},
$$
where $\bK_*$ is the cone dual to $\bK$. Consequently, (\ref{7newexer1}) reads
$$
\Opt[Q]=\min\limits_{\Theta,\Upsilon,\theta}\left\{\theta+\Tr(\Theta):\begin{array}{ll}\Upsilon_\ell\succeq0,1\leq \ell\leq L&(a)\\
\left[\begin{array}{c|c}\Theta&\half ZM\cr \hline \half M^TZ^T&\sum_\ell\cS_\ell^*[\Upsilon_\ell]\cr\end{array}\right]\succeq0&(b)\\
{[-\lambda[\Upsilon];\theta]}\in\bK_*&(c)\\
\end{array}\right\}.\eqno{(P)}
$$
\bigskip\par\noindent\textbf{4$^o$.} Now let us prove
 that there exists matrix $W\in\bS^q_+$ and $r\in\cR$ such that
\begin{equation}\label{7newexer3}
\cS_\ell[W]\preceq r_\ell I_{f_\ell},\,\ell\leq L,
\end{equation}
and
\begin{equation}\label{7newexer4}
\Opt[Q]{\leq} \sum_i\sigma_i(ZMW^{1/2}),
\end{equation}
where $\sigma_1(\cdot)\geq\sigma_2(\cdot)\geq...$ are singular values.
\par
To get the announced result, let us pass from problem $(P)$ to its conic dual.  Applying Lemma \ref{verylastlemma} we conclude that $(P)$ is strictly feasible; in addition, $(P)$ clearly is bounded, so that the dual to $(P)$ problem $(D)$ is solvable with optimal value $\Opt[Q]$. Let us build $(D)$. Denoting by
$\Lambda_\ell\succeq0,\ell\leq L$, $\left[\begin{array}{c|c}G&-R\cr
\hline -R^T&W\cr\end{array}\right]\succeq 0$, $[r;\tau]\in\bK$ the Lagrange multipliers for the respective constraints in $(P)$, and aggregating these constraints, the multipliers being the aggregation weights,
we arrive at the following aggregated constraint:
$$
\begin{array}{l}
\Tr(\Theta G)+\Tr(W\sum_\ell\cS_\ell^*[\Upsilon_\ell])+\sum_\ell\Tr(\Lambda_\ell\Upsilon_\ell)-\sum_\ell r_\ell\Tr(\Upsilon_\ell)+\theta\tau \geq \Tr(ZMR^T).
\end{array}
$$
To get the dual problem, we impose on the Lagrange multipliers, in addition to the initial conic constraints like $\Lambda_\ell\succeq0$, $1\leq \ell\leq L$, the restriction that the left hand side in the aggregated constraint, identically in $\Theta$, $\Upsilon_\ell$ and $\theta$, is  equal to the objective of $(P)$, that is,
$$
G=I,\;\cS_\ell[W]+\Lambda_\ell -r_\ell I_{f_\ell}=0,\;1\;\leq \ell\leq L,\;\tau=1,
$$
and maximize, under the resulting restrictions, the right-hand side of the aggregated constraint. After immediate simplifications, we arrive at
$$
\Opt[Q]=\max\limits_{W,R,r}\left\{\Tr(ZMR^T):\; W\succeq R^TR,r\in\cR, \cS_\ell[W]\preceq r_\ell I_{f_\ell},\,1\leq \ell\leq L\right\}
$$
(note that $r\in\cR$ is equivalent to $[r;1]\in\bK$, and $W\succeq R^TR$ is the same as $\left[\begin{array}{c|c}I&-R\cr\hline
-R^T&W\cr\end{array}\right]\succeq0$). Now, to say that  $R^TR\preceq W$ is exactly the same as to say that $R=SW^{1/2}$ with the spectral norm $\|S\|_{\Sh,\infty}$ of $S$ not exceeding 1, so that
{\small$$
\Opt[Q]=\max\limits_{W,S,r}\bigg\{\underbrace{\Tr([ZM[SW^{1/2}]^T)}_{=\Tr([ZMW^{1/2}]S^T)}: W\succeq0,\|S\|_{\Sh,\infty}\leq1,r\in\cR, \cS_\ell[W]\preceq r_\ell I_{f_\ell},\,\ell\leq L\bigg\}
$$}\noindent
and we can immediately eliminate the $S$-variable, using the well-known fact that for every $p\times q$ matrix $J$, it holds
$$
\max\limits_{S\in\bR^{p\times q},\|S\|_{\Sh,\infty}\leq1}\Tr(JS^T)=\|J\|_{\Sh,1},$$
where $\|J\|_{\Sh,1}$ is the nuclear norm (the sum of singular values) of $J$.
We arrive at
$$
\Opt[Q]=\max\limits_{W,r}\left\{\|ZMW^{1/2}\|_{\Sh,1}:r\in\cR, W\succeq0, \cS_\ell[W]\preceq r_\ell I_{d_\ell},\,\ell\leq L\right\}.
$$
The resulting problem clearly is solvable, and its optimal solution $W$ ensures the target relations (\ref{7newexer3}) and (\ref{7newexer4}).
\bigskip\par\noindent\textbf{5$^o$.} Given $W$ satisfying (\ref{7newexer3}) and (\ref{7newexer4}), let $UJV=W^{1/2}M^TZ^T$ be the singular value decomposition of $W^{1/2}M^TZ^T$, so that $U$ and $V$ are, respectively, $q\times q$ and $N\times N$ orthogonal matrices, $J$ is $q\times N$ matrix with diagonal $\sigma=[\sigma_1;...;\sigma_p],\;p=\min[q,N],$ and zero off-diagonal entries; the diagonal entries $\sigma_i$, $1\leq i\leq p$ are the singular values of  $W^{1/2}M^TZ^T$, or, which is the same, of $ZMW^{1/2}$. Therefore, we have
\begin{equation}\label{7newexercise15}
\sum_i\sigma_i\geq\Opt[Q].
\end{equation}
Now consider the following construction. Let $\eta\sim\cN(0,I_N)$; we denote by $\upsilon$ the  vector comprised of the first $p$ entries in $V\eta$; note that $\upsilon\sim \cN(0,I_p)$, since $V$ is orthogonal. We then augment, if necessary, $\upsilon$ by $q-p$ independent of each other and of $\eta$ $\cN(0,1)$ random variables to obtain a $q$-dimensional normal vector $\upsilon'\sim\cN(0,I_q)$, and set $\chi=U\upsilon'$; because $U$ is orthogonal we also have $\chi\sim\cN(0,I_q)$. Observe that
       \begin{equation}\label{7newexercise10}
       \chi^TW^{1/2}M^TZ^T\eta=\chi^TUJV\eta=[\upsilon']^TJ\upsilon=\sum_{i=1}^p\sigma_i\upsilon_i^2.
       \end{equation}
To continue we need two simple observations.
       \begin{enumerate}
       \item[\em (i)] \label{newexeritem11} {{\sl One has
       \begin{equation}\label{7newexercise11}
       \alpha:=\Prob\left\{\sum_{i=1}^p\sigma_i{\upsilon}_i^2<\four\sum_{i=1}^p\sigma_i\right\}\leq {\mathrm{e}^{3/8}\over 2}\;[=0.7275...].
       \end{equation}}\noindent
       The claim is evident when $\sigma:=\sum_i\sigma_i=0$. Now let $\sigma>0$, and let us apply the Cramer bounding scheme. Namely, given $\gamma>0$, consider the random variable
       $$
       \omega=\exp\left\{\four\gamma\sum_i\sigma_i -\gamma\sum_i\sigma_i{\upsilon}_i^2\right\}.
       $$
       Note that $\omega>0$ a.s., and is $>1$ when $\sum_{i=1}^p\sigma_i{\upsilon}_i^2<\four\sum_{i=1}^p\sigma_i$, so that
       $\alpha\leq \bE\{\omega\}$, or, equivalently, thanks to $\upsilon\sim \N(0,I_p)$,
       {\small$$
       \ln(\alpha)\leq \ln(\bE\{\omega\})=\four\gamma\sum_i\sigma_i+\sum_i\ln\left(\bE\{\exp\{-\gamma\sigma_i{\upsilon}_i^2\}\}\right)
       \leq \four\gamma\sigma -\half\sum_i\ln(1+2\gamma\sigma_i).
       $$}\noindent
       Function $-\sum_i\ln(1+2\gamma\sigma_i)$ is convex in $[\sigma_1;...;\sigma_p]\geq0$, therefore, its maximum over the simplex $\left\{\sigma_i\geq0,i\leq p,\sum_i\sigma_i=\sigma\right\}$
       is attained at a vertex, and we get
       $$
        \ln(\alpha)\leq \four\gamma\sigma -\half\ln(1+2\gamma\sigma).
$$
       Minimizing the right hand side in $\gamma>0$, we arrive at (\ref{7newexercise11}).}

       \item[\em (ii)]\label{newexeritem22} {\sl Whenever ${\varkappa}\geq{1}$, one has
       \begin{equation}\label{7newexercise12}
       \Prob\{\|MW^{1/2}\chi\|_*>{\varkappa}\}\leq {2F\exp\{-{\varkappa}^2/2\}},
       \end{equation}
       with $F$ given by { (\ref{797eqFD})}.}
       \par
        Indeed, setting ${\rho}=1/{\varkappa}^2\leq{1}$ and $\omega=\sqrt{{\rho}}W^{1/2}\chi$, we get $\omega\sim\cN(0,{\rho} W)$. Let us apply Lemma \ref{corcorspect}  to $Q={\rho} W$,  $\cR$ in the role of $\cT$, $L$ in the role of $K$, and $\cS_\ell[\cdot]$ in the role of $\cR_k[\cdot]$. Denoting
        \[\cY:=\{y:\exists r\in\cR: S^2_\ell[y]\preceq r_\ell I_{f_\ell},\ell\leq L\},
        \] we have
       $\cS_\ell[Q]={\rho}\cS_\ell[W]\preceq {\rho} r_\ell I_{f_\ell}$, $\ell\leq L$, with $r\in\cR$ (see (\ref{7newexer3})), so we are under the premise of Lemma \ref{corcorspect} (with $\cY$ in the role of $\cX$ and therefore with $F$ in the role of $D$). Applying
       the lemma, we conclude that
       $$
       \begin{array}{l}
       \Prob\left\{\chi:{\varkappa}^{-1}W^{1/2}\chi\not\in\cY\right\}\leq {2F\exp\{-1/(2{\rho})\}}={2F\exp\{-{\varkappa}^2/2\}}.
       \end{array}
       $$
       Recalling that $\cB_*=M\cY$, we see that $\Prob\{\chi:{\varkappa}^{-1}MW^{1/2}\chi\not\in\cB_*\}$ is indeed upper-bounded by the right hand side of (\ref{7newexercise12}), and (\ref{7newexercise12}) follows.
       \end{enumerate}
{Now, for   ${\varkappa}\geq{1}$, let
       $$
       E_{\varkappa}=\left\{(\chi,\eta):\; \|MW^{1/2}\chi\|_*\leq {\varkappa},\; \sum_i\sigma_i\upsilon_i^2\geq \aic{{1\over 2}}{\mbox{$1\over 4$}}\sum_i\sigma_i\right\}.
       $$
       For $(\chi,\eta)\in E_{\varkappa}$ we have
 $$
 {\varkappa}\|Z^T\eta\|\geq \|MW^{1/2}\chi\|_*\|Z^T\eta\|\geq \chi^TW^{1/2}M^TZ^T\eta=\sum_i\sigma_i \upsilon_i^2\geq \mbox{$1\over 4$}\sum_i\sigma_i\geq\four\Opt[Q],
 $$
 (we have used (\ref{7newexercise10}) and (\ref{7newexercise15})). On the other hand,
due to (\ref{7newexercise11}) and (\ref{7newexercise12}),
       \begin{equation}\label{newexer13}
       \Prob\{E_{\varkappa}\}\geq \beta_{\varkappa},
       \end{equation}
 and we arrive at \rf{797newexer2proba}. The latter relation clearly implies (\ref{797newexer2proba1}) which, in turn, implies the right inequality in (\ref{797newexer2}). \qed
}
}

{
\subsection{Proof of Proposition \ref{newoptimalityprop}}\label{newproofs}

In what follows, we use the assumptions and the notation of Proposition \ref{newoptimalityprop}.
\bigskip\par\noindent\textbf{1$^0$.} Let
\[
\Phi(H,\Lambda,\Upsilon,\Upsilon',\Theta;Q)=\phi_{\cT}(\lambda[\Lambda])+\phi_{\cR}(\lambda[\Upsilon])+\phi_\cR(\lambda[\Upsilon'])+\Tr(Q\Theta):
\cM\times\Pi\to\bR,
\]
where
$$
\begin{array}{rcl}
\cM&=&\bigg\{(H,\Lambda,\Upsilon,\Upsilon',\Theta):\;
\begin{array}{l}\Lambda=\{\Lambda_k\succeq0,k\leq K\},\\\Upsilon=\{\Upsilon_\ell\succeq0,\ell\leq L\},\;\Upsilon'=\{\Upsilon'_\ell\succeq0,\ell\leq L\}\end{array}\\
&&
~~~~~~~~~~~~~~~~~~~~~~~\left.\begin{array}{l}
\left[\begin{array}{c|c}\sum_k\R_k^*[\Lambda_k]&\half [B^T-A^TH]M\cr\hline\half M^T[B-H^TA]&\sum_\ell\S_\ell^*[\Upsilon_\ell]\cr\end{array}\right]\succeq0\\
\left[\begin{array}{c|c}\Theta&\half HM\cr \hline \half M^TH^T&\sum_\ell \cS_\ell^*[\Upsilon'_\ell ]\cr\end{array}\right]\succeq0
\end{array}\right\}
\end{array}
$$
Looking at (\ref{777eq1}), we conclude immediately that the optimal value $\Opt$ in (\ref{777eq1}) is nothing but
\begin{equation}\label{saddlepointOpt}
\Opt=\min\limits_{(H,\Lambda,\Upsilon,\Upsilon',\Theta)\in\cM}\left[\overline{\Phi}(H,\Lambda,\Upsilon,\Upsilon',\Theta):=\max\limits_{Q\in\Pi}
\Phi(H,\Lambda,\Upsilon,\Upsilon',\Theta;Q)\right].
\end{equation}
Note that the sets $\cM$ and $\Pi$ are closed and convex,  $\Pi$ is compact, and $\Phi$ is a continuous convex-concave function on $\cM\times\Pi$. In view of these observations, the fact that $\Pi\subset\inter \bS^m_+$ combines with the Sion-Kakutani Theorem to imply that $\Phi$ possesses saddle point $(H_*,\Lambda_*,\Upsilon_*,\Upsilon_*^\prime,\Theta_*;Q_*)$ ($\min$ in $(H,\Lambda,\Upsilon,\Upsilon',\Theta)$, $\max$ in $Q$) on $\cM\times \Pi$, whence  $\Opt$ is the saddle point value of $\Phi$ by (\ref{saddlepointOpt}). We conclude that for properly selected $Q_*\in\Pi$ it holds
{\small\begin{equation}\label{forproperlyselected}
\begin{array}{l}
\Opt=\min\limits_{(H,\Lambda,\Upsilon,\Upsilon',\Theta)\in\cM} \Phi(H,\Lambda,\Upsilon,\Upsilon',\Theta;Q_*)
\\
=\min\limits_{H,\Lambda,\Upsilon,\Upsilon',\Theta}\bigg\{\phi_{\cT}(\lambda[\Lambda])+\phi_{\cR}(\lambda[\Upsilon])
+\phi_\cR(\lambda[\Upsilon'])+\Tr(Q_*\Theta):\\
\hskip20pt\left.{\quad\begin{array}{l}\Lambda=\{\Lambda_k\succeq0,k\leq K\},\;\Upsilon=\{\Upsilon_\ell\succeq0,\ell\leq L\},\;\Upsilon'=\{\Upsilon'_\ell\succeq0,\ell\leq L\}\\
\left[\begin{array}{c|c}\sum_k\R_k^*[\Lambda_k]&\half [B^T-A^TH]M\cr\hline\half M^T[B-H^TA]&\sum_\ell\S_\ell^*[\Upsilon_\ell]\cr\end{array}\right]\succeq0,\\
\left[\begin{array}{c|c}\Theta&\half HM\cr \hline \half M^TH^T&\sum_\ell \cS_\ell^*[\Upsilon'_\ell ]\cr\end{array}\right]\succeq0
\end{array}}\right\}\\
=\min\limits_{H,\Lambda,\Upsilon,\Upsilon',G}\bigg\{\phi_{\cT}(\lambda[\Lambda])
+\phi_{\cR}(\lambda[\Upsilon])+\phi_\cR(\lambda[\Upsilon'])+\Tr(G):\\
\hskip20pt\left.{\quad\begin{array}{l}\Lambda=\{\Lambda_k\succeq0,k\leq K\},\;\Upsilon=\{\Upsilon_\ell\succeq0,\ell\leq L\},\;\Upsilon'=\{\Upsilon'_\ell\succeq0,\ell\leq L\}\\
\left[\begin{array}{c|c}\sum_k\R_k^*[\Lambda_k]&\half [B^T-A^TH]M\cr\hline\half M^T[B-H^TA]&\sum_\ell\S_\ell^*[\Upsilon_\ell]\cr\end{array}\right]\succeq0,\\
\left[\begin{array}{c|c}G&\half Q_*^{1/2}HM\cr \hline \half M^TH^TQ_*^{1/2}&\sum_\ell \cS_\ell^*[\Upsilon'_\ell ]\cr\end{array}\right]\succeq0
\end{array}}\right\}\\
=\min\limits_{H,\Lambda,\Upsilon}\bigg\{\phi_{\cT}(\lambda[\Lambda])+\phi_{\cR}(\lambda[\Upsilon])+\overline{\Psi}(H):\\
\hskip20pt\left.\begin{array}{l}\Lambda=\{\Lambda_k\succeq0,k\leq K\},\;\Upsilon=\{\Upsilon_\ell\succeq0,\ell\leq L\}\\
\left[\begin{array}{c|c}\sum_k\R_k^*[\Lambda_k]&\half [B^T-A^TH]M\cr\hline
\half M^T[B-H^TA]&\sum_\ell\S_\ell^*[\Upsilon_\ell]\cr\end{array}\right]\succeq0\end{array}\right\},\\
\overline{\Psi}(H):=\min\limits_{G,\Upsilon'}\bigg\{\phi_\cR(\lambda[\Upsilon'])+\Tr(G):
\Upsilon'=\{\Upsilon'_\ell\succeq0,\ell\leq L\},\;\\
\hskip80pt\left.\left[\begin{array}{c|c}G&\half Q_*^{1/2}HM\cr \hline \half M^TH^TQ_*^{1/2}&\sum_\ell \cS_\ell^*[\Upsilon'_\ell ]\cr\end{array}\right]\succeq0
\right\}\\
\end{array}
\end{equation}}\noindent
where $\Opt$ is given by (\ref{777eq1}), and the equalities are due to (\ref{797newexer1}) and (\ref{797newexer11}).

From now on we assume that the observation noise $\xi$ in observation (\ref{eq1obs}) is $\xi\sim\cN(0,Q_*)$. Besides this, we assume that $B\neq0$, since otherwise the conclusion of Proposition \ref{newoptimalityprop} is evident.
\bigskip\par\noindent\textbf{2$^0$. $\epsilon$-risk.} In Proposition \ref{newoptimalityprop}, we are speaking about $\|\cdot\|$-risk of an estimate -- the maximal, over signals $x\in \cX$, expected norm $\|\cdot\|$ of the error in recovering $Bx$; what we need to prove that the minimax optimal risk
$\RiskPinormopt[\cX]$ as given by (\ref{gaussrisk}) can be lower-bounded by a quantity ``of order of'' $\Opt$. To this end, of course, it suffices to build such a lower bound for the quantity
$$
\RiskOpt_{\|\cdot\|} := \inf\limits_{\widehat{x}(\cdot)}\left[\sup_{x\in\cX}\bE_{\xi\sim\cN(0,Q_*)}\{\|Bx-\widehat{x}(Ax+\xi)\|\}\right],
$$
since this quantity is a lower bound on $\RiskPinormopt$. Technically, it is more convenient to work with the {\sl $\epsilon$-risk} defined in terms of ``$\|\cdot\|$-confidence intervals'' rather than in terms of the expected norm of the error. Specifically, in the sequel we will heavily use the {\sl minimax $\epsilon$-risk} defined as
\begin{equation}\label{epsriskopt}
\RiskOpt_\epsilon=\inf_{\widehat{x},\rho}\left\{\rho: \Prob_{\xi\sim\cN(0,Q_*)} \{\|Bx-\widehat{x}(Ax+\xi)\|\leq\epsilon\,\,\forall x\in\cX\right\}
\end{equation}
When $\epsilon\in(0,1)$ is once for ever fixed (in the sequel, we use $\epsilon={1\over 8}$), $\epsilon$-risk lower-bounds $\RiskOpt_{\|\cdot\|}$, since by evident reasons
\begin{equation}\label{evidentreasons}
\RiskOpt_{\|\cdot\|}\geq \epsilon\RiskOpt_\epsilon.
\end{equation}
Consequently, all we need in order to prove Proposition \ref{newoptimalityprop} is to lower-bound $\RiskOpt_{{1\over 8}}$ by a ``not too small'' multiple of $\Opt$, and this is what we are achieving below.

\bigskip\par\noindent\textbf{3$^o$.} Let $W$ be a positive semidefinite $n\times n$ matrix, let $\eta\sim\cN(0,W)$ be random signal, and let $\xi\sim\cN(0,Q_*)$ be independent of $\eta$; vectors $(\eta,\xi)$ induce random vector
$$
\omega=A\eta+\xi\sim\N(0,AWA^T+Q_*).
$$
{Consider the Bayesian version of the estimation problem where given $\omega$ we are interested to recover $B\eta$.
Recall that, because $[\omega;B\eta]$ is zero mean Gaussian, the conditional expectation
$\bE_{|\omega}\{B\eta\}$ of $B\eta$ given $\omega$ is linear in $\omega$: $\bE_{|\omega}\{B\eta\}=\bar{H}^T\omega$
for some $\bar{H}$ depending on $W$ only\footnote{We have used the following standard fact: {\sl let $\zeta=[\omega;\eta]\sim\cN(0,S)$, the covariance matrix of the marginal distribution of $\omega$  being nonsingular. Then the conditional, $\omega$ given, distribution of $\eta$ is Gaussian with mean linearly depending on $\omega$ and covariance matrix  independent of $\omega$}.}. Therefore, denoting by $P_{|\omega}$ conditional, $\omega$ given, probability distribution, for any $\rho>0$ and estimate $\wh{x}(\cdot)$ one has
\begin{equation}\label{797eq1n}
\begin{array}{rcl}
\Prob_{\eta,\xi}\{\|B\eta-\widehat{x}(A\eta+\xi)\|\geq\rho\}
&=&\bE_{\omega}\big\{\Prob_{|\omega}\{\|B\eta-\widehat{x}(\omega)\|\geq\rho\}\big\}\\
&\geq&
\bE_{\omega}\big\{\Prob_{|\omega}\{\|B\eta-\bE_{|\omega}\{B\eta\}\|\geq \rho\}\big\}\\
&=&\Prob_{\eta,\xi}\{\|B\eta-\bar{H}^T(A\eta+\xi)\|\geq \rho\},
\end{array}
\end{equation}
with the inequality given by the  Anderson Lemma \cite{anderson55} as applied to the shift of the Gaussian distribution $P_{|\omega}$ by its mean.
Applying the Anderson Lemma again we get
\bse
\Prob_{\eta,\xi}\{\|B\eta-\bar{H}^T(A\eta+\xi)\|\geq \rho\}&=&
\bE_{\xi}\big\{\Prob_{\eta}\{\|(B-\bar{H}^TA)\eta-\bar{H}^T\xi\|\geq \rho\}\big\}\\
&\geq&  \Prob_{\eta}\{\|(B-\bar{H}^TA)\eta\|\geq \rho\},
\ese
and, by ``symmetric'' reasoning,
\[
\Prob_{\eta,\xi}\{\|B\eta-\bar{H}^T(A\eta+\xi)\|\geq \rho\}\geq \Prob_{\xi}\{\|\bar{H}^T\xi\|\geq \rho\}.
\]
We conclude that for any $\wh{x}(\cdot)$
\begin{equation}\label{n_lower2}
\begin{array}{l}
\Prob_{\eta,\xi}\{\|B\eta-\widehat{x}(\omega)\|\geq\rho\}\\
\quad \geq \max\Big\{
\Prob_{\eta}\{\|(B-\bar{H}^TA)\eta\|\geq \rho\},\,\Prob_{\xi}\{\|\bar{H}^T\xi\|\geq \rho\}\Big\}.\\
\end{array}
\end{equation}
\bigskip\par\noindent\textbf{{4}$^o$.}
Let $H$ be $m\times \nu$ matrix. Applying  Lemma \ref{lemtightnew} to $N=m$, $Y=\bar H$, $Q=Q_*$, we get from (\ref{797newexer2proba})
\begin{equation}\label{797eq101}
 \Prob_{\xi\sim\cN(0,Q_*)}\{\|H^T\xi\|\geq [4\varkappa]^{-1}\overline{\Psi}(\bar H)\}\geq \beta_{\varkappa}
 \end{equation}
 where $\overline\Psi(H)$ is defined by \rf{forproperlyselected}.
Similarly, applying Lemma \ref{lemtightnew} to $N=n$, $Y=(B-\bar H^TA)^T$, $Q=W$, we obtain
\begin{equation}\label{n_lower3}
 \Prob_{\eta\sim \N(0,W)}\{\|(B-\bar{H}^TA)\eta\|\geq  [4\varkappa]^{-1}\overline{\Phi}(W,\bar H)\}\geq  \beta_\varkappa
\end{equation}
where
\begin{equation}\label{n_lower4}
\begin{array}{rl}
\overline{\Phi}(W,H)=\min\limits_{\Upsilon=\{\Upsilon_\ell,\ell\leq L\},\Theta}
&\bigg\{\Tr\left(W\Theta\right) + \phi_\cR(\lambda[\Upsilon]):\;\Upsilon_\ell\succeq0\,\forall \ell,\\
&\multicolumn{1}{r}{ \left[\begin{array}{c|c}\Theta&\half [B^T-A^TH]M\cr\hline\half M^T[B-H^TA]&\sum_\ell \cS^*_\ell[\Upsilon_\ell]\cr\end{array}\right]\succeq0\bigg\}.}\\
 \end{array}
\end{equation}
Let us put
$
\rho(W,\bar H)=[8\varkappa]^{-1}  [\overline{\Psi}(\bar H)+\overline{\Phi}(W,\bar H)];$
when combining \rf{n_lower3} with \rf{797eq101}
we conclude that
\[
\max \Big\{
\Prob_{\eta}\{\|(B-\bar{H}^TA)\eta\|\geq \rho(W,\bar H)\},\,\Prob_{\xi}\{\|\bar H^T\xi\|\geq
\rho(W,\bar H)\}\Big\}\geq \beta_\varkappa,
\]
and the same inequality holds if $\rho(W,\bar H)$ is replaced with the smaller quantity
\[
\bar\rho(W)=[8\varkappa]^{-1}\inf_{H} [\overline{\Psi}(H)+\overline{\Phi}(W,H)].
\]
Now, the latter bound combines with \rf{n_lower2} to imply the following result:
\begin{lemma}\label{797lem2} Let $W$ be a positive semidefinite $n\times n$ matrix, and $\varkappa\geq 1$. Then for any estimate $\wh{x}(\cdot)$ of $B\eta$ given observation $\omega=A\eta+\xi$, one has
\[
\Prob_{\eta,\xi}\{\|B\eta-\widehat{x}(\omega)\|\geq [8\varkappa]^{-1}  \inf_H[\overline{\Psi}(H)+\overline{\Phi}(W,H)]\}\geq \beta_\varkappa=1-{\e^{3/8}\over 2}-2F\e^{-\varkappa^2/2}
\]
where $\overline\Psi(H)$ and $\overline\Phi(W,H)$ are defined, respectively, by {\em (\ref{forproperlyselected})} and
{\em (\ref{n_lower4})}.\\
In particular, for
\begin{equation}\label{barvarkappa}
\varkappa=\bar{\varkappa}:=\sqrt{2\ln F+10\ln 2}
\end{equation} the latter probability is $>3/16$.
\end{lemma}
}
\bigskip\par\noindent\textbf{{5}$^o$.} For $0<\kappa\leq1$, let us set
\begin{equation}\label{797eq242}
\begin{array}{ll}
(a)&\cW_\kappa=\{W\in\bS^n_+:\; \exists t\in\cT: \,\cR_k[W]\preceq \kappa t_kI_{d_k},1\leq k\leq K\},\\
(b)&\cZ=\hbox{\footnotesize$\left\{(\Upsilon=\{\Upsilon_\ell,\ell\leq L\},\Theta,H):\;\begin{array}{l}
\Upsilon_\ell\succeq0\,\forall\ell,\\
\left[\begin{array}{c|c}\Theta&{1\over 2}[B^T-A^TH]M\cr\hline{1\over 2} M^T[B-H^TA]&\sum_\ell \cS^*_\ell[\Upsilon_\ell]\cr\end{array}\right]\succeq0\end{array}\right\}.
$}
\end{array}
\end{equation}
Note that $\cW_\kappa$ is a nonempty convex compact (by Lemma \ref{verylastlemma}) set such that $\cW_\kappa=\kappa\cW_1$, and $\cZ$ is a nonempty closed convex set. Consider the parametric saddle point problem
\begin{equation}\label{797eq110}
\Opt(\kappa)=\max\limits_{W\in\cW_\kappa}\min\limits_{(\Upsilon,\Theta,H)\in\cZ} \left[E(W;\Upsilon,\Theta,H):=\Tr(W\Theta)+\phi_{\cR}(\lambda[\Upsilon])+\overline{\Psi}(H)\right].
\end{equation}
This problem is convex-concave; utilizing the fact that $\cW_\kappa$ is compact and contains positive definite matrices, it is immediately seen that the Sion-Kakutani theorem ensures the existence of a saddle point whenever $\kappa\in(0,1]$. We claim that
{  \begin{equation}\label{797eqsqrt}
0<\kappa\leq 1\Rightarrow  \Opt(\kappa)\geq\sqrt{\kappa}\Opt(1).
\end{equation}}\noindent
Indeed, $\cZ$ is invariant w.r.t. scalings
$$(\Upsilon=\{\Upsilon_\ell,\ell\leq L\},\Theta,H)\mapsto (\theta\Upsilon:=\{\theta\Upsilon_\ell,\ell\leq L\},\theta^{-1}\Theta,H),\eqno{[\theta>0]}.$$
When taking into account that $\phi_{\cR}(\lambda[\theta\Upsilon])=\theta\phi_{\cR}(\lambda[\Upsilon])$, we get
$$
\begin{array}{rcl}
\underline{E}(W)&:=&\min\limits_{(\Upsilon,\Theta,H)\in\cZ}E(W;\Upsilon,\Theta,H)=\min\limits_{(\Upsilon,\Theta,H)\in\cZ}\inf\limits_{\theta>0}E(W;\theta\Upsilon,\theta^{-1}\Theta,H)\\
&=&
\min\limits_{(\Upsilon,\Theta,H)\in\cZ}\left[2\sqrt{\Tr(W\Theta)\phi_{\cR}(\lambda[{\Upsilon}])} +\overline{\Psi}(H)\right].
\\
\end{array}
$$
Because $\overline{\Psi}$ is nonnegative we conclude that whenever $W\succeq0$ and $\kappa\in(0,1]$, one has
$$
\underline{E}(\kappa W)\geq \sqrt{\kappa}\underline{E}(W),
$$
which combines with $\cW_\kappa=\kappa\cW_1$ to imply that
$$
\Opt(\kappa)=\max_{W\in\cW_\kappa}\underline{E}(W)=\max_{W\in\cW_1}\underline{E}(\kappa W)\geq \sqrt{\kappa}\max_{W\in\cW_1}\underline{E}(W)=\sqrt{\kappa}\Opt(1),
$$
and (\ref{797eqsqrt}) follows.
\bigskip\par\noindent\textbf{{6}$^o$.} We claim that
{  \begin{equation}\label{797eq201}
\Opt(1)=\Opt,
\end{equation}
where $\Opt$ is given by (\ref{777eq1}) (and, as we have seen, by (\ref{forproperlyselected}) as well). Note that (\ref{797eq201}) combines with (\ref{797eqsqrt}) to imply that
\begin{equation}\label{797eq2000}
0<\kappa\leq1 \Rightarrow \Opt(\kappa)\geq\sqrt{\kappa}\Opt.
\end{equation}}\noindent
 Verification of (\ref{797eq201}) is given by the following computation. By the Sion-Kakutani Theorem,
$$
\begin{array}{l}
\Opt(1)=\max\limits_{W\in\cW_1}\min\limits_{(\Upsilon,\Theta ,H)\in\cZ} \left\{\Tr(W\Theta )+\phi_{\cR}(\lambda[\Upsilon])+\overline{\Psi}(H)\right\}\\
=\min\limits_{(\Upsilon,\Theta ,H)\in\cZ}\max\limits_{W\in\cW_1}\left\{\Tr(W\Theta )+\phi_{\cR}(\lambda[\Upsilon])+\overline{\Psi}(H)\right\}\\
=\min\limits_{(\Upsilon,\Theta ,H)\in\cZ}\bigg\{\overline{\Psi}(H)+\phi_{\cR}(\lambda[\Upsilon])+\max\limits_W\bigg\{\Tr(\Theta W):\\
\multicolumn{1}{r}{ W\succeq0,\exists t\in\cT:\cR_k[W]\preceq t_kI_{d_k},k\leq K\bigg\}\bigg\}}\\
=\min\limits_{(\Upsilon,\Theta ,H)\in\cZ}\bigg\{\overline{\Psi}(H)+\phi_{\cR}(\lambda[\Upsilon])+\max\limits_{W,t}\bigg\{\Tr(\Theta W):\\ \multicolumn{1}{r}{W\succeq0,[t;1]\in\bT,\cR_k[W]\preceq t_kI_{d_k},k\leq K\bigg\}\bigg\},}
\end{array}
$$
where $\bT$ is the closed conic hull of $\cT$.
Now, using Conic Duality combined with the fact that $\bT_*=\{[g;s]:s\geq\phi_{\cT}(-g)\}$ we obtain
\bse
\lefteqn{\max\limits_{W,t}\left\{\Tr(\Theta W): \;W\succeq0,[t;1]\in\bK[\cT],\,\cR_k[W]\preceq t_kI_{d_k},\,k\leq K\right\}}\\
&=&\min\limits_{Z,[g;s],\Lambda=\{\Lambda_k\}}\left\{s:\,\left\{\begin{array}{l}Z\succeq0,[g;s]\in(\bK[\cT])_*,\,\Lambda_k\succeq0,k\leq K\\
-\Tr(ZW)-g^Tt+\sum_k\Tr(\cR_k^*[\Lambda_k]W)\\
\multicolumn{1}{r}{-\sum_kt_k\Tr(\Lambda_k)= \Theta}\\
\forall (W\in\bS^n,t\in\bR^K)\\
\end{array}\right.\right\}
\\
&=&\min\limits_{Z,[g;s],\Lambda=\{\Lambda_k\}}\left\{s:\,\left\{\begin{array}{l}Z\succeq0,\,s\geq\phi_{\cT}(-g),\,\Lambda_k\succeq0,\,k\leq K\\
\Theta=\sum_k\cR_k^*[\Lambda_k]-Z,\,g=-\lambda[\Lambda]\\
\end{array}\right.\right\}\\
&=&\min\limits_{\Lambda}\left\{\phi_{\cT}(\lambda[\Lambda]):\, \Lambda=\{\Lambda_k\succeq0,k\leq K\},\,\Theta\preceq \sum_k\cR_k^*[\Lambda_k]\right\},
\ese
and we arrive at
$$
\begin{array}{rl}
\Opt(1)
=\min\limits_{\Upsilon,\Theta,H,\Lambda}&\bigg\{\overline{\Psi}(H)+\phi_{\cR}(\lambda[\Upsilon])+\phi_{\cT}(\lambda[\Lambda]):\\
&\multicolumn{1}{r}{\begin{array}{l}
\Upsilon=\{\Upsilon_\ell\succeq0,\ell\leq L\},\Lambda=\{\Lambda_k\succeq0,k\leq K\},\\
\Theta\preceq \sum_k\cR^*_k[\Lambda_k],\\
\left[\begin{array}{c|c}\Theta&{1\over 2}[B^T-A^TH]M\cr\hline {1\over 2}M^T[B-H^TA]&\sum_\ell \cS^*_\ell[\Upsilon_\ell]\cr\end{array}\right]\succeq0
\end{array}\bigg\}}\\
=\min\limits_{\Upsilon,H,\Lambda}&\bigg\{\overline{\Psi}(H)+{\phi_{\cR}(\lambda[\Upsilon])}+\phi_{\cT}(\lambda[\Lambda]):\\
&\multicolumn{1}{r}{\begin{array}{l}
\Upsilon=\{\Upsilon_\ell\succeq0,\,\ell\leq L\},\Lambda=\{\Lambda_k\succeq0,k\leq K\}\\
\left[\begin{array}{c|c}\sum_k\cR^*_k[\Lambda_k]&{1\over 2}[B^T-A^TH]M\cr\hline {1\over 2}M^T[B-H^TA]&\sum_\ell \cS^*_\ell[\Upsilon_\ell]\cr\end{array}\right]\succeq0\\
\end{array}\bigg\}}\\
=\Opt&\hbox{\ [see (\ref{forproperlyselected})].}
\end{array}
$$

\bigskip\par\noindent\textbf{{{7$^o$.}}} Now we can complete the proof. For $\kappa\in(0,{1}]$, let $W_\kappa$ be the $W$-component of a saddle point solution to the saddle point problem  (\ref{797eq110}). Then, by \rf{797eq2000},
\begin{equation}\label{797eq252}
\begin{array}{rcl}
\sqrt{\kappa}\Opt\leq \Opt(\kappa)
&=&\min\limits_{(\Upsilon,\Theta,H)\in\cZ}\Big\{\Tr(W_\kappa \Theta)+\phi_{\cR}(\lambda[\Upsilon])+\overline{\Psi}(H)\Big\}\\
&=&\min\limits_H\big\{\overline{\Phi}(W_\kappa, H)+\overline{\Psi}(H)\big\}\\
\end{array}
\end{equation}
On the other hand, when applying Lemma \ref{corcorspect} to $Q=W_\kappa$ and $\rho=\kappa$, we obtain, in view of relations $0<\kappa\leq{1}$, $W_\kappa\in\cW_\kappa$,
\begin{equation}\label{797eq253}
\delta(\kappa):=\Prob_{\zeta\sim\cN(0,I_n)}\{W_\kappa^{1/2}\zeta\not\in\cX\}\leq 2D\e^{-{1\over 2\kappa}},
\end{equation}
with $D$ given by (\ref{797eqFD}). In particular, when setting
\begin{equation}\label{barkappa}
\bar{\kappa}={1\over 2\ln D+10\ln 2}
\end{equation} we obtain $\delta_\kappa\leq 1/16$.
Therefore,
\be
\Prob_{\eta\sim\cN(0,W_{\bar{\kappa}})}\{\eta\not\in\cX\}\leq \mbox{$1\over 16$}.
\ee{1/16}
Now let
\begin{equation}\label{varrhostar}
\varrho_*:= {\Opt\over 8\sqrt{(2\ln F+10\ln 2)(2\ln D+10\ln 2)}}.
\end{equation}
All we need in order to achieve our goal, that is, to justify (\ref{eqhurra}), is to show that
\begin{equation}\label{allweneedtoshow}
\RiskOpt_{{1\over 8}}\geq\varrho_*,
\end{equation}
since given the latter relation, (\ref{eqhurra}) will be immediately given by (\ref{evidentreasons}) as applied with $\epsilon={1\over 8}$.
\par
To prove (\ref{allweneedtoshow}), assume, on the contrary to what should be proved,
that the ${1\over8}$-risk is  $<\varrho_*$, and let $\bar{x}(\cdot)$ be an estimate
with ${1\over 8}$-risk $\leq \varrho_*$.  We can utilize $\bar{x}$ to estimate $B\eta$, in the Bayesian problem of recovering $B\eta$ from observation $\omega=A\eta+\xi$, $(\eta,\xi)\sim\cN(0,\Sigma)$ with $\Sigma=\Diag\{W_{\bar{\kappa}},Q_*\}$. From \rf{1/16} we conclude that
\begin{equation}\label{3/16-1}
\begin{array}{l}
\Prob_{(\eta,\xi)\sim\cN(0,\Sigma)}\{\|B\eta-\bar{x}(A\eta+\xi)\|> \varrho_*\}\\
\leq \Prob_{(\eta,\xi)\sim\cN(0,\Sigma)}
\{\|B\eta-\bar{x}(A\eta+\xi)\|> \varrho_*,\;\eta\in
\X\}+\Prob_{\eta\sim\cN(0,W_{\bar{\kappa}})}\{\eta\not\in\cX\}\\
\leq \mbox{$1\over 8$}+\mbox{$1\over 16$}=\mbox{$3\over 16$}.
\\
\end{array}
\end{equation}
On the other hand,  by \rf{797eq252} we have
\[
\min_H\, [\overline{\Phi}(W_{\bar{\kappa}},H)+\overline{\Psi}(H)]=\Opt(\bar{\kappa})\geq\sqrt{\bar{\kappa}}\Opt=
[8\bar{\varkappa}]\varrho_*
\]
with $\bar{\varkappa}$ given by (\ref{barvarkappa}),
so by Lemma \ref{797lem2}, for any estimate $\hat{x}(\cdot)$ of $B\eta$ via observation $\omega=Ax+\xi$ it holds
\[
\Prob_{\eta,\xi}\{\|B\eta-\widehat{x}(A\eta+\xi)\|\geq \varrho_*\}\geq \beta_{\bar{\varkappa}}>3/16;
\]
in particular, this relation should hold true for $\widehat{x}(\cdot)\equiv \bar{x}(\cdot)$, but the latter is impossible: the ${1\over 8}$-risk of $\bar{x}$ is $\leq\varrho_*$, see  \rf{3/16-1}. \qed
}

\subsection{Proof of Proposition \ref{propbnoiseopt}}\label{proofofpropbnoiseopt}
\paragraph{1$^o$.} Item (i) is a direct consequence of Proposition \ref{propmaxqf}, modulo the claim that problem (\ref{considerprob}) is solvable, and we start with justifying this claim. Let $F={\rm Im} A$. Clearly, feasibility of a candidate solution $(H,\Lambda,\Upsilon)$ to the problem depends solely on the restriction of the linear mapping $z\mapsto H^Tz$ onto $F$, so that adding to the constraints of the problem the requirement that the restriction of this linear mapping on the orthogonal complement of $F$ in $\bR^m$ is identically zero, we get an equivalent problem. It is immediately seen that in the resulting problem, the feasible solutions with the value of the objective $\leq a$ for every $a\in\bR$ form  a compact set, so that the latter problem (and thus  -- the original one)  indeed is solvable.\par
Let us prove the near-optimality result of (ii).
\paragraph{2$^o$.} Observe that setting
\begin{equation}\label{mathfrakR}
{\varrho}=\max_x\left\{\|Bx\|:x\in \X,Ax=0\right\},
\end{equation}
we ensure that
\begin{equation}\label{ensurethat}
\Riskopt[\cX] \geq{\varrho}.
\end{equation}
Indeed, let $\bar{x}$ be an optimal solution to the (clearly solvable) optimization problem in  (\ref{mathfrakR}). Then observation $\omega=0$ can be obtained from both the signals $x=\bar{x}$ and $x=-\bar{x}$, and therefore the risk of any (deterministic) recovery routine is at least $\|B\bar{x}\|={\varrho}$, as claimed.
\paragraph{3$^o$.} It may happen that $\Ker A=\{0\}$. In this case the situation is trivial: specifying $A^\dag$ as a partial inverse to $A$: $A^\dag A=I_n$ and setting $H^T=BA^\dag$ (so that $B-H^TA=0$), $\Upsilon_\ell=0_{f_\ell\times f_\ell}$, $\ell\leq L$, $\Lambda_k=0_{dk\times d_k}$, $k\leq K$, we get a feasible solution to the optimization problem in (\ref{considerprob}) with zero value of the objective, implying
that $\Opt_\#=0$; consequently, the linear estimate induced by an optimal solution to the problem is with zero risk, and the conclusion of Proposition \ref{propbnoiseopt} is clearly true. With this in mind, we assume from now on that $\Ker A\neq\{0\}$. Denoting ${\kappa}=\dim\Ker A$, we can build  an $n\times {\kappa}$ matrix $E$ of rank $\kappa$ such that $\Ker A$ is the image space of $E$.
\paragraph{4$^o$.} Setting
$$
\begin{array}{rcl}
\cZ&:=&\{z\in\bR^{{\kappa}}: Ez\in\cX\}=\left\{z\in\bR^{{\kappa}}: \exists (t\in\cT): \bar{R}_k^2[z]\preceq t_kI_{d_k},\,k\leq K \right\},\,\,\bar{R}_k[z]=R_k[Ez],
\\
C&=&\left[\begin{array}{c|c}&{1\over 2}BE\cr\hline{1\over 2}E^TB^T&\cr\end{array}\right],\\
\end{array}
$$
note that when $z$ runs trough the spectratope $\cZ$, $Ez$ runs exactly through the entire set $\{x\in \cX:Ax=0\}$. With this in mind,
invoking Proposition \ref{propmaxqf}, we arrive at
\be
&&{\varrho}=\max_{g:\|g\|_*\leq1}\max_{z\in\cZ} g^TBEz=\max_{[u;z]\in\cB_*\times\cZ}[u;z]^TC[u;z]\nn
&\leq&\Opt:=\min\limits_{\Upsilon=\{\Upsilon_\ell:\,\ell\leq L\},\atop
\Lambda=\{\Lambda_k,\,k\leq K\}}\bigg\{\phi_{\cR}(\lambda[\Upsilon])+\phi_{\cT}(\lambda[\Lambda]):\Upsilon_\ell{\succeq}0,\;\Lambda_k{\succeq}0,\;\forall (\ell,k)\nn
&&~~~~~~~~~~~~~~~~~~~~~~~~~~~~~\left.
\begin{array}{l}\\
\left[\begin{array}{c|c} \sum_\ell S_\ell^*[\Upsilon_\ell]&{1\over 2}BE\cr\hline {1\over 2}E^TB^T&E^T\left[\sum_k\cR_k^*[\Lambda_K]\right]E
\end{array}\right]\succeq0\\
\end{array}\right\}
\ee{wehave16}
(we have used the straightforward identity $\bar{\cR}^*_k[\Lambda_k]=E^T\cR^*_k[\Lambda_k]E$). By the same Proposition \ref{propmaxqf}, the optimization problem in (\ref{wehave16}) specifying $\Opt$ is solvable, and
\begin{equation}\label{estimatetiught}
\varrho\leq \Opt\leq {2} \ln({2}D) {\varrho},\;\;D=\sum_kd_k+\sum_\ell f_\ell.
\end{equation}
\paragraph{5$^o$.} Let $\bar{\Upsilon}=\{\bar{\Upsilon}_\ell\}$, $\bar{\Lambda}=\{\bar{\Lambda}_k\}$ be an optimal solution to the optimization problem specifying $\Opt$, see
(\ref{wehave16}), and let
\[
{\mathbf{\Upsilon}}=\sum_\ell \cS^*_\ell[\bar{\Upsilon}_\ell],\;{\mathbf{\Lambda}}=\sum_k \cR^*_k[\bar{\Lambda}_k],
\]
so that
\begin{equation}\label{sothat23}
\Opt=\phi_{\cR}(\lambda[\bar{\Upsilon}])+\phi_{\cT}(\lambda[\bar{\Lambda}])\ \&\ \left[\begin{array}{c|c}{\mathbf{\Upsilon}}&{1\over 2}BE\cr\hline {1\over 2}E^TB^T&E^T{\mathbf{\Lambda}}E\cr\end{array}\right]\succeq0.
\end{equation}
We claim that for properly selected $m\times\nu$ matrix $H$ it holds
\begin{equation}\label{target}
\left[\begin{array}{c|c}{\mathbf{\Upsilon}}&{1\over 2}(B-H^TA)\cr\hline{1\over 2}(B-H^TA)^T&{\mathbf{\Lambda}}\cr\end{array}\right]\succeq0.
\end{equation}
This claim implies the conclusion of Proposition \ref{propbnoiseopt}: by the claim, we have $\overline\Opt\leq \Opt$, which combines with (\ref{estimatetiught}) and (\ref{ensurethat}) to imply   (\ref{newupper}).\par
In order to justify the claim, assume that it fails to be true, and
 let us lead this assumption to a contradiction.
To this end, consider the semidefinite program
\begin{equation}\label{sdprogram}
\tau_*=\min_{\tau,H}\left\{\tau:\left[\begin{array}{c|c}{\mathbf{\Upsilon}}&{1\over 2}(B-H^TA)\cr\hline{1\over 2}(B-H^TA)^T&{\mathbf{\Lambda}}\cr\end{array}\right]+\tau I_{\nu+n}\succeq0\right\}.
\end{equation}
The problem clearly is strictly feasible, and the value of the objective at every feasible solution is positive. In addition, the problem is solvable (by exactly the same argument as in item 1$^o$ of the proof).
\paragraph{4$^0$.b.} As we have seen, (\ref{sdprogram}) is a strictly feasible solvable problem with positive optimal value $\tau_*$, so that the problem dual to (\ref{sdprogram}) is solvable with positive optimal value. Let us build the dual problem. Denoting by $\left[\begin{array}{c|c}U&V\cr\hline V^T&W\cr\end{array}\right]\succeq0$ the Lagrange multipliers for the semidefinite constraint in (\ref{sdprogram}) and taking inner product of the left hand side of the  constraint with the multiplier, we get the aggregated constraint
$$
\begin{array}{l}
\Tr(U{\mathbf{\Upsilon}})+\Tr(W{\mathbf{\Lambda}})+\tau[\Tr(U)+\Tr(W)]+\Tr((B-H^TA)V^T)\geq0.
\end{array}
$$
The equality constraints of the dual problem should make the homogeneous in $\tau,H$ part of the left hand side in the aggregated constraint identically equal to $\tau$, which amounts to
\begin{equation}\label{amountsto}
\Tr(U)+\Tr(W)=1,\;VA^T=0,
\end{equation}
so the aggregated constraint reads
$$
\tau\geq-\left[\Tr(U{\mathbf{\Upsilon}})+\Tr(W{\mathbf{\Lambda}})+\Tr(BV^T)\right].
$$
The dual problem is to maximize the right hand side of the latter constraint over Lagrange multiplier $\left[\begin{array}{c|c}U&V\cr\hline V^T&W\cr\end{array}\right]\succeq0$ satisfying (\ref{amountsto}), and its optimal value is $\tau_*>0$, that is, there exists $\left[\begin{array}{c|c}\bar{U}&\bar{V}\cr\hline \bar{V}^T&\bar{W}\cr\end{array}\right]\succeq0$ such that $A\bar{V}^T=0$ and
\begin{equation}\label{toifel}
\Tr(\bar{U}{\mathbf{\Upsilon}})+\Tr(\bar{W}{\mathbf{\Lambda}})+\Tr(B\bar{V}^T)<0.
\end{equation}
Adding to $\bar{U}$ a small positive multiple of the unit matrix, we can assume, in addition, that $\bar{U}\succ0$. Now, the relation $A\bar{V}^T=0$ combines with the definition of $E$ to imply that $\bar{V}^T=EF$ for properly selected matrix $F$, so that
$$
\left[\begin{array}{c|c}\bar{U}&F^TE^T\cr\hline EF&\bar{W}\cr\end{array}\right]\succeq0.
$$
Hence, by Schur Complement Lemma,
$$
\bar{W}\succeq EF\bar{U}^{-1}F^TE^T,
$$
and (\ref{toifel}) combines with ${\mathbf{\Lambda}}\succeq0$ to imply that
$$
\begin{array}{l}
0>\Tr(\bar{U}{\mathbf{\Upsilon}})+\Tr(\bar{W}{\mathbf{\Lambda}})+\Tr(B\bar{V}^T)=\Tr(\bar{U}{\mathbf{\Upsilon}})+\Tr(\bar{W}{\mathbf{\Lambda}})+\Tr(BEF)\\
\geq
\Tr(\bar{U}{\mathbf{\Upsilon}})+\Tr(EF\bar{U}^{-1}F^TE^T{\mathbf{\Lambda}})+\Tr(BEF)=\Tr\left(\left[\begin{array}{c|c}{\mathbf{\Upsilon}}&{1\over 2}BE\cr\hline
{1\over 2}E^TB^T&E^T{\mathbf{\Lambda}}E\cr\end{array}\right]\left[\begin{array}{c|c}\bar{U}&F^T\cr\hline F&F\bar{U}^{-1}F^T\cr\end{array}\right]\right)\\
\end{array}
$$
Both matrix factors in the concluding the chain $\Tr(\cdot)$ are positive semidefinite (the first one due to (\ref{sothat23}), and the second -- by Schur Complement Lemma); consequently, the concluding
quantity in the chain is nonnegative, which is impossible. We have arrived at a desired contradiction.
 \qed

\appendix

\section{Calculus of spectratopes} \label{app:Scalc}
The principal rules of the calculus of spectratopes are as follows:
\begin{itemize}
\item {[finite intersections]} If
\[\X_\ell=\{x\in\bR^\nu:\;\exists (y^\ell\in\bR^{n_\ell},t^\ell\in\T_\ell):\;x=P_\ell y^\ell,\,R_{k\ell}^2[y^\ell]\preceq t^\ell_kI_{d_{k\ell}},\,k\leq K_\ell\}, \;\;1\leq\ell\leq L,
\] are spectratopes, so is $\X=\bigcap\limits_{\ell\leq L} \X_\ell$. Indeed, let
\[E=\{[y=[y^1;...;y^L]\in\bR^{n_1}\times...\times\bR^{N_L}: \;P_1y^1=P_2y^2=...=P_Ly^L\}.
\] When $E=\{0\}$, we have $\X=\{0\}$, so that $\X$ is a spectratope; when $E\neq\{0\}$, we have
    $$
    \begin{array}{l}
    \X=\{x\in\bR^\nu:\exists (y=[y^1;...;y^L]\in E,t=[t^1;...;t^L]\in\T:=\T_1\times...\times\T_L):\\
    \multicolumn{1}{r}{ x=Py:=P_1y^1,R_{k\ell}^2[y^\ell]\preceq t^\ell_kI_{d_{k\ell}},1\leq k\leq K_\ell,1\leq \ell\leq L\};}\\
    \end{array}
    $$
    identifying $E$ and appropriate $\bR^n$, we arrive at a valid spectratopic representation of $\X$.
\item {[direct product]} If
\[\X_\ell=\{x^\ell\in\bR^{\nu_\ell}:\exists (y^\ell\in\bR^{n_\ell},t^\ell\in\T_\ell): x^\ell=P_\ell y^\ell,R_{k\ell}^2[y^\ell]\preceq t^\ell_kI_{d_{k\ell}},k\leq K_\ell\}, \;\;1\leq\ell\leq L,
\] are spectratopes, so is $\X=\X_1\times...\times \X_L$:
$$
\begin{array}{l}
\X_1\times...\times \X_L=\{x=[x^1;...;x^L]:\exists (y=[y^1;...;y^L],t=[t^1;...;t_L]\in\T=\T_1\times...\times \T_\ell):\\
\multicolumn{1}{r}{x=Py:=[P_1y^1;...;P_Ly^L],R_{k\ell}^2[y^\ell]\preceq t^\ell_kI_{d_{k\ell}},1\leq k\leq K_\ell,1\leq \ell\leq L\};}\\
\end{array}
$$
\item {[linear image]} If
\[
\X=\{x\in\bR^{\nu}:\exists (y\in\bR^{n},t\in\T): x=Py,R_{k}^2[y]\preceq t_kI_{d_{k}},k\leq K\}
\] is a spectratope and $S$ is a $\mu\times\nu$ matrix, the set $S\X=\{z=Sx:x\in \X\}$ is a spectratope:
$$
S\X=\{z\in\bR^\mu:\exists (y\in\bR^n,t\in\T): z=SPy,R_k^2[y]\preceq t_kI_{d_k},k\leq K\}.
$$
\item {[inverse linear image under embedding]}  If
\[\X=\{x\in\bR^{\nu}:\exists (y\in\bR^{m},t\in\T): x=Py,R_{k}^2[y]\preceq t_kI_{d_{k}},k\leq K\}\]
is a spectratope, and $S$ is a $\nu\times\mu$ matrix with trivial kernel, the set $S^{-1}\X=\{z:Sz\in \X\}$ is a spectratope. Indeed, setting $E=\{y\in\bR^m:Py\in\hbox{Im} S\}$, we get a linear subspace of $\bR^n$; if $E=\{0\}$, $S^{-1}\X=\{0\}$ is a spectratope, otherwise we have
$$
S^{-1}\X=\{z\in\bR^\mu:\exists (y\in E,t\in\T): z=Qy,R_k^2[y]\preceq t_kI_{d_k},k\leq K\},
$$
where linear mapping $y\mapsto Qy:E\to\bR^\mu$ is uniquely defined by the relation $Py=SQy$. When identifying $E$ with appropriate $\bR^n$, we get a valid spectratopic representation of $S^{-1}\X$.
\item {[arithmetic sum]} If $\X_\ell$, $\ell\leq L$, are spectratopes in $\bR^\nu$, so is the arithmetic sum $\X=\X_1+...+\X_L$ of $\X_\ell$. Indeed, $\X$ is the image of $\X_1\times...\times \X_L$ under the linear mapping $[x^1;...;x^L]\mapsto x^1+...+x^L$, and taking direct products and linear images preserve spectratopes.
\end{itemize}

\section{Processing covariance estimation problem in the diagonal case}\label{Processing}
We start with  setting some additional notation to be used when operating with Euclidean space $\bS^n$.
 \begin{itemize}
 \item We denote $\bar{n}={n(n+1)\over2}=\dim\bS^n$, $\cI=\{(i,j): 1\leq i\leq j\leq n\}$, and for $(i,j)\in\cI$ we set $$e^{ij}=\left\{\begin{array}{ll}e_ie_i^T,&i=j,\\
 {1\over\sqrt{2}}[e_ie_j^T+e_je_i^T],&i<j,\cr\end{array}\right.$$ where $e_i$ are the standard basic orths in $\bR^n$. Note that $\{e^{ij}:(i,j)\in\cI\}$ is the standard orthonormal basis in $\bS^n$. Given $v\in\bS^n$, we denote by
 $\rX(v)$ the vector of coordinates of $v$ in this basis: $$\rX_{ij}(v)=\Tr(ve^{ij})=\left\{\begin{array}{ll}v_{ii},&i=j,\\
 \sqrt{2}v_{ij},&i<j,\cr\end{array}\right.\;\; (i,j)\in\cI.$$
 Similarly, for $x\in\bR^{\bar{n}}$, we index the entries in $x$ by pairs $ij$, $(i,j)\in\cI$, and set $\rV(x)=\sum_{(i,j)\in \cI_p} x_{ij}e^{ij}$, so that $v\mapsto \rX(v)$ and $x\mapsto\rV(x)$ are inverse to each other linear isometries identifying the Euclidean spaces $\bS^n$ and $\bR^{\bar{n}}$ (recall that the inner products on these spaces are, respectively, the Frobenius and the standard one).
 \item Recall that $\cV$ is the matrix box $\{v\in\bS^n: v^2\preceq I_n\}=\{v\in\bS^n: \exists t\in\cT:=[0,1]: v^2\preceq tI_n\}$. We denote by $\cX$ the image of $\cV$ under the mapping $\rX$:
 $$
 \cX=\{x\in\bR^{\bar{n}}: \exists t\in[0,1]: \rV^2[x]\preceq tI_n\}.
 $$
 Note that $\cX$ is a {basic} spectratope.
 \end{itemize}
 Now we can assume that the signal underlying our observations is $x\in\cX$, and the observations themselves are
 $$
 w_t=\rX(\omega_t) = \underbrace{\rX({\half}\rV(x))}_{=\half x}+z_t,\,\,\,z_t=\rX(\zeta_t).
 $$
 Note that $z_t\in\bR^{\bar{n}}$, $1\leq t\leq T$, are zero mean i.i.d. random vectors with covariance matrix $D[x]$ satisfying, in view of (\ref{conclude}), the relation
 $$
 x\in\cX\;\Rightarrow\;D[x]\preceq 2I_{\bar{n}}
 $$
 (recall that we are in the case of $A=I_n$, $\vartheta_*=I_n$).
 Our goal is to estimate $B\vartheta[v]B^T={\half}BB^T+{\half}BvB^T$, or, what is the same, to recover
 $$
 \overline{B}x:={1\over 2}\rX(B\rV(x)B^T).
 $$
 Recall that we are in the situation where the norm in which the recovery error is measured is either the Frobenius, or the nuclear norm on $\bS^n$; we ``transfer'' this norm from $\bS^n$ to $\bR^{\bar{n}}$. In the situation in question, Proposition \ref{proprepeated} supplies the linear estimate
 $$
 \widehat{x}(w^{(T)})={1\over T}H^T_*\sum_{t=1}^Tw_t
 $$
of $\overline{B}x$ with ${H}_*$ stemming from the optimal solution to the convex optimization problem presented in Proposition \ref{proprepeated}. It is immediately seen that under the circumstances, this optimization problem reads
\\ in the Frobenius norm case:
\begin{equation}\label{weconsideragainFN-F}
\overline{\Opt}=\hbox{\footnotesize$\min\limits_{\begin{array}{c}\bar{H},\Lambda\\\upsilon,\upsilon',\Theta\end{array}}\left\{\Tr(\Lambda)+\upsilon+\upsilon'+{2\over T}\Tr(\Theta):
\left\{\begin{array}{r}\bar{H}\in\bR^{\bar{n}\times\bar{n}},\Lambda\in\bS^n_+,\upsilon\in\bR_+,\upsilon'\in\bR_+,\Theta\in\bS^{\bar{n}}\\
\left[\begin{array}{c|c}\R^*[\Lambda]&\half [\overline{B}^T-\half\bar{H}]\cr\hline\half [\overline{B}-\half \bar{H}^T]&\upsilon I_{\bar{n}}\cr\end{array}\right]\succeq0,\\
\left[\begin{array}{c|c}\Theta&\half \bar{H}\cr \hline \half \bar{H}^T&\upsilon^\prime I_{\bar{n}}\cr\end{array}\right]\succeq0;
\end{array}\right.\right\}$}
\end{equation}
in the nuclear norm case:
\begin{equation}\label{weconsideragainFN-N}
\overline{\Opt}=\hbox{\footnotesize$\min\limits_{\begin{array}{c}\bar{H},\Lambda\\\Upsilon,\Upsilon',\Theta\end{array}}\left\{\Tr(\Lambda)+\Tr(\Upsilon)+\Tr(\Upsilon')+{2\over T}\Tr(\Theta):
\left\{\begin{array}{r}\bar{H}\in\bR^{\bar{n}\times\bar{n}},\Lambda\in\bS^n_+,\Upsilon\in\bS^n_+,\Upsilon'\in\bS^n_+,\Theta\in\bS^{\bar{n}}\\
\left[\begin{array}{c|c}\R^*[\Lambda]&\half [\overline{B}^T-\half\bar{H}]\cr\hline\half [\overline{B}-\half\bar{H}^T]&\R^*[\Upsilon]\cr\end{array}\right]\succeq0,\\
\left[\begin{array}{c|c}\Theta&\half \bar{H}\cr \hline \half \bar{H}^T&\R^*[\Upsilon^\prime]\cr\end{array}\right]\succeq0
\end{array}\right.\right\}$}
\end{equation}
where
\[\R^*[S]=\left[\Tr(e^{ij}Se^{k\ell})\right]_{{(i,j)\in\cI\atop(k,\ell)\in\cI}}\in\bS^{\bar{n}},\,S\in\bS^n.
\] So far, we did not use the fact that $B$ is diagonal; this is what we intend to utilize now. Specifically, it is immediately seen that with diagonal $B$,
\begin{enumerate}
\item[1)] The $\bar{n}\times\bar{n}$ matrix $\overline{B}$ is diagonal, with diagonal entries $\overline{B}_{ij,ij}={1\over 2}B_{ii}B_{jj}$.
\item[2)] Let $\cE$ be the multiplicative group comprised of $n\times n$ diagonal matrices with diagonal entries $\pm 1$. Every matrix $E\in\cE$ induces diagonal $\bar{n}\times\bar{n}$ matrix $F_E$ with diagonal entries $\pm1$ such that
    $$
    \R^*[ESE]=F_E\R^*[S]F_E\,\,\forall (E\in\cE,S\in\bS^n).
    $$
    Indeed, when $E\in\cE$ and $S\in\bS^n$, we have
    $$
    \begin{array}{l}
    [\cR^*[ESE]]_{ij,k\ell}=\Tr(e^{ij}[ESE]e^{k\ell})=\Tr(S[Ee^{k\ell}E][Ee^{ij}E])=\Tr(S[E_{ii}E_{jj}E_{kk}E_{\ell\ell}e^{k\ell}e^{ij}])\\
    =
    [F_E]_{ij,ij}[F_E]_{k\ell,k\ell}\Tr(e^{ij}Se^{k\ell})= [F_E\cR^*[S]F_E]_{ij,k\ell},\\
    \end{array}
    $$
    where $F_E$ is the diagonal $\bar{n}\times\bar{n}$ matrix with diagonal entries $[F_E]_{pq,pq}=E_{pp}E_{qq}$, $(p,q)\in\cI$.
    \item[3)] When $S\in\bS^n$ is diagonal, $\cR^*[S]$ is diagonal as well.\\
     Indeed, when $S$ is diagonal and $(i,j)\in\cI$, $(k,\ell)\in\cI$ with $(i,j)\neq(k,\ell)$, the supports of matrices $e^{ij}$ and $(Se^{k\ell})^T$  (the sets of cells where the entries of the respective matrices are nonzero) do not intersect, whence $[\cR^*[S]]_{ij,k\ell}=\Tr(e^{ij}Se^{k\ell})=\sum_{p,q}[e^{ij}]_{pq}[(Se^{k\ell})^T]_{pq}=0$.
\end{enumerate}
Observe that by 1) -- 3) problems \rf{weconsideragainFN-F} and \rf{weconsideragainFN-N} have optimal solutions with diagonal matrix components. To see this, let us start with the nuclear norm case. Let $\bar{H}_*,\Lambda_*,\Upsilon_*,\Upsilon^\prime_*,\Theta_*$ be an optimal solution to the problem. We claim that when $E\in \cE$, the collection $F_E\bar{H}_*F_E,E\Lambda_*E,E\Upsilon_*E,E\Upsilon^\prime_*E,F_E\Theta_*F_E$ is an optimal solution as well. Indeed, the values of the objective at the original and the transformed solution clearly are the same, so that all we need in order to justify our claim is to check that the transformed solution is feasible, which boils down to verifying that it satisfies the LMI constraints of the problem. We have
{\small $$
\begin{array}{l}
\left[\begin{array}{c|c}\R^*[E\Lambda_* E]&\half [\overline{B}^T-\half F_E\bar{H}_*F_E]\cr\hline\half [\overline{B}-\half [F_E\bar{H}_*F_E]^T]&\R^*[E\Upsilon_* E]\cr\end{array}\right]=\left[\begin{array}{c|c}F_E\R^*[\Lambda]F_E&\half F_E[\overline{B}^T-\half \bar{H}_*]F_E\cr\hline\half F_E[\overline{B}-\half \bar{H}_*^T]F_E&F_E\R^*[\Upsilon_*]F_E\cr\end{array}\right]\\
=\Diag\{F_E,F_E\}\left[\begin{array}{c|c}\R^*[\Lambda]&\half [\overline{B}^T-\half \bar{H}_*]\cr\hline\half [\overline{B}-\half \bar{H}_*^T]&\R^*[\Upsilon_*]\cr\end{array}\right]\Diag\{F_E,F_E\}\succeq0,
\end{array}
$$}\noindent
where the first equality is due to 3) combined with diagonality of $\overline{B}$ stated in 1). Similarly,
$$
\begin{array}{l}
\left[\begin{array}{c|c}F_E\Theta_*F_E&\half F_E\bar{H}_*F_E\cr \hline \half F_E\bar{H}_*^TF_E&\R^*[E\Upsilon^\prime_* E]\cr\end{array}\right]=
\left[\begin{array}{c|c}F_E\Theta_*F_E&\half F_E\bar{H}_*F_E\cr \hline \half F_E\bar{H}_*^TF_E&F_E\R^*[\Upsilon^\prime_*]F_E\cr\end{array}\right] \\
=\Diag\{F_E,F_E\}\left[\begin{array}{c|c}\Theta_*&\half \bar{H}_*\cr \hline \half \bar{H}_*^T&\R^*[\Upsilon^\prime_*]\cr\end{array}\right]\Diag\{F_E,F_E\}\succeq0.
\end{array}
$$
Thus, the transformed solution indeed is feasible, as claimed. Since the problem is convex, the average, over $E\in\cE$, of the above transformations of an optimal solution again is an optimal solution, let it be denoted $\bar{H}_{\#},\Lambda_{\#},\Upsilon_{\#},\Upsilon^\prime_{\#},\Theta_{\#}$. By construction, $\Lambda_{\#},\Upsilon_{\#},\Upsilon^\prime_{\#}$ are diagonal, whence by 3) $\cR^*[\Lambda_{\#}],\cR^*[\Upsilon_{\#}],\cR^*[\Upsilon^\prime_{\#}]$ are diagonal as well. This combines with diagonality of $\overline{B}$ to imply, similarly to the above, that if $L$ is a diagonal $\bar{n}\times\bar{n}$ matrix with diagonal entries $\pm1$, the collection $L\bar{H}_{\#}L,\Lambda_{\#},\Upsilon_{\#},\Upsilon^\prime_{\#},L\Theta_{\#}L$ is an optimal solution to \rf{weconsideragainFN-N}. Averaging these optimal solutions over $L$'s, we conclude that the problem has an optimal solution comprised of diagonal matrices, as claimed. The same reasoning, with evident simplifications, works in the case of Frobenius norm.
\par
We see that when solving \rf{weconsideragainFN-F} and \rf{weconsideragainFN-N}, we lose nothing when restricting ourselves to candidate solutions with diagonal matrix components, which, by 1) and  3) automatically ensures the diagonality of blocks in the LMI constraints of the problem. As a result, we, first, reduce dramatically the design dimension of the problem, and, second, can now replace ``large-scale'' LMI constraints (which now state that some $2\times 2$ block matrices with diagonal
$\bar{n}\times\bar{n}$ blocks should be $\succeq0$) with a bunch of small -- just $2\times 2$ -- LMI's, thus making the problems easily solvable by the existing software, e.g. {\tt CVX} \cite{cvx2014}, provided $n$ is in the range of hundreds.

\section{Conic duality}\label{conicd}
A conic problem is an optimization problem of the form
$$
\Opt(P)=\max_x\left\{c^Tx: A_ix-b_i\in\bK_i,i=1,...,m, Px=p\right\}\eqno{(P)}
$$
where $\bK_i$ are regular (i.e., closed, convex, pointed and with a nonempty interior) cones in Euclidean spaces $E_i$. Conic dual of $(P)$ is ``responsible'' for upper-bounding the optimal value in $(P)$ and is built as follows: selecting somehow  {\sl Lagrange multipliers} $\lambda_i$ for the conic constraints $A_ix-b_i\in\bK_i$ in the cones dual to $\bK_i$:
$$
\lambda_i\in \bK_i^*:=\{\lambda:\langle\lambda,y\rangle\geq0\,\forall y\in\bK_i\},
$$
and a Lagrange multiplier $\mu\in\bR^{\dim p}$ for the equality constraints,
every feasible solution $x$ to $(P)$ satisfies the linear inequalities $\langle\lambda_i,A_ix\rangle\geq \langle \lambda_i,b_i\rangle$, $i\leq m$, same as the inequality $\mu^TPx\geq \mu^Tp$, and thus satisfies the aggregated inequality
$$
\sum_i\langle\lambda_i,A_ix\rangle +\mu^TPx\geq \sum_i\langle \lambda_i,b_i\rangle +\mu^Tp.
$$
If the left hand side of this inequality is, {\sl identically in $x$}, equal to $-c^Tx$ (or, which is the same, $-c=\sum_iA_i^*\lambda_i+P^T\mu$, where $A_i^*$ is the conjugate of $A_i$), the inequality produces an upper bound $-\langle \lambda_i,b_i\rangle-p^T\mu$ on $\Opt(P)$. The dual problem
$$
\Opt(D)=\min_{\lambda_1,...,\lambda_m,\mu}\left\{-\sum_i\langle \lambda_i,b_i\rangle-p^T\mu:\lambda_i\in\bK_i^*,i\leq m, \sum_iA_i^*\lambda_i+P^T\mu=-c\right\}\eqno{(D)}
$$
is the problem of minimizing this upper bound. Note that $(D)$ is a conic problem along with $(P)$ -- it is a problem of optimizing a linear objective under a bunch of linear equality constraints and conic inclusions of the form ``affine function of the decision vector should belong to a given regular cone.'' Conic problem, like $(P)$, is called {\sl strictly feasible}, if it admits a feasible solution $x$ for which all conic inclusions are satisfied strictly: $A_ix-b_i\in\inter K_i$ for all $i$. Conic Duality Theorem (see, e.g., \cite{LMCO}) states that when one of the problems $(P)$, $(D)$ is bounded\footnote{for a maximization (minimization)  problem, boundedness means that the objective is bounded from above (resp., from below) on the feasible set.} and strictly feasible, then the other problem in the pair is solvable, and $\Opt(P)=\Opt(D)$.

\end{document}